\documentclass[UTF-8,reqno]{amsart}
\usepackage{enumerate, bbm}
\usepackage{amssymb,url,color, xcolor, booktabs}
\usepackage{tikz}
\usepackage{mathrsfs}
\usepackage{dutchcal}
\usepackage{threeparttable}
\usepackage{makecell}
\usepackage{color}
\usepackage[colorlinks=true]{hyperref}
\hypersetup{
    linkcolor=blue,          
    citecolor=red,        
    filecolor=blue,      
    urlcolor=cyan
}

\usepackage{color}
\usepackage{ulem}

\definecolor{MyDarkBlue}{cmyk}{0.8,0.3,0.8,0.4}
\definecolor{yellow}{rgb}{0.99,0.99,0.70}
\definecolor{white}{rgb}{1.0,1.0,1.0}
\definecolor{black}{rgb}{0.00,0.00,0.00}

\newcommand{\red}{\color{red}}

\numberwithin{equation}{section}

\newcommand{\be}{\begin{eqnarray}}
\newcommand{\ee}{\end{eqnarray}}
\newcommand{\ce}{\begin{eqnarray*}}
\newcommand{\de}{\end{eqnarray*}}
\newtheorem{theorem}{Theorem}[section]
\newtheorem{lemma}[theorem]{Lemma}
\newtheorem{remark}[theorem]{Remark}
\newtheorem{definition}[theorem]{Definition}
\newtheorem{proposition}[theorem]{Proposition}
\newtheorem{Examples}[theorem]{Example}
\newtheorem{corollary}[theorem]{Corollary}

\def\bx{\begin{Examples}}
\def\ex{\end{Examples}}

\def\1{{\mathbf{1}}}

\def\eps{\varepsilon}

\def\e{\mathrm{e}}

\def\p{\partial}

\def\[{{\Big[}}
\def\]{{\Big]}}
\def\<{{\langle}}
\def\>{{\rangle}}
\def\({{\big(}}
\def\){{\big)}}

\def\dif{{\mathord{{\rm d}}}}

\def\min{{\mathord{{\rm min}}}}

\def\bb2{{\boldsymbol{2}}}
\def\no{\nonumber}
\def\={&\!\!=\!\!&}

\def\bB{{\mathbf B}}

\def\bC{{\mathbf C}}
\def\bI{{\mathbf I}}

\def\cB{{\mathcal B}}

\def\cH{{\mathcal H}}

\def\cP{{\mathcal P}}

\def\cR{{\mathcal R}}
\def\cS{{\mathcal S}}

\def\mB{{\mathbb B}}
\def\mC{{\mathbb C}}

\def\mE{{\mathbb E}}

\def\mI{{\mathbb I}}

\def\mH{{\mathbb H}}
\def\mL{{\mathbb L}}

\def\mN{{\mathbb N}}

\def\mP{{\mathbb P}}

\def\mR{{\mathbb R}}

\def\mT{{\mathbb T}}

\def\mZ{{\mathbb Z}}

\def\bP{{\mathbf P}}

\def\bE{{\mathbf E}}

\def\b1{{\mathbbm 1}}

\def\sA{{\mathscr A}}
\def\sB{{\mathscr B}}

\def\sF{{\mathscr F}}

\def\sI{{\mathscr I}}

\def\sL{{\mathscr L}}

\def\sN{{\mathscr N}}

\def\sS{{\mathscr S}}
\def\sT{{\mathscr T}}
\def\sU{{\mathscr U}}
\def\sV{{\mathscr V}}

\def\geq{\geqslant}
\def\leq{\leqslant}
\def\ge{\geqslant}
\def\le{\leqslant}

\def\div{\mathord{{\rm div}}}

\def\eps{\varepsilon}

\def\e{\mathrm{e}}

\def\p{\partial}

\def\[{{\Big[}}
\def\]{{\Big]}}
\def\<{{\langle}}
\def\>{{\rangle}}

\def\dif{{\mathord{{\rm d}}}}

\def\min{{\mathord{{\rm min}}}}

\def\no{\nonumber}
\def\={&\!\!=\!\!&}
\def\bt{\begin{theorem}}
\def\et{\end{theorem}}
\def\bl{\begin{lemma}}
\def\el{\end{lemma}}
\def\br{\begin{remark}}
\def\er{\end{remark}}
\def\bd{\begin{definition}}
\def\ed{\end{definition}}
\def\bp{\begin{proposition}}
\def\ep{\end{proposition}}
\def\bc{\begin{corollary}}
\def\ec{\end{corollary}}

\def\geq{\geqslant}
\def\leq{\leqslant}
\def\ge{\geqslant}
\def\le{\leqslant}

\def\div{\mathord{{\rm div}}}

\def\bH{{\mathbf H}}

\def\bX{{\mathbf X}}

 \def\R{\mathbb R}
 \def\R{\mathbb R}

\def\<{\langle} \def\>{\rangle}

\def\wt{\widetilde}

\def\ww{{w}}

\def\bb2{{\boldsymbol{2}}}
\def\bbb1{\boldsymbol{1}}

\def\no{\nonumber}
\def\={&\!\!=\!\!&}

\allowdisplaybreaks

\begin{document}
	
\title{SDEs with supercritical  distributional drifts}
\author{Zimo Hao and Xicheng Zhang}



\address{Zimo Hao:
Fakult\"at f\"ur Mathematik, Universit\"at Bielefeld,
33615, Bielefeld, Germany\\
Email: zhao@math.uni-bielefeld.de}

\address{Xicheng Zhang:
School of Mathematics and Statistics, Beijing Institute of Technology, Beijing 100081, China\\
Faculty of Computational Mathematics and Cybernetics, Shenzhen MSU-BIT University, 518172 Shenzhen, China\\
Email: XichengZhang@gmail.com
 }

\thanks{
This work is partially supported by National Key R\&D program of China (No. 2023YFA1010103) and  NNSFC grants of China (Nos. 12131019), and the German Research Foundation (DFG) through the 
Collaborative Research Centre(CRC) 1283/2 2021 - 317210226 ``Taming uncertainty and profiting from randomness and low regularity in analysis, stochastics and their applications".}

\begin{abstract}
Let $d\geq 2$.  In this paper, we investigate the following  stochastic differential equation (SDE) in ${\mathbb R}^d$ driven by Brownian motion 
$$
{\rm d} X_t=b(t,X_t){\rm d} t+\sqrt{2}{\rm d} W_t,
$$
where $b$ belongs to the space ${\mathbb L}_T^q \mathbf{H}_p^\alpha$ with $\alpha \in [-1, 0]$ and $p,q\in[2, \infty]$, which is a distribution-valued and divergence-free vector field.
In the subcritical case $\frac dp+\frac 2q<1+\alpha$, we establish the existence and uniqueness of a weak solution to the integral equation:
$$
X_t=X_0+\lim_{n\to\infty}\int^t_0b_n(s,X_s){\rm d} s+\sqrt{2} W_t.
$$
Here, $b_n:=b*\phi_n$ represents the mollifying approximation, and the limit is taken in the $L^2$-sense.
In the critical and supercritical case $1+\alpha\leq\frac dp+\frac 2q<2+\alpha$, assuming the initial distribution has an $L^2$-density, we show the existence of weak solutions and associated Markov processes.
Moreover, under the additional assumption that $b=b_1+b_2+\div a$, where $b_1\in {\mathbb L}^\infty_T{\mathbf B}^{-1}_{\infty,2}$, $b_2\in {\mathbb L}^2_TL^2$, 
and $a$ is a bounded antisymmetric matrix-valued function, we establish the convergence of mollifying approximation solutions without the need to subtract a subsequence.
To illustrate our results, we provide examples of Gaussian random fields and singular interacting particle systems, including the two-dimensional vortex models.
\end{abstract}

\maketitle
\setcounter{tocdepth}{2}
\tableofcontents

\section{Introduction}
Throughout this paper, we fix $d \geq 2$ and consider the following stochastic differential equation (SDE) in $\mathbb{R}^d$, driven by a standard $d$-dimensional Brownian motion $W$:  
\begin{align}\label{in:SDE}
\mathrm{d} X_t = b(t, X_t) \, \mathrm{d}t + \sqrt{2} \, \mathrm{d}W_t, \quad t \geq 0,
\end{align}  
where the drift $b$ is a time-dependent distribution.  
Since $b$ is a distribution, the drift term is not meaningful in the classical sense, as it is not possible to assign a value to a distribution at the point $X_t$. To define solutions, a natural approach is to use mollifying approximations. Let $\phi_n(x) = n^d \phi(nx)$ be a family of mollifiers, where $\phi \in C_c^\infty(\mathbb{R}^d)$ is a smooth probability density function with compact support. The smooth approximation of $b$ is then defined by convolution as follows:  
\begin{align}\label{BN1}
b_n(t, x) := b(t, \cdot) * \phi_n(x).
\end{align}

Let $\mathcal{P}(\mathbb{R}^d)$ denote the space of all probability measures on $\mathbb{R}^d$.  
We begin by introducing the concept of weak solutions for the SDE \eqref{in:SDE} with a distributional drift.

\begin{definition}[Weak solutions]\label{Def1}  
Let $\mathfrak{F} := (\Omega, \sF, (\sF_t)_{t \geq 0}, \mathbb{P})$ be a stochastic basis, and let $(X, W)$ be a pair of $\mathbb{R}^d$-valued, continuous, $\sF_t$-adapted processes on $\mathfrak{F}$. We call $(\mathfrak{F}, X, W)$ a weak solution of the SDE \eqref{in:SDE} with initial distribution $\mu \in \mathcal{P}(\mathbb{R}^d)$ if $W$ is an $\sF_t$-Brownian motion, $\mathbb{P} \circ X_0^{-1} = \mu$, and for all $t \in [0, T]$,  
$$
X_t = X_0 + A^b_t + \sqrt{2} W_t, \quad \text{a.s.},
$$  
where $A^b_t := \lim_{n \to \infty} \int_0^t b_n(s, X_s) \, \mathrm{d}s$ exists in the $L^2$-sense.  
\end{definition}

To the best of the authors' knowledge, the concept of weak solutions introduced above first appeared in \cite{BC01} (see also \cite{ZZ17}). It is important to note that the definition of a weak solution in Definition \ref{Def1} depends on the choice of mollifiers $\phi_n$. The primary objective of this paper is to identify the weakest possible conditions on $b$ under which the SDE \eqref{in:SDE} admits a solution in the sense of Definition \ref{Def1} and to construct the corresponding Markov process.  

To present our results, we begin with a basic scaling analysis, which allows us to distinguish between subcritical, critical, and supercritical drifts.  
Let $\bH^\alpha_p$ (resp. $\dot{\bH}^\alpha_p$) denote the (resp. homogeneous) Bessel potential space, where $\alpha \in \mathbb{R}$ and $p \in [1, \infty]$ (see \cite{BL76} or Subsection 2.1 below for precise definitions). Suppose that $b \in L^q(\mathbb{R}_+; \dot{\bH}^\alpha_p)$ for some $\alpha \in \mathbb{R}$ and $q, p \in [1, \infty]$,    
and the SDE \eqref{in:SDE} admits a solution denoted by $X$. For $\lambda > 0$, define  
$$
X^\lambda_t := \lambda^{-1} X_{\lambda^2 t}, \quad W^\lambda_t := \lambda^{-1} W_{\lambda^2 t}, \quad b_\lambda(t, x) := \lambda b(\lambda^2 t, \lambda x).
$$  
Formally, it follows that  
$$
\mathrm{d} X^\lambda_t = b_\lambda(t, X^\lambda_t) \, \mathrm{d}t + \sqrt{2} \, \mathrm{d}W^\lambda_t.
$$  
Moreover, by a change of variables, we have  
$$
\|b_\lambda\|_{L^q(\mathbb{R}_+; \dot{\bH}^\alpha_p)} = \lambda^{1 + \alpha - \frac{d}{p} - \frac{2}{q}} \|b\|_{L^q(\mathbb{R}_+; \dot{\bH}^\alpha_p)}.
$$  
As $\lambda \to 0$, we classify the SDE \eqref{in:SDE} into the following three cases:  
$$
\text{\bf Subcritical: } \tfrac{d}{p} + \tfrac{2}{q} < 1 + \alpha; \quad  
\text{\bf Critical: } \tfrac{d}{p} + \tfrac{2}{q} = 1 + \alpha; \quad  
\text{\bf Supercritical: } \tfrac{d}{p} + \tfrac{2}{q} > 1 + \alpha.
$$

To solve the SDE \eqref{in:SDE}, a crucial step involves establishing a more precise regularity estimate for the associated Kolmogorov equation:  
\begin{align}\label{PDE0}
\p_t u=\Delta u+b\cdot\nabla u+f.
\end{align}
Let us perform a straightforward analysis of the differentiability index $\alpha$ when $b \in \bC^\alpha$ with $\alpha < 0$  
(see Definition \ref{Def25} below for a definition of $\bC^\alpha$).  
According to the Schauder theory for the heat equation, $u$ belongs to at most $\mathbf{C}^{2 + \alpha}$. To ensure the product $b \cdot \nabla u$ is well-defined, we need the condition $1 + 2\alpha > 0$, which implies $\alpha > -1/2$.  

We summarize some well-known results regarding the SDE \eqref{in:SDE} in three cases, categorized based on the value of $\alpha$. For clarity, we introduce the following abbreviations:  

\begin{center}
{\color{blue} \textsc{Seu}}: Strong existence-uniqueness; 
{\color{blue} \textsc{Weu}}: Weak existence-uniqueness;\\
{\color{blue} \textsc{We}}: Weak existence; 
{\color{blue} \textsc{Eue}}: Existence-uniqueness of energy solution.
\end{center}

\begin{table}[ht]
\begin{threeparttable}
\begin{tabular}{c|c|c|c}\toprule
Value of $\alpha$ & Subcritical & Critical &Supercritical \\ \midrule
$\alpha=0$  & {\color{blue} \textsc{Seu}}: V$^{\tiny 79}_{\mbox{\tiny\cite{Ve79}}}$, KR$^{\tiny 05}_{\mbox{\tiny\cite{KR05}}}$, 
Z$^{\tiny 05,10}_{\mbox{\tiny\cite{Z05, Z11}}}$  & \thead{
 {\color{blue}\textsc{Weu\&Seu}}: BFGM$^{\tiny 19}_{\mbox{\tiny\cite{BFGM}}}$, K$^{\tiny 21}_{\mbox{\tiny\cite{Kr21}}}$,\\ RZ$^{\tiny 21}_{\mbox{\tiny\cite{RZ21s, RZ23}}}$, KM$^{\tiny 23}_{\mbox{\tiny\cite{KM23}}}$ } & 
{\color{blue} \textsc{We}}: ZZ$^{\tiny 21}_{\mbox{\tiny\cite{ZZ21}}}$\\ \midrule

$\alpha\in[-\frac12,0)$ & {\color{blue} \textsc{Weu}}: BC$^{\tiny 01}_{\mbox{\tiny\cite{BC01}}}$, FIR$^{\tiny 17}_{\mbox{\tiny\cite{FIR17}}}$, ZZ$^{\tiny 17}_{\mbox{\tiny\cite{ZZ17}}}$ & -- & -- \\ \midrule

$\alpha\in[-1,-\frac12)$  & {\color{blue} \textsc{Eue}}: GP$^{\tiny 23}_{\mbox{\tiny\cite{GP23}}}$  & -- &  {\color{blue} \textsc{Eue}}: GP$^{\tiny 24}_{\mbox{\tiny\cite{GDJP}}}$,  G$^{\tiny 24}_{\mbox{\tiny\cite{Gr24}}}$   \\ \midrule
\bottomrule
\end{tabular}
\end{threeparttable}
\end{table}

\begin{itemize}
    \item[$\bullet$] 
    When $\alpha = 0$ and $p = q = \infty$, Veretenikov \cite{Ve79} first established the strong well-posedness of the SDE \eqref{in:SDE} for any starting point $X_0$.

    \item[$\bullet$] 
    In the subcritical case with $\alpha = 0$, Krylov and R\"ockner \cite{KR05} demonstrated the strong well-posedness of the SDE \eqref{in:SDE} using Girsanov's technique. Related results obtained through Zvonkin's transformation are available in \cite{Z05, Z11}.

    \item[$\bullet$] 
    In the critical case with $\alpha = 0$, the weak and strong well-posedness were studied in recent works \cite{BFGM, KM23, RZ21s, RZ23}. See also Krylov's series of works \cite{Kr21b, Kr21}.

    \item[$\bullet$] 
    In the supercritical case, when $\alpha = 0$ and $\div b = 0$, a weak solution was constructed in \cite{ZZ21} using the maximum principle proven via De-Giorgi's iteration technique. For the case of multiplicative and possibly degenerate noise, see \cite{Z23}. In the case where $\alpha=0$ and $1/q+d/p<1$, the existence of a weak solution was obtained in \cite{Kr21a} (see \cite{BG23} for the fractional Brownian motion case).

    \item[$\bullet$] 
    When $\alpha \in (-\frac{1}{2}, 0)$ and $b \in \mathbf{C}^\alpha$ is a time-independent vector field, the existence and uniqueness of a solution, termed a ``virtual solution,'' were proved in \cite{FIR17}. When $\alpha \in [-\frac{1}{2}, 0)$ and $b \in \bH^\alpha_p$ with $\alpha \in [-\frac{1}{2}, 0]$ and $p > \frac{d}{1 + \alpha}$, the authors in \cite{ZZ17} showed that there is a unique weak solution $(\mathfrak{F}, X, W)$ to the SDE \eqref{in:SDE} in the class such that the following Krylov estimate holds:
    for some $\theta, T > 0$ and any $m \in \mathbb{N}$, $f \in C_c(\mathbb{R}^d)$, and $0 \leq t_0 < t_1 \leq T$, 
    \begin{align}\label{Kry1}
        \left\|\int_{t_0}^{t_1} f(X_s) \, \dif s \right\|_{L^m(\Omega)} 
        \leq C (t_1 - t_0)^{\frac{1 + \theta}{2}} \|f\|_{\bH^\alpha_p},
    \end{align}
    where $C = C(m, T, d, \alpha, p, \|b\|_{\bH^\alpha_p}) > 0$. In particular, the above Krylov estimate implies that $A^b_t := \lim_{n \to \infty} \int_0^t b_n(X_s) \, \dif s$ exists, and $A^b_\cdot$ is a zero-energy process. In the one-dimensional case, the zero-energy solution was explored by Bass and Chen in \cite{BC01}.

    \item[$\bullet$] 
    In the subcritical case, with $\alpha \in (-1, 0]$ and a time-independent, divergence-free $b$, the existence of a unique energy solution on the torus is established in \cite{GP23} for the SDE \eqref{in:SDE}, assuming initial data whose probability distribution has an $L^2$-density with respect to the Lebesgue measure.  
     In the supercritical case, when $b = b_1 + b_2$, where $b_1$ is divergence free and belongs to $L^p(\mathbb{R}_+; \bB^{-\gamma}_{p,1})$ for some $\gamma \in [0,1]$ and $p \geq 2 / (1 - \gamma)$,  and $b_2$ represents some scaling-critical, non-divergence-free perturbations, the uniqueness of the energy solution is established in \cite{GDJP} for any initial data whose law is absolutely continuous with respect to the Lebesgue measure. See also \cite{Gr24} for related results.

    \item[$\bullet$] 
    Without assuming that $b$ is divergence-free, when $b \in \mathbf{C}^\alpha$ for some $\alpha \in \left(-\frac{2}{3}, -\frac{1}{2}\right)$, the authors of \cite{DD16} and \cite{CC18} independently utilized rough path theory and paracontrolled distribution techniques to establish the existence of a unique martingale solution in the renormalized sense. Additional references include \cite{KP20}, which considers the case of Lévy processes, and \cite{HZZZ21}, which addresses the degenerate kinetic case.  
     It is important to note, however, that in their framework, not all distributions $b \in \mathbf{C}^\alpha$ with $\alpha \in \left(-\frac{2}{3}, -\frac{1}{2}\right)$ can be treated. In fact, Theorem 6.7 in \cite{KP23} provides a counterexample showing that, for any $\alpha \leq -\frac{1}{2}$, there exists a distribution $b \in \mathbf{C}^\alpha$ for which uniqueness fails under their definition of solutions.
\end{itemize}

\subsection{Main results}
In this subsection, we present our main results in two cases: the subcritical case and the supercritical case. Our first result establishes the well-posedness of SDE \eqref{in:SDE} in the subcritical case. 

\bt[Subcritical case]\label{Th1}
Let $m \geq 2$, $\alpha_b \in (-1, -\frac{1}{2}]$, and $p_b, q_b \in [2, \infty]$ satisfy 
\begin{align}\label{MT}
\theta_0 := 1 + \alpha_b - \tfrac{d}{p_b} + \tfrac{2}{q_b} > 0, \quad m\theta_0 > 2.
\end{align}
Suppose that $b = b_1 + b_2$, where $b_1, \div b_1 \in \cap_{T>0} \mL^{q_b}_T \bH^{\alpha_b}_{p_b}$ and $|b_2(t,x)| \leq c_0 + c_1|x|$ for some $c_0, c_1 \geq 0$. 
\begin{enumerate}[(i)]
\item For any $\mu \in \cP(\mR^d)$ with finite $m$-order moment, there exists a weak solution $(\mathfrak{F}, X, W)$ to the SDE \eqref{in:SDE} 
with initial distribution $\mu$ in the sense of Definition \ref{Def1}, which is also unique 
in the class such that the following Krylov estimate holds:
for any $(\alpha, p, q) \in [\alpha_b, 0] \times [2, \infty]^2$ with 
\begin{align}\label{BQ1}
\tfrac{q_b}{2} \leq q \leq q_b, \quad p \leq p_b, \quad \alpha - \tfrac{d}{p} - \tfrac{2}{q} \geq \alpha_b - \tfrac{d}{p_b} - \tfrac{2}{q_b},
\end{align}
and for any $T > 0$, there exists $\theta > 0$ such that for any $f \in \mL^q_T \bH^{\alpha}_p \cap \mL^q_T C^\infty_b$ and $0 \leq t_0 < t_1 \leq T$, 
\begin{align}\label{Kry2}
\left\|\int^{t_1}_{t_0} f(s, X_s) \dif s \right\|_{L^m(\Omega)} \leq C (t_1 - t_0)^{\frac{1 + \theta}{2}} \|f\|_{\mL^q_T \bH^{\alpha}_p},
\end{align}
where the constant $C$ is independent of $f$ and $t_0,t_1$.

\item If $c_1 = 0$, then the moment assumption on $\mu$ can be dropped;  
and for each $t > 0$ and $X_0 = x \in \mR^d$, the law of $X_t$ admits a density $\rho_t(x, x')$ called the heat kernel of $\Delta + b \cdot \nabla$, 
which satisfies the following two-sided Gaussian estimate:  
for fixed $T > 0$, all $x, x' \in \mR^d$, and $t \in (0, T]$,  
\begin{align*}
C_0 t^{-d/2} \e^{-\gamma_0 |x - x'|^2 / t} \leq \rho_t(x, x') \leq C_1 t^{-d/2} \e^{-\gamma_1 |x - x'|^2 / t},
\end{align*}
where $C_0, C_1, \gamma_0, \gamma_1 > 0$ depend only on the parameters $T, d, \alpha_b, p_b, q_b, c_0$, and $\|b\|_{\mL^{q_b}_T \bH^{\alpha_b}_{p_b}}$.  
\end{enumerate}
\et
\br
Under the linear growth assumption of $b_2$, we require a moment assumption on the initial distribution $\mu$,  
which arises from the Krylov estimates and the Young integral 
(see Proposition \ref{Pro44} and Lemma \ref{Le42} below). Moreover, by Krylov's estimate \eqref{Kry2}, one easily sees that
$$
X_t = X_0 + A^{b_1}_t + \int_0^t b_2(s, X_s) \, \dif s + \sqrt{2} W_t, \quad \text{a.s.},
$$
where the process $A^{b_1}$ defined in Definition \ref{Def1} does not depend on the choice of the mollifiers.
\er

\br
In the subcritical case, our result improves upon \cite[Theorem 2.10]{GP23} since our initial distribution is not required to have an $L^2$-density; it can be a Dirac measure. Moreover, we are working in the whole space, not in the torus. Our proof is based on the Schauder estimate of the heat equation in Besov spaces and Zvonkin's transformation.  
\er

To present our main result in the supercritical case, we introduce the following class of distributions for later use:
\begin{align}\label{CB1}
\cB:=\left\{b\in\sS'(\mR^d;\mR^d): \|b\|_\cB:=\sup_{\varphi\in C^\infty_c(\mR^d)}\|b\cdot\nabla \varphi\|_{\bH^{-1}_{2}}/\|\varphi\|_{\bH^{1}_{2}}<\infty\right\},
\end{align}
where $\sS'(\mR^d;\mR^d)$ stands for the class of $\mR^d$-valued Schwartz distributions over $\mR^d$.

{\bf Examples:} 
(i) Suppose that $b=\div a$, where $a:\mR^d\to\mR^d\otimes\mR^d$ is an antisymmetric matrix-valued bounded measurable function. Then $b\in\cB$.
In fact, for any $h\in \bH^1_2$ and $\varphi\in C^\infty_c(\mR^d)$, 
\begin{align*}
\<b\cdot\nabla \varphi,h\>=\<\div a\cdot\nabla \varphi,h\>=-\<a,\nabla \varphi\otimes\nabla h\>\leq\|a\|_\infty\|\nabla\varphi\|_2\|\nabla h\|_2,
\end{align*}
where we have used that $\sum_{i,j}a_{ij}\p_i\p_j\varphi=0$. 

(ii) If $b, \div b\in \bB^{-1}_{\infty,2}$, then $b\in \cB$.
We shall prove it after introducing  Besov spaces (see Lemma \ref{Le27} below).

In the following we fix a time level $T>0$.
Let $\mC_T:=C([0,T];\mR^d)$ be the Banach space of all continuous functions.
 A path in $\mC_T$ is denoted by $\omega$. The canonical process is denoted by 
$w_t(\omega)=\omega(t).$
Let $\sB_t:=\sigma(w_s,s\le t)$ be  the natural filtration.

Our second main result of this paper is:

\bt[Supercritical case]\label{Th2}
Let $\alpha_b\in[-1,0]$ and $p_b,q_b\in[2,\infty]$ satisfy $\frac2{q_b}+\frac d{p_b}<2+\alpha_b$.
Suppose that $b\in  \mL^{q_b}_T\bH^{\alpha_b}_{p_b}$ is divergence-free. 
For any initial distribution $\mu\in\cP(\mR^d)$ having $L^2$-density $\rho_0$, there is a weak solution $(\mathfrak{F}, X,W)$ to SDE \eqref{in:SDE}
in the sense of Definition \ref{Def1}, which has a density $\rho_t\in \mL^\infty_TL^2\cap \mL^2_T\bH^1_2$ that satisfies the Fokker-Planck equation
$\p_t\rho=\Delta\rho-\div(b\rho)$ in the distributional sense.
Moreover, consider the approximation SDE
$$
X_t^n=X_0+\int_0^tb_n(s,X_s^n)\dif s+\sqrt2 W_t,
$$
where $b_n\in \mL^{q_b}_TC^\infty_b$ is defined by \eqref{BN1}. Let $\mP_n$ be the law of $X^n_\cdot$ in $\mC_T$.
\begin{enumerate}[{\bf (A)}]
\item  For any subsequence $n_k$, there is a subsubsequence $n_k'$ such that  for {\bf any} initial distribution $\mu\in\cP(\mR^d)$ with $L^2$-density, 
$\mP_{n'_k}$ weakly converges to a solution of SDE \eqref{in:SDE} starting from $\mu$.
\item Let $\mP_\mu$ be the law of the solution constructed  in {\bf (A)}. The following almost surely Markov property holds: there is a Lebesgue zero set $\sN\subset(0,T)$ such that for all $s\in [0,T)\backslash\sN$,
$$
\mE^{\mP_\mu} [f(w_t)|\sB_s]=\mE^{\mP_\mu} [f(w_t)|w_s],\ \forall t\in[s,T], \ f\in C_b(\mR^d).
$$
\item If in addition that $b=b_1+b_2$, where $b_1\in \mL^\infty_T\cB$ and $b_2\in \mL^2_TL^2$,
then without subtracting subsequence, $\mP_n$ weakly converges to $\mP_\mu$ in $\mC_T$ as $n\to\infty$,
where $\mP_\mu$ does not depend on the choice of the mollifiers.
\end{enumerate}
\et

\br
In the supercritical case, we are unable to demonstrate the uniqueness of weak solutions a priori. However, the assertion {\bf (A)} affirms that we can identify a subsequence such that, for any initial value, the corresponding approximation solution converges weakly to a solution. This enables us to select a Markov process.
In the case where $b=b_1+b_2$, with  $b_1\in \mL^\infty_T\cB$ and $b_2\in \mL^2_TL^2$, 
it is unique in the approximation sense since it is not necessary to choose a subsequence.
\er

The proof of Theorem \ref{Th2} depends on the solvability of PDE. Specifically, we consider the following parabolic equation:
\begin{align}\label{PDE9}
\p_t u=\Delta u+b\cdot\nabla u+f.
\end{align}
When $f=0$ and $b=\div a$, where $a:\mathbb{R}^d\to\mathbb{R}^d\otimes\mathbb{R}^d$ is an antisymmetric matrix-valued measurable function, 
the Harnack inequality was established in \cite{Os87} for $a\in L^\infty$  and in \cite{SSSZ} for $a\in BMO$.
In the case where $u_0=0$ and $b\in L^p$ with $\div b=0$, through De-Giorgi's argument, the maximum estimate of weak solution $u$ was obtained in \cite{ZZ21} for $p>d/2$
and in \cite[Corollary 1.5]{Z23} for $p>(d-1)/2$. In this work, for drifts that are distribution-valued and under certain conditions on the divergence of $b$, we show the weak well-posedness of PDE \eqref{PDE9} using the energy method.
To establish the existence of a weak solution to SDE \eqref{in:SDE}, 
we rely on two types of Krylov's estimates. One is used to prove tightness (see Lemma \ref{Le51}), while the other is utilized to take limits (see Lemma \ref{Le54}).
 
To establish conclusions {\bf (A)} and {\bf (C)}, we will employ the concept of generalized martingale solutions. 
In scenarios involving distributional drifts, to the best of our knowledge, no existing results address the equivalence between weak solutions and generalized martingale solutions for SDEs with scaling supercritical distributional drifts. Nevertheless, in our proofs, both types of solutions are obtained as subsequential limits of smooth approximation solutions. It is worth emphasizing that in the supercritical case, even in the absence of uniqueness, we retain the flexibility to select a Markov process.

As an application, we consider the following singular interaction particle system  in $\mR^{Nd}$:
\begin{align}\label{HG1}
\dif X^{N,i}_t=\sum_{j\not=i}\gamma_j K(X^{N,i}_t-X^{N,j}_t)\dif t+\sqrt{2}\dif W^{N,i}_t,\ \ i=1,\cdots,N,
\end{align}
where $K\in\bH^{-1}_\infty(\mR^d;\mR^d)$ is divergence free, $W^{N,i}_t, i=1,\cdots,N$ are $N$-independent standard $d$-dimensional Brownian motions, $\gamma_j\in\mR$ and initial value $\bX^{N}_0:=(X^{N,1}_0,\cdots, X^{N,N}_0)$ 
has an $L^2$-density $\rho^N_0$. Define $b(\mathbf{x})$ for $\mathbf{x}=(x_1,\cdots,x_N)$ as
$$
b({\bf x}):=\left(\sum_{j\not=1}\gamma_j K(x_1-x_j),\cdots, \sum_{j\not=N}\gamma_j K(x_N-x_j)\right).
$$
It can be verified that $b\in \bH^{-1}_\infty(\mathbb{R}^{Nd};\mathbb{R}^{Nd})$ is divergence-free. Consequently, by Theorem \ref{Th2} above, 
the SDE \eqref{HG1} has a solution $\mathbf{X}^{N}_t$ with a density $\rho^N_t$ that satisfies the Liouville equation:
$$
\p_t\rho^N_t=\Delta\rho^N_t-\sum_{i\not=j}\gamma_jK(x_i-x_j)\p_{x_i}\rho^N_t.
$$
In \cite[Proposition 1]{JW18}, the entropy solution of the above Liouville equation is directly constructed. However, the existence of a solution to the SDE \eqref{HG1} appears to be an open question.
In the context of the two-dimensional vortex model, where
$$
K(x)=\frac{(x_2,-x_1)}{|x|^2}=\Big(-\p_{x_1}\arctan\big(\frac{x_2}{x_1}\big),\p_{x_2}\arctan\big(\frac{x_1}{x_2}\big)\Big) \in\bH^{-1}_\infty(\mR^2;\mR^2)
$$ 
known as the Biot-Savart law, and when $\gamma_j$  have the same sign, Takanobu \cite{Ta85} established, 
through a purely probabilistic argument, the existence and uniqueness of a solution that avoids collisions starting from a 
point $\mathbf{x}=(x_1,\cdots,x_N)$ with $x_i\neq x_j$. For general $\gamma_j\in\mathbb{R}$, 
Osada \cite{Os86} showed the same result based on heat kernel estimates for generators in a generalized divergence form, obtained in \cite{Os87}, and on potential theoretical results. 
Note that for $i \neq j$ and $M > 0$, using Krylov's estimate \eqref{SW22} below with $f({\bf x}) = |x_i - x_j|^{-1/N} \1_{|{\bf x}| \leq M}$, we obtain for $p > \frac{Nd}{2}$,
$$
\mathbb{E} \left( \int_0^T |X^{N,i}_t - X^{N,j}_t|^{-1/N} 
\1_{\{|\mathbf{X}^{N}_t| \leq M\}} \, \mathrm{d}t \right) \lesssim \|f\|_p < \infty.
$$
This implies that $(\mathbb{P} \times \mathrm{d}t)(\{(\omega, t) : X^{N,i}_t(\omega) \neq X^{N,j}_t(\omega)\}) = 1$. 
Unfortunately, within our general framework, it remains unclear whether two particles can meet.

We summarize our main contributions as follows:
\begin{itemize}
    \item In  Theorem \ref{Th1}, we establish the weak well-posedness of SDEs with drift $b = b_1 + b_2$, where $b_1 \in \mL^{q_b}_T \bH^{\alpha_b}_{p_b}$ with $\alpha_b \in (-1, -\frac{1}{2})$ and $|b_2(t,x)| \lesssim 1 + |x|$. This result enables us to address cases where $b$ includes Gaussian fields (see Example \ref{KE2} and Subsection \ref{sec63}). Potential applications include SDEs in random environments, such as random directed polymers \cite{DD16}, and self-attracting Brownian motion in a random medium \cite{CC18}.
    \item In Lemma \ref{in:thm01}, we demonstrate the tightness of the laws of solutions to approximating SDEs with arbitrary initial data and divergence-free drifts in $\mL^{q_b}_T \bH^{-1}_{p_b}$, where $2/q_b + d/p_b < 1$. Subsequently, in Theorem \ref{in:main}, we prove the well-posedness of the generalized martingale problem for scaling-supercritical SDEs with divergence-free drift and initial data possessing an $L^2$ density. As an application, we discuss $N$-particle systems with singular interaction kernels, such as the one described in equation \eqref{HG1}, which is closely related to the study of propagation of chaos in vortex systems (see \cite{JW18}).
\end{itemize}

\subsection{Comparison with related works} 
In this subsection we make a detailed comparison with the related well-known results.

{\bf Subcritical case: $\alpha\in(-\frac12,0]$.} 
Our result, Theorem \ref{23thm:48} below, not only encompasses but also extends the work presented in \cite{FIR17, ZZ17} for the range $\alpha\in(-\frac12,0]$. In particular, when $\alpha=\frac12$, 
 our conditions for the coefficient $b$  are not  comparable to the conditions  stated in \cite{ZZ17}, where $b\in \mathbf{H}^{-1/2}_p$ with $p>2d$ is time-independent. Instead, our result allows $b$ to 
 belong to $\mL^q_T\mathbf{B}^{-1/2}_{p,1}$ with $\frac{d}{p}+\frac{2}{q}<\frac12$, which encompasses cases with time integrability index $q\in(4,\infty]$.
 
{\bf Energy method comparison for $\alpha\in(-1,-\frac12]$.} In contrast to the energy method employed in \cite{GP23, GDJP}, our approach differs significantly. We will outline these differences within the context of the following three aspects.
\begin{itemize}
\item  {\bf $\mT^d$ vs $\mR^d$:} In this paper, we study the SDE in the full space $\mathbb{R}^d$. By contrast, the work in \cite{GP23} focuses on the $d$-dimensional torus. Similarly, \cite{GDJP} begins with the $d$-dimensional torus and subsequently uses an $m$-periodic torus to approximate $\mathbb{R}^d$ as $m\to\infty$. The approach in \cite{GP23, GDJP} is motivated by the use of the Lebesgue measure as the initial distribution for the solution. This choice facilitates the application of the forward-backward martingale argument, which is crucial for deriving the It\^o trick introduced later in this part. 
\item {\bf Dimensional dependence:} 
Our supercritical condition, which requires $b\in\mL^q_T\mathbf{H}^{-1}_p$ with $\frac{2}{q}+\frac{d}{p}<1$, is dimension-dependent. This is in contrast to the assumption made in \cite{GP23, GDJP}, where they require $b\in \bB^{-\gamma}_{p,1}$ with $\gamma\in[0,1]$, $\frac{2}{p}\le 1-\gamma$, 
and this assumption is independent of dimension.  It is worth noting that when 
$b$ is time-independent, the case $b=\div A$ with some $A\in L^2(K)\cap L^\infty(K^C)$, where $K$ is a small compact set, is also studied in \cite{GDJP}. This framework is further extended to the time-dependent case in \cite{Gr24}. 

\item {\bf Methodology:} Our methodologies differ fundamentally from each other. 
In \cite{GP23, GDJP},  the key aspect of the energy method is to find a stationary distribution $\mu$ and show the following It\^o trick:
\begin{align*}
\mE_{\mu}\left[\sup_{t\in[0,T]}\left|\int_0^tf(X_s)\dif s\right|^m\right]\lesssim T^{m/2}\|(-\sL_{\rm s})^{-1/2}f\|_{L^p(\mT^d)}^p,
\end{align*}
  where $\sL_{\rm s}:=\frac12(\sL+\sL^*)$ is the symmetric part of the generator $\sL:=\Delta+b\cdot\nabla$, which is exactly the Laplacian $\Delta$ when $\div b=0$.
In our work, we make use of a classical Krylov-type estimate, where the condition $\frac2q+\frac dp<1$ is sharp due to scaling. More precisely, 
let $p',q'$ be the conjugate index of $p,q$, and
$\varphi_t(x)=(2\pi t)^{-d/2}\e^{-|x|^2/(2t)}$ be the Gaussian heat kernel. Clearly, we have
\begin{align*}
    \left|\mE\int_0^T b(t,W_t)\dif t\right|= \left|\int_0^T\!\!\int_{\mR^d} b(t,x)\varphi_t(x)\dif x\dif t\right|
    \leq \|b\|_{\mL^q_T\bH^{-1}_{p}}\|\varphi\|_{\mL^{q'}_T\bH^{1}_{p'}}.
\end{align*} 
Since $\|\nabla \varphi_t\|_{{p'}}=t^{-\frac12-\frac{d}{2p}}\|\nabla \varphi_1\|_{{p'}}$, it is easy to see that $\|\varphi\|_{\mL^{q'}_T\bH^{1}_{p'}}$ is finite if and only if
$\frac2{q'}+\frac{d}{p'}>d+1$. In other words,  the condition $\frac2q+\frac dp<1$ is sharp for  the above Krylov estimate.
Moreover, our method does not necessitate the existence of a stationary distribution, which can be challenging in unbounded spaces.  Thus, on the contrary, in Lemma \ref{in:thm01} below, our method achieves  the tightness without any assumption on the initial data. Additionally, our approach allows us to handle time-dependent drifts efficiently.
Furthermore, our method could be used to deal with the multiplicative noise as well as the degenerate kinetic SDEs, which seems not possible by the energy method developed from \cite{GP23, GDJP}. 
\end{itemize}

{\bf Rough path method.} Recently, numerous studies have been dedicated to investigating singular SDEs driven by fractional Brownian motions. Due to the absence of Markov and martingale properties for fractional Brownian motions, traditional stochastic analysis arguments are no longer applicable. 
In \cite{GG23}, path-by-path well-posedness is established by using the Young integral in rough path theory and the stochastic sewing lemma. It's important to note that this method only works with subcritical drifts. For related results on existence and pathwise uniqueness, see also \cite{BLM23}.

\subsection{Organization of the paper}

In Section 2, we provide a brief overview of classical Besov and Bessel potential spaces and establish a precise Schauder estimate for classical heat equations.

Section 3 is dedicated to the study of the solvability of the linear PDE \eqref{PDE9} with a distributional drift $b$. In the subcritical case, 
we utilize Duhamel's formulation and fixed point theorem, while in the supercritical case, we apply maximum principle and the energy method.

In Section 4, we focus on proving Theorem \ref{Th1}, primarily through the Zvonkin transformation. A crucial step in this process is to establish the following 
uniform Krylov estimate for the approximation solution $X^n$:
for some $m\geq 1$ and any  $f\in \mL^q_T\bH^{\alpha}_{p}$ and $0\leq t_0<t_1\leq T$,
$$
\sup_n\left\|\int^{t_1}_{t_0}f(s,X^n_s)\dif s\right\|_{L^m(\Omega)}\leq C(t_1-t_0)^{\frac{1+\theta_0}2}
\|f\|_{\mL^q_T\bH^{\alpha}_p},
$$
where $\theta_0>0$ and $(\alpha,p,q)$ satisfies \eqref{BQ1}. This estimate ensures that $A^b_t$ is a zero-energy process,
which, in turn, enables us to apply the generalized It\^o's formula and derive the effectiveness of Zvonkin's transformation. 
Subsequently, the uniqueness and heat kernel estimates follow by the corresponding results of the transformed SDEs.

Section 5 is dedicated to the proof of Theorem \ref{Th2}. To establish the existence of a weak solution, we need to establish two types of Krylov's estimates (as presented in Lemma \ref{Le51} and Lemma \ref{Le54}). For the well-posedness of the associated generalized martingale problem, which is crucial for {\bf (A)}-{\bf (C)}, our approach primarily revolves around solving PDE \eqref{PDE9}

In Section 6, we apply the results to the diffusion in random environments. In particular, vector-valued Gaussian random fields are provided to illustrate our results.

Finally, in the appendix, we present a decomposition result for a weighted distribution, originally established by Gubinelli and Hofmanov\`a \cite[Lemma 2.4]{GH19}. For the convenience of readers, we include a detailed proof.

Throughout this paper, we use $C$ with or without subscripts to denote constants, whose values
may change from line to line. We also use $:=$ to indicate a definition and set
$$
a\vee b:=\max(a,b),\ \ a\wedge b:=\min(a,b).
$$ 
By $A\lesssim B$, we mean that for some constant $C\geq 1$, $A\leq C B$.

\section{Preliminaries}

In this section, we recall the definition and some properties of Besov and Bessel potential spaces. Additionally, we establish a Schauder estimate for the heat equation in Besov spaces. 
This estimate will be used to solve the Kolmogorov equation in the subcritical case.

Let $\sS({\mR^d})$ be the Schwartz space of all rapidly decreasing functions on ${\mR^d}$, and $\sS'({\mR^d})$
the dual space of $\sS({\mR^d})$ called Schwartz generalized function (or tempered distribution) space. Given $f\in\sS({\mR^d})$,
the Fourier transform $\hat f$ and inverse Fourier transform $\check f$ are defined, respectively, by
\begin{align*}
\hat f(\xi):=\frac{1}{(2\pi)^{d/2}}\int_{{\mR^d}} \e^{-{\rm i}\xi\cdot x}f(x)\dif x, \quad
\check f(x):=\frac{1}{(2\pi)^{d/2}}\int_{{\mR^d}} \e^{{\rm i}\xi\cdot x}f(\xi)\dif\xi.
\end{align*}
For $q\in[1,\infty]$ and a normed space $\mB$, we shall simply denote
$$
\mL^q_T\mB:=L^q([0,T];\mB),\ \ \mL^q_T:=\mL^q_T L^q,
$$
where $(L^p,\|\cdot\|_p)$ is the usual $L^p$-space in $\mR^d$. Moreover, we use $C_b=C_b(\mR^d)$ to denote the Banach space of all bounded continuous functions on $\mR^d$,
and $C^\infty_b=C^\infty_b(\mR^d)$ the space of all smooth functions with bounded derivatives of all orders. Here, we clarify that the notation $\mathbb{L}^q_T L^q = \mathbb{L}^q_T$ is used throughout the rest of the paper purely for conciseness, making the norm $\|\cdot\|_{\mathbb{L}^q_T}$ more compact.

\subsection{Besov and Bessel potential spaces} To 
{define Besov spaces}, we first introduce dyadic partitions of unity.  Let $\phi_{-1}$ be a symmetric
nonnegative $C^{\infty}-$function on $\mathbb{R}^d$ with
$$
\phi_{-1}(\xi)=1\ \mathrm{for}\ \xi\in B_{1/2}\ \mathrm{and}\ \phi_{-1}(\xi)=0\ \mathrm{for}\ \xi\notin B_{2/3}.
$$
For   $j\geq 0$, we define
\begin{align}\label{Phj}
\phi_j(\xi):=\phi_{-1}(2^{- (j+1)}\xi)-\phi_{-1}(2^{-j}\xi).
\end{align}
By 
definition, one sees that for $j\geq 0$, $\phi_j(\xi)=\phi_0(2^{-j}\xi)$ and
$$
\mathrm{supp}\,\phi_j\subset B_{2^{j+2}/3}\setminus B_{2^{j-1}},\quad\sum^n_{j=-1}\phi_j(\xi)=\phi_{-1}(2^{-(n+1)}\xi)\to 1,\quad n\to\infty.
$$
\bd
For  $j\geq -1$, the Littlewood-Paley block operator $\cR_j$ is defined on $\sS'({\mR^d})$ by
$$
\cR_jf(x):=(\phi_j\hat f)\check{\,\,}(x)=\check\phi_j* f(x),
$$
with the convention $\cR_j\equiv0$ for $j\leq-2$.
In particular, for $j\geq 0$,
\begin{align}\label{Def2}
\cR_jf(x)=2^{jd}\int_{{\mR^d}}\check\phi_0(2^{j}y) f(x-y)\dif y.
\end{align}
\ed

For $j\geq -1$, by definition it is easy to see that
\begin{align}\label{KJ2}
\cR_j=\cR_j\widetilde\cR_j,\ \mbox{where }\ \widetilde\cR_j:=\cR_{j-1}+\cR_{j}+\cR_{j+1},
\end{align}
and $\cR_j$ is symmetric in the sense that
$$
\<g, \cR_j f\>=\< f,\cR_jg\>,\ \ f,g\in\sS'(\mR^d),
$$
where $\<\cdot,\cdot\>$ stands for the dual pair between $\sS'(\mR^d)$ and $\sS(\mR^d)$.

Now we recall the definition of Besov spaces (see \cite{BCD11}).
\bd\label{Def25}
Let $p,q\in[1,\infty]$ and $s\in\mR$. The Besov space $\bB^{s}_{p,q}$ is defined by
$$
\bB^{s}_{p,q}:=\left\{f\in\sS'(\mR^d): \|f\|_{\bB^{s}_{p,q}}:=
\left(\sum_{j\ge -1} 2^{sjq}\|\cR_j f\|_p^q\right)^{1/q}<\infty\right\}.
$$
For $p=q=\infty$, we simply denote $\bC^s:=\bB^{s}_{\infty,\infty}$.
\ed
\br
For $s\in(0,1)$, an equivalent characterization of $\bB^s_{p,q}$ is given by (see \cite[p74, Theorem 2.36]{BCD11} or \cite[Theorem 2.7]{HZZZ21})
$$
\|f\|_{\bB^{s}_{p,q}}\asymp\left(\int_{|h|\leq 1}
\left(\frac{\|f(\cdot+h)-f(\cdot)\|_p}{|h|^{s}}\right)^q\frac{\dif h}{|h|^d}\right)^{1/q}
+\|f\|_p.
$$
In particular,  for any $s\in(0,1)$ and $p\in[1,\infty]$, there is a constant $C=C(s,d,p)>0$ such that
$$
\|f(\cdot+h)-f(\cdot)\|_p\leq C\|f\|_{\bB^{s}_{p,\infty}} (|h|^s\wedge 1),
$$
and for any $s_0\in\mR$,
\begin{align}\label{Ho1}
\|f(\cdot+h)-f(\cdot)\|_{\bB^{s_0}_{p,\infty}}\leq C\|f\|_{\bB^{s_0+s}_{p,\infty}} (|h|^s\wedge 1).
\end{align}
From this estimate, one sees that $\bC^s$ coincides with the classical H\"older space.
\er
For $(\alpha,p)\in(\mR\setminus\mZ)\times[1,\infty]$, let $\bH^\alpha_p$ be the Bessel potential space defined by
$$
\bH^\alpha_p:=\big\{f\in \sS'(\mR^d): \|f\|_{\bH^\alpha_p}:=\|(\mI-\Delta)^{\alpha/2}f\|_p<\infty\big\},
$$
where  $(\mI-\Delta)^{\alpha/2}f$ is defined through Fourier's  transform
$$
(\mI-\Delta)^{\alpha/2}f:=\big((1+|\cdot|^2)^{\alpha/2}\hat f\big)^{\check{}}.
$$
For $\alpha=0,1,2,\cdots$ and $p\in[1,\infty]$, we define 
$$
\bH^{\alpha}_p:=\big\{f\in\sS'(\mR^d): \|f\|_{\bH^{\alpha}_p}:=\|f\|_p+\|\nabla^\alpha f\|_{p}<\infty\big\},
$$
and for $\alpha=-1,-2,\cdots$ and $p\in[1,\infty]$,
$$
\bH^{\alpha}_p:=(\bH^{-\alpha}_{p/(p-1)})'.
$$
Note that we do not use the space $\bH^\alpha_1$ for $\alpha<0$.
Moreover, for integer $\alpha\in\mN$ and $p\in(1,\infty)$, an equivalent norm in $\bH^\alpha_p$ is given by (cf. \cite[p135, Theorem 3]{St})
\begin{align}\label{AB22}
\|f\|_{\bH^\alpha_p}\asymp\|(I-\Delta)^{\alpha/2}f\|_p.
\end{align}

Below we recall some well-known facts about Besov spaces and Bessel potential spaces, where \eqref{AB2} and \eqref{EW-2} below are easily derived by definition.
\bl\label{lemB1}
\begin{enumerate}[(i)]
 \item For any $p\in[1,\infty]$,  $s'>s$ and $q\in[1,\infty]$, it holds that
\begin{align}\label{AB2}
\bB^{0}_{p,1}\hookrightarrow L^p\hookrightarrow\bB^{0}_{p,\infty},\ \ \bB^{s'}_{p,\infty}\hookrightarrow \bB^{s}_{p,1}\hookrightarrow \bB^{s}_{p,q}\hookrightarrow \bB^{s}_{p,\infty}.
\end{align}
\item For any  $\alpha\in\mR$ and $p\in(1,\infty)$, we have the following embedding (see \cite[Theorem 6.4.4]{BL76})
\begin{align}\label{Con}
\left\{
\begin{aligned}
&\bB^\alpha_{p,2}\hookrightarrow \bH^\alpha_p,\ \ p\in[2,\infty),\\
&\bH^\alpha_p\hookrightarrow\bB^\alpha_{p,2} ,\ \ p\in(1,2],
\end{aligned}
\right.
\end{align}
and for  $\alpha\in\mR$ and $p\in(1,\infty]$ (see \cite[Theorem 6.2.4]{BL76}),
$$
\bB^{\alpha}_{p,1}\subset\bH^\alpha_p\subset\bB^{\alpha}_{p,\infty}.
$$
\item For $1\leq p_1\leq p\leq\infty$, $q\in[1,\infty]$ and $\alpha\leq\alpha_1-\tfrac{d}{p_1}+\tfrac{d}{p}$, it holds that (see \cite[Theorem 6.5.1]{BL76})
\begin{align}\label{Sob1}
\|f\|_{\bB^{\alpha}_{p,q}}\lesssim\|f\|_{\bB^{\alpha_1}_{p_1,q}},
\end{align}
and for $1<p_1\leq p<\infty$,
\begin{align}\label{Sob11}
\|f\|_{\bH^{\alpha}_{p}}\lesssim\|f\|_{\bH^{\alpha_1}_{p_1}}.
\end{align}
\item Let $\alpha\in[-1,0)$ and $p\in(1,\infty)$. For any $p_1,p_2\in(1,\infty)$ satisfying
$$
p_1\geq p,\ p_2\geq\tfrac{p_1}{p_1-1},\ \ \tfrac{1}{p}\leq\tfrac{1}{p_1}+\tfrac{1}{p_2}<\tfrac{1}{p} -\tfrac{\alpha}{d},
$$ 
there is a constant $C>0$ such that
for all $f\in \bH^{\alpha}_{p_1}$ and $g\in \bH^{-\alpha}_{p_2}$ (see \cite[Lemma 2.2]{ZZ17}),
\begin{align}\label{EW2}
\|fg\|_{\bH^{\alpha}_p}\le C\|f\|_{\bH^{\alpha}_{p_1}} \|g\|_{\bH^{-\alpha}_{p_2}},
\end{align}
and
\begin{align}\label{EW-2}
\|fg\|_{\bH^{-1}_\infty}\leq\|f\|_{\bH^{-1}_\infty} \|g\|_{\bH^1_\infty}.
\end{align}
\end{enumerate}
\el

Next we recall the well-known Bony decomposition and related paraproduct estimates.
Let $S_k$ be the cut-off low frequency operator defined by
\begin{align}\label{EM9}
S_kf:=\sum_{j=-1}^{k-1}\cR_j f\to f,\ \ k\to\infty.
\end{align}
For $f,g\in\sS'(\R^d)$, define
$$
f\prec g:=\sum_{k\geq -1} S_{k-1}f\cR_k g,\ \ f\circ g:=\sum_{|i-j|\leq1}\cR_i f\cR_jg.
$$
The Bony decomposition of $fg$ is formally given by (cf. \cite{BCD11})
\begin{align}\label{Bony}
fg=f\prec g+ f\circ g+g\prec f=:f\preceq g+f\succ g.
\end{align}

We recall the following paraproduct estimates (cf. \cite[Theorem 2.82 and Theorem 2.85]{BCD11}).
\bl\label{Le24}
Let $p,p_1,p_2,q, q_1,q_2\in[1,\infty]$ with  $\frac1p=\frac1{p_1}+\frac1{p_2}$ and $\frac1q=\frac1{q_1}+\frac1{q_2}$ and $\alpha,\beta\in\mR$.
\begin{enumerate}[(i)]
\item If $\beta<0$, then there is a constant $C=C(\alpha,\beta,d,p,q,p_1,q_1,p_2,q_2)>0$ such that
$$
\|f\prec g\|_{\bB^{\alpha+\beta}_{p,q}}\le C\|f\|_{\bB^{\beta}_{p_1,q_1}} \|g\|_{\bB^{\alpha}_{p_2,q_2}}.
$$
Moreover, for $\beta=0$, we have
$$
\|f\prec g\|_{\bB^{\alpha}_{p,q}}\le C\|f\|_{p_1} \|g\|_{\bB^{\alpha}_{p_2,q}}.
$$

\item If $\alpha+\beta>0$,  then there is a constant $C=C(\alpha,\beta,d,p,q,p_1,q_1,p_2,q_2)>0$ such that
$$
\|f\circ g\|_{\bB^{\alpha+\beta}_{p,q}}\le C\|f\|_{\bB^{\beta}_{p_1,q_1}} \|g\|_{\bB^{\alpha}_{p_2,q_2}}.
$$
Moreover, when $\alpha+\beta=0$ and $q=1$, we have
\begin{align}\label{Cir1}
\|f\circ g\|_{\bB^{0}_{p,\infty}}\le C\|f\circ g\|_p\le C\|f\|_{\bB^{-\alpha}_{p_1,q_1}} \|g\|_{\bB^{\alpha}_{p_1,q_2}}.
\end{align}
\end{enumerate}
\el
\begin{proof}
We only prove \eqref{Cir1}. By definition and H\"older's inequality, we have
\begin{align*}
\|f\circ g\|_{\bB^{0}_{p,\infty}}&\lesssim\|f\circ g\|_{p}
\leq\sum_{|i-j|\leq1}\|\cR_i f\cR_jg\|_p\leq\sum_{|i-j|\leq1}\|\cR_if\|_{p_1}\|\cR_jg\|_{p_2}\\
&\lesssim\left(\sum_{i}2^{-q_1\alpha i}\|\cR_if\|^{q_1}_{p_1}\right)^{1/q_1}
\left(\sum_{j}2^{q_2\alpha j}\|\cR_if\|^{q_2}_{p_2}\right)^{1/q_2}=\|f\|_{\bB^{-\alpha}_{p_1,q_1}} \|g\|_{\bB^{\alpha}_{p_1,q_2}}.
\end{align*}
The proof is complete.
\end{proof}

Now we discuss a decomposition of $b\cdot\nabla u$, that arises in the study of PDE \eqref{PDE0}, 
where $b\in\sS'(\mR^d)$ is a distribution-valued vector field and $u\in\sS(\mR^d)$.
By Bony's decomposition, we can write
\begin{align}\label{08:21}
b\cdot\nabla u=\{\div(b\preceq u)+ b \succ \nabla u\}-\div b\preceq u=: b\odot \nabla u-\div b\preceq u.
\end{align}
In particular, if $\div b=0$, then
$$
b\cdot\nabla u=b\odot \nabla u.
$$
We remark that similar decomposition \eqref{08:21} has been used in \cite[Lemma 2.6]{GP23}. 

We have the following regularity estimate.
\bl\label{Le26}
For $\alpha\in(-1,0]$ and $p,q\in[1,\infty]$, there is a constant $C=C(d,\alpha,p,q)>1$ such that
\begin{align}\label{08:20}
\|b\odot\nabla u\|_{\bB^{\alpha}_{p,q}}\leq C\|b\|_{\bB^{\alpha}_{p,q}}\|u\|_{\bH^1_{\infty}}\leq 
C\|b\|_{\bB^{\alpha}_{p,q}}\|u\|_{\bB^{1+d/p}_{p,1}}
\end{align}
and
\begin{align}\label{08:210}
\|\div b\preceq u\|_{\bB^{0}_{p/2,\infty}}\le C\|\div b\|_{\bB^{-\alpha-2}_{p,q/(q-1)}}\|u\|_{\bB^{2+\alpha}_{p,q}}.
\end{align}
\el
\begin{proof}
Since $\alpha\in(-1,0]$, by Lemma \ref{Le24} we have
$$
\|\div(b\preceq u)\|_{\bB^{\alpha}_{p,q}}
\lesssim \|b\preceq u\|_{\bB^{1+\alpha}_{p,q}}\lesssim\|b\|_{\bB^{\alpha}_{p,q}}\|u\|_{\bB^{1}_{\infty,\infty}}
$$
and
$$
\|b\succ \nabla u\|_{\bB^{\alpha}_{p,q}}\lesssim\|b\|_{\bB^{\alpha}_{p,q}}\|\nabla u\|_{\infty}\lesssim \|b\|_{\bB^{\alpha}_{p,q}}\|u\|_{\bH^{1}_{\infty}}.
$$
Hence, by embedding  \eqref{Con} and \eqref{Sob1},
$$
\|b\odot\nabla u\|_{\bB^{\alpha}_{p,q}}\lesssim \|b\|_{\bB^{\alpha}_{p,q}}
\|u\|_{\bH^{1}_{\infty}}\lesssim\|b\|_{\bB^{\alpha}_{p,q}}\|u\|_{\bB^{1}_{\infty,1}}
\lesssim\|b\|_{\bB^{\alpha}_{p,q}}\|u\|_{\bB^{1+d/p}_{p,1}}.
$$
For \eqref{08:210}, by \eqref{AB2} and  Lemma \ref{Le24}-(i) with $p_1=p_2=p$ and $(q_1,q_2)=(q/(q-1),q)$, we have
$$
    \|\div b\prec u\|_{\bB^0_{p/2,\infty}}\lesssim \|\div b\prec u\|_{\bB^0_{p/2,1}}\lesssim \|\div b\|_{\bB^{-\alpha-2}_{p,q/(q-1)}}\|u\|_{\bB^{2+\alpha}_{p,q}}.
$$
For the resonant part, by \eqref{Cir1} with $p_1=p_2=p$ and $(q_1,q_2)=(q/(q-1),q)$, we have
$$
    \|\div b\circ u\|_{\bB^0_{p/2,\infty}}\lesssim \|\div b\|_{\bB^{-\alpha-2}_{p,q/(q-1)}}\|u\|_{\bB^{2+\alpha}_{p,q}}.
$$
Combining the above two estimates, we complete the proof of \eqref{08:210}. 
\end{proof}

We also have the following easy result.
\bl\label{Le27}
If $b, \div b\in\bB^{-1}_{\infty,2}$, then $b\in\cB$, where $\cB$ is defined in \eqref{CB1}. 
\el
\begin{proof}
By definition \eqref{08:21}, \eqref{Con} and Lemma \ref{Le24}, we have
\begin{align*}
\|b\odot\nabla u\|_{\bH^{-1}_2}&\leq\|\div(b\preceq u)\|_{\bH^{-1}_2}+\|b\succ\nabla u\|_{\bH^{-1}_2}\lesssim\|b\preceq u\|_2+\|b\succ\nabla u\|_{\bB^{-1}_{2,2}}\\
&\lesssim \|b\prec u\|_{\bB^{0}_{2,1}}+\|b\circ u\|_2+\|b\|_{\bB^{-1}_{\infty,2}}\|\nabla u\|_2\lesssim \|b\|_{\bB^{-1}_{\infty,2}} \big(\|u\|_{\bB^{1}_{2,2}}+\|u\|_{\bH^1_2}\big),
\end{align*}
and
\begin{align*}
\|\div b\preceq u\|_{\bH^{-1}_2}\lesssim \|\div b\prec  u\|_{\bB^{-1}_{2,2}}+\|\div b\circ u\|_{\bB^{-1}_{2,2}}
\lesssim \|\div b\|_{\bB^{-1}_{\infty,2}}\|u\|_{\bB^{1}_{2,2}}.
\end{align*}
Combining the above two estimates and by the definition of $\cB$ and \eqref{Con}, we get $b\in\cB$.
\end{proof}

\subsection{Schauder's estimate for heat equation in Besov spaces}
In this subsection we study the regularity estimates for heat equation in Besov spaces. 
Let  $(P_t)_{t\geq 0}$ be the Gaussian heat semigroup given by
\begin{align*}
P_tf=\varphi_t*f\quad \text{where}\quad \varphi_t(x)=(4\pi t)^{-\frac{d}{2}}\e^{-\frac{|x|^2}{4t}}.
\end{align*}
For $\lambda\ge0$ and a time-dependent distribution $f_t(\cdot):\mR_+\to\sS'(\mR^d)$, we define
\begin{align}\label{LA1}
u^\lambda_t(x):=\sI^\lambda_t(f)(x):=\int_0^t \e^{-\lambda(t-s)}P_{t-s}f_s(x)\dif s,\ \ t>0.
\end{align}
In particular, in the distributional sense,
\begin{align}\label{LA4}
\p_t u^\lambda=(\Delta-\lambda)u^\lambda+f.
\end{align}
We now show the following Schauder estimate, which shall play a basic role in the subcritical case.
\bl\label{08:202}
Let $\alpha\in\mR$, $1\le q\leq q'\le \infty$ and $1\le p\le p'\le \infty$. Define
$$
\alpha':=\tfrac2q-\tfrac2{q'}+\tfrac{d}{p}-\tfrac{d}{p'}.
$$ 
For any $T>0$ and  $\gamma\in[q,q']$, there is a constant $C=C(T,\alpha,d,p,q,p',q', \gamma)>0$ such that 
\begin{align}\label{08:203}
\|\sI^\lambda(f)\|_{\mL^{q'}_T\bB^{2+\alpha-\alpha'}_{p',\gamma}}\le C\|f\|_{\mL^q_T\bB^{\alpha}_{p,\gamma}},\ \ \lambda\ge0,
\end{align}
and for $\theta\in(0,2+\frac2{q'}-\frac2q)$,
\begin{align}\label{08:204}
\|\sI^\lambda(f)\|_{\mL^{q'}_T\bB^{2+\alpha-\alpha'-\theta}_{p',1}}\lesssim_{C,\theta} (1+\lambda)^{-\frac\theta2}\|f\|_{\mL^q_T\bB^{\alpha}_{p,\infty}},\quad \text{$\lambda\ge0$}.
\end{align}
\el
\begin{proof}
By virtue of the embedding \eqref{Sob1} and the equivalence (see, for instance, \cite[Proposition 2.78]{BCD11}),
$$
\|(1-\Delta)^{\alpha/2}f\|_{\bB^{0}_{p,q}} \asymp \|f\|_{\bB^{\alpha}_{p,q}},
$$
we may, without loss of generality, restrict our proof of \eqref{08:203} and \eqref{08:204} to the case $ p' = p $ and $ \alpha = 0 $.

For any $j\ge-1$, by Young's inequality and Bernstein's inequality (see \cite[Lemma 2.1, p52]{BCD11}), one sees that
\begin{align*}
\|\cR_j\sI^\lambda_t(f)\|_{p}&\le \int_0^t \e^{-\lambda(t-s)}\|\cR_j (\varphi_{t-s}*f_s)\|_{p}\dif s\\
&\!\!\stackrel{\eqref{KJ2}}{=} \int_0^t \e^{-\lambda(t-s)}\|\cR_j \varphi_{t-s}*\wt\cR_j f_s\|_{p}\dif s\\
&\le \int_0^t \e^{-\lambda(t-s)}\|\cR_j \varphi_{t-s}\|_{1}\|\widetilde{\cR}_j f_s\|_{p}\dif s\\
&\lesssim  \int_0^t \e^{-\lambda(t-s)}\|\cR_j \varphi_{t-s}\|_{1}\|\widetilde{\cR}_j f_s\|_{p}\dif s.
\end{align*}
Let $r\in[1,\infty]$ be defined by $1+\frac1{q'}=\frac 1{r}+\frac 1 q$. 
By Young's inequality for time variable, we have
\begin{align*}
\left(\int_0^T\|\cR_j\sI^\lambda_t(f)\|_{p}^{q'}\dif t\right)^{\frac1{q'}} \lesssim 2^{(d-\frac dr)j}\left(\int_0^T \e^{-r'\lambda s}\|\cR_j \varphi_s\|_1^{r}\dif s\right)^{\frac{1}{r}}\left(\int_0^T \|\widetilde{\cR}_jf_s\|_p^{q}\dif s\right)^{\frac{1}{q}}.
\end{align*}
It is well-known that for any $\theta\geq 0$, there is a constant $C=C(\theta,d)>0$ such that  (see \cite{HWZ})
$$
\|\cR_j \varphi_s\|_1\le C[1\wedge (2^{2j} s)^{-\theta}],\ \ j\geq -1,\ s\in(0,T].
$$
If $r<\infty$, then by the change of variable, we have
\begin{align*}
\left(\int_0^T \e^{-r\lambda s}\|\cR_j \varphi_s\|_1^{r}\dif s\right)^{\frac1{r}}
&\lesssim \left(\int_0^\infty \e^{-r\lambda s} [1\wedge (2^{2j}s)^{-2}]\dif s\right)^{\frac1{r}}\\
&=\left(2^{-2j}\int_0^\infty \e^{-r2^{-2j}\lambda s} [1\wedge s^{-2}]\dif s\right)^{\frac1{r}}\\
&\leq\left(2^{-2j}(1\wedge(r2^{-2j}\lambda)^{-1}) \right)^{\frac1{r}}.
\end{align*}
If $r=\infty$, then we directly have
$$
\sup_{j\geq -1}\sup_{s\geq 0}\e^{-\lambda s}\|\cR_j \varphi_s\|_1\lesssim 1.
$$
Thus combining the above calculations, we obtain that for any $j\geq -1$ and $\lambda\geq 0$,
\begin{align*}
&\|\cR_j\sI^\lambda_\cdot(f)\|_{\mL^{q'}_T L^{p}}=\left(\int_0^T\|\cR_j\sI^\lambda_t(f)\|_{p}^{q'}\dif t\right)^{\frac1{q'}} \\
&\qquad\lesssim 2^{-\frac2{r}j}(1\wedge(2^{-2j}\lambda)^{-\frac1{r}})\left(\int_0^T \|\widetilde{\cR}_jf_s\|_p^{q}\dif s\right)^{\frac{1}{q}}\\
&\qquad=2^{(\alpha'-2)j}(1\wedge(2^{-2j}\lambda)^{-\frac1{r}})\left(\int_0^T \|\widetilde{\cR}_jf_s\|_p^{q}\dif s\right)^{\frac{1}{q}},
\end{align*}
where $\alpha'=2-\frac2{r}=\tfrac2q-\tfrac2{q'}$. 

Now, for any {$\gamma\in[1,q']$} and $\theta\geq 0$, by the definition of Besov norm and Minkowskii's inequality,
\begin{align*}
\|\sI^\lambda(f)\|_{\mL^{q'}_T\bB^{2-\alpha'-\theta}_{p,\gamma}}
&=\left\|\left(\sum_{j\geq-1} 2^{\gamma(2-\alpha'-\theta)j}\|\cR_j\sI^\lambda_\cdot(f)\|^\gamma_{p}\right)^{\frac1\gamma}\right\|_{\mL^{q'}_T}\\
&\leq \left(\sum_{j\geq-1} 2^{\gamma(2-\alpha'-\theta)j}\|\cR_j\sI^\lambda_\cdot(f)\|^\gamma_{\mL^{q'}_T L^{p}}\right)^{\frac1\gamma}\\
&\lesssim  \left(\sum_{j\geq-1} 2^{\gamma(\alpha-\theta) j}(1\wedge(2^{-2j}\lambda)^{-\frac{\gamma}{r}}) \left(\int_0^T \|\widetilde{\cR}_jf_s\|_p^{q}\dif s\right)^{\frac{\gamma}{q}}\right)^{\frac1\gamma}\\
&=\left(\sum_{j\geq-1} 2^{-\gamma\theta j}(1\wedge(2^{-2j}\lambda)^{-\frac{\gamma}{r}}) \left(\int_0^T 2^{q\alpha j}\|\widetilde{\cR}_jf_s\|_p^{q}\dif s\right)^{\frac{\gamma}{q}}\right)^{\frac1\gamma}.
\end{align*}
In particular, for $\theta=0$ and  $\gamma\in[q,q']$, 
by $\left(\sum_{j}|a_j|^{\gamma}\right)^{1/{\gamma}}\le\left(\sum_{j}|a_j|^q\right)^{1/q}$, we have
$$
\|\sI^\lambda(f)\|_{\mL^{q'}_T\bB^{2-\alpha'-\theta}_{p,\gamma}}
\lesssim \left(\sum_{j\geq-1} 
\left(\int_0^T 2^{q\alpha j}\|{\cR}_jf_s\|_p^{q}\dif s\right)^{\frac{\gamma}{q}}\right)^{\frac1\gamma}
\leq
\|f\|_{\mL^q_T\bB^{\alpha}_{p,\gamma}}.
$$
For $\theta>0$, if we take $\gamma=1$, then
\begin{align*}
\|\sI^\lambda(f)\|_{\mL^{q'}_T\bB^{2-\alpha'-\theta}_{p,1}}
\lesssim \sum_{j\geq-1} 2^{-\theta j}(1\wedge(2^{-2j}\lambda)^{-\frac{1}{r}})\|f\|_{\mL^q_T\bB^{0}_{p,\infty}}.
\end{align*}
Note that for $\theta\in(0,\frac{2}{r})=(0,2+\frac2{q'}-\frac2q)$,
$$
\sum_{j\geq-1} 2^{-\theta j}(1\wedge(2^{-2j}\lambda)^{-\frac{1}{r}})\lesssim 1\wedge\lambda^{-\frac\theta2}\lesssim(1+\lambda)^{-\frac\theta2}.
$$
The proof is complete.
\end{proof}

\br
Note that \eqref{08:204} does not hold for $\theta=0$. Indeed, $\bB^2_{2,1}\subsetneqq \bB^2_{2,2}$, $\bB^0_{2,2}\subsetneqq\bB^0_{2,\infty}$, and by \eqref{Con} and Fourier's transform,
$$
\|\sI^\lambda(f)\|_{\mL^2_T\bB^2_{2,2}}\asymp \|\sI^\lambda(f)\|_{\mL^2_T\bH^2_2}\asymp\|f\|_{\mL^2_T\bH^0_{2}}\asymp\|f\|_{\mL^2_T\bB^0_{2,2}},
$$
where $\asymp$ means that both sides are comparable up to a constant.
\er
\br
Let $1\leq q\leq q'\leq\infty$ and $1\le p\le p'\le \infty$ and $\theta\in(0,2-\frac2q)$.
By \eqref{08:203}, \eqref{08:204} and $\bB^{\alpha}_{p,q}\hookrightarrow \bB^{\alpha}_{p,\infty}$, we immediately have
\begin{align}\label{08:205}
\|\sI^\lambda(f)\|_{\mL^{q'}_T\bB^{2+\alpha-(\frac2q-\frac2{q'}+\frac{d}{p}-\frac{d}{p'})}_{p',q}}+
 (1+\lambda)^{\frac{\theta}2}\|\sI^\lambda(f)\|_{\mL^\infty_T\bB^{2+\alpha-(\frac2q+\frac{d}{p}-\frac{d}{p'})-\theta}_{p',1}}
\lesssim\|f\|_{\mL^q_T\bB^{\alpha}_{p,q}}.
\end{align}
\er

\section{A study of PDEs with distribution drifts}\label{sec:PDE}

In this section,  our goal is to investigate the solvability of the nonhomogeneous PDE:
\begin{align}\label{S1:PDE}
\p_t u=\Delta u-\lambda u+b\cdot \nabla u+f,\quad u_0\equiv0,
\end{align}
where $\lambda\geq 0$, $b$ and $f$ are distributions. 
We will establish the well-posedness of PDE \eqref{S1:PDE} in two cases: the subcritical case and the supercritical case, utilizing different methods.
In the subcritical case, we employ the Duhamel formula, while in the supercritical case, we utilize the maximal principle and the energy method.
Throughout this section we fix a terminal time $T>0$.

\subsection{Subcritical case: $\frac{d}{p}+\frac{2}{q}<1+\alpha$}

In this section, we show the existence and uniqueness of weak solutions to PDE \eqref{S1:PDE} by using Lemma \ref{Le26}
under the following subcritical conditions:
\begin{itemize}
\item[{\bf(H$^{\rm sub}_b$)}] Let $(\alpha_b,p_b,q_b)\in(-1,-\frac12]\times [2,\infty]^2$ with
$\tfrac d{p_b}+\tfrac2{q_b}< 1+\alpha_b$. Suppose that
\begin{align*}
\kappa_b:=\|b\|_{\mL^{q_b}_T\bB^{\alpha_b}_{p_b,q_b}}<\infty\quad \text{and}\quad \tilde \kappa_b:=\|\div b\|_{\mL^{q_b}_T\bB^{-2-\alpha_b}_{p_b,q_b/(q_b-1)}}<\infty.
\end{align*}
\end{itemize}

\br\label{KP1}
Since $\frac d{p_b}+\frac2{q_b}<1+\alpha_b$ and $\alpha_b\in(-1,-\frac12]$, it is easy to see that for any $q\geq\tfrac{q_b}2$,
$$
2+\alpha_b+\tfrac d{p_b}-\tfrac{2}{q_b}<2-\tfrac4{q_b}\leq 2-\tfrac2q.
$$
\er

\br
Suppose that $b\in \mL^{q_b}_T\bB^{-1/2}_{p_b,1}$ with $\frac d{p_b}+\frac2{q_b}<\frac12$. Then {\bf(H$^{\rm sub}_b$)} holds for $\alpha_b=-\frac12$. Indeed, by \eqref{AB2}
we have
$$
\kappa_b=\|b\|_{\mL^{q_b}_T\bB^{-1/2}_{p_b,q_b}}\lesssim\|b\|_{\mL^{q_b}_T\bB^{-1/2}_{p_b,1}},\ \ \tilde \kappa_b=\|\div b\|_{\mL^{q_b}_T\bB^{-3/2}_{p_b,q_b/(q_b-1)}}\lesssim\|b\|_{\mL^{q_b}_T\bB^{-1/2}_{p_b,1}}.
$$
If $b$ is time-independent, then it was covered by the condition $b\in\bH^{-1/2}_{p_b}$ with $p_b>2d$ in \cite{ZZ17}.
\er
For convenience of notations,  we introduce the following parameter set for later use:
\begin{align}\label{Pa1}
\Theta:=(T,d,\alpha_b,p_b,q_b,\kappa_b,\tilde \kappa_b).
\end{align}
Now we show the following result.

\bt\label{S5:thm01}
Let $(\alpha,p,q)\in [-1,0]\times[1,\infty]^2$ satisfy
$$
\tfrac{q_b}2\leq q\leq q_b,\ \ p\leq p_b,\ \ \alpha-\tfrac dp-\tfrac2q\geq\alpha_b-\tfrac d{p_b}-\tfrac2{q_b}.
$$
Under {\bf(H$^{\rm sub}_b$)},
there is a $\lambda_0=\lambda_0(\Theta)>0$ such that for any $\lambda\ge\lambda_0$ and $f\in \mL^q_T\bB^{\alpha}_{p,q}$, 
there is a unique $u^\lambda$
solving the integral equation $u^\lambda_t=\sI^\lambda_t(b\cdot \nabla u^\lambda+f),$
where $\sI^\lambda$ and $b\cdot\nabla u^\lambda$ are defined by \eqref{LA1} and \eqref{08:21} respectively,
and such that
for any $\theta\in(0,2+\alpha_b+\frac d{p_b}-\frac{2}{q_b}]$, 
\begin{align}\label{S5:AA07}
\|u^\lambda\|_{\mL^{q_b}_T\bB^{2+\alpha_b}_{p_b,q_b}}+(1+\lambda)^{\frac\theta2}
\|u^\lambda\|_{\mL^\infty_T\bB^{2+\alpha_b-2/q_b-\theta}_{p_b,1}}
\leq C(\Theta,\theta,\alpha,p,q)\|f\|_{\mL^q_T\bB^{\alpha}_{p,q}}.
\end{align}
\et
\begin{proof}
We only  show the a priori estimate \eqref{S5:AA07}. The existence follows by standard Picard's iteration.
For simplicity of notations, we drop the superscript $\lambda$ over $u^\lambda$, and write for $\lambda,\theta\geq 0$,
$$
\|g\|_{\cS^\lambda_\theta}:=\|g\|_{\mL^{q_b}_T\bB^{2+\alpha_b}_{p_b,q_b}}+(1+\lambda)^{\frac\theta2}
\|g\|_{\mL^\infty_T\bB^{2+\alpha_b-2/q_b-\theta}_{p_b,1}}.
$$
By Schauder's estimate \eqref{08:205}, \eqref{AB2} and \eqref{08:20}, we have for $\theta\in(0, 2-\frac2{q_b})$,
\begin{align}
\|\sI^\lambda(b\odot\nabla u)\|_{\cS^\lambda_\theta}
\lesssim\|b\odot\nabla u\|_{\mL^{q_b}_T\bB^{\alpha_b}_{p_b,q_b}}\lesssim\kappa_b\|u\|_{\mL^\infty_T\bB^{1}_{\infty,1}}.\label{08:207}
\end{align}
On the other hand,  by \eqref{08:204}, we have for any $\theta'\in(0,2-\frac2{q_b})$,
\begin{align}\label{LA2}
\|\sI^\lambda(\div b\preceq u)\|_{\mL^{q_b}_T\bB^{2-d/p_b-2/q_b-\theta'}_{p_b,1}}
\lesssim(1+\lambda)^{-\frac{\theta'}{2}}\|\div b\preceq u\|_{\mL^{q_b/2}_T\bB^{0}_{p_b/2,\infty}},
\end{align}
and for any $\theta''\in(0,2-\frac4{q_b})$,
\begin{align}\label{LA3}
\|\sI^\lambda(\div b\preceq u)\|_{\mL^{\infty}_T\bB^{2-d/p_b-4/q_b-\theta''}_{p_b,1}}
\lesssim(1+\lambda)^{-\frac{\theta''}{2}}\|\div b\preceq u\|_{\mL^{q_b/2}_T\bB^{0}_{p_b/2,\infty}}.
\end{align}
Since $\alpha_b\in(-1,-\frac12]$ and $\tfrac{d}{p_b}+\tfrac2{q_b}<1+\alpha_b$, we have
\begin{align*}
0<\vartheta_0:=-\alpha_b-\tfrac{d}{p_b}-\tfrac2{q_b}<2-\tfrac2{q_b}.
\end{align*}
So, for fixed $\theta\in(0,2+\alpha_b+\frac{d}{p_b}-\tfrac2{q_b}]$ and $\vartheta\in(0,\vartheta_0)$, it holds that
$$
\vartheta+\theta<2-\tfrac4{q_b}.
$$
If we take $\theta'=\vartheta$ in \eqref{LA2} and $\theta''=\vartheta+\theta$, and $\theta''=\vartheta$ 
in \eqref{LA3}, then by $\alpha_b<-d/q_b-\vartheta$ and \eqref{08:210},
\begin{align}\label{HB1}
\|\sI^\lambda(\div b\preceq u)\|_{\cS^\lambda_\theta}
\lesssim(1+\lambda)^{-\frac{\vartheta}{2}}\|\div b\preceq u\|_{\mL^{q_b/2}_T\bB^{0}_{p_b/2,\infty}}
\lesssim (1+\lambda)^{-\frac{\vartheta}{2}}\tilde \kappa_b\|u\|_{\mL^{q_b}_T\bB^{2+\alpha_b}_{p_b,q_b}},
\end{align}
Moreover, by $\alpha-\tfrac{d}{p}-\tfrac2q\geq \alpha_b-\tfrac{d}{p_b}-\tfrac{2}{q_b}$, \eqref{08:205} and
\eqref{AB2},  for $\theta\in(0, 2-\frac2{q})$,
\begin{align}\label{HB2}
\|\sI^\lambda(f)\|_{\cS^\lambda_\theta}\lesssim \|f\|_{\mL^q_T\bB^{\alpha}_{p,q}}.
\end{align}
 Combining \eqref{08:207}, \eqref{HB1}, \eqref{HB2} and by Remark \ref{KP1}, for $\theta\in(0,2+\alpha_b+\frac{d}{p_b}-\tfrac2{q_b}]$ and $\vartheta\in(0,\vartheta_0)$,
 there are constants
 $C_1=C_1(\Theta,\theta)>0$, $C_2=C_2(\Theta,\theta,\vartheta)>0$ and $C_3=C_3(\Theta,\theta,\alpha,p,q)>0$ such that
\begin{align}\label{08:208}
\|u\|_{\cS^\lambda_\theta}
\leq C_1\kappa_b\|u\|_{\mL^\infty_T\bB^{1}_{\infty,1}}+C_2(1+\lambda)^{-\frac{\vartheta}{2}}\tilde \kappa_b\|u\|_{\mL^{q_b}_T\bB^{2+\alpha_b}_{p_b,q_b}}+C_3\|f\|_{\mL^q_T\bB^{\alpha}_{p,q}}.
\end{align}
Since $\frac2{q_b}+\frac d{p_b}<1+\alpha_b$, by \eqref{Sob1}, 
one can choose $\theta_0=\theta_0(d,\alpha_b,p_b,q_b)$ being close to zero so that
$$
\|u\|_{\mL^\infty_T\bB^{1}_{\infty,1}}\lesssim_C\|u\|_{\mL^\infty_T\bB^{2+\alpha_b-2/q_b-\theta_0}_{p_b,1}}.
$$
Thus by \eqref{08:208} with $\theta=\vartheta=\theta_0$, 
one can take $\lambda_0=\lambda_0(\Theta)>0$ large enough so that for any $\lambda\ge\lambda_0$,
$$
\|u\|_{\cS^\lambda_{\theta_0}}=\|u\|_{\mL^{q_b}_T\bB^{2+\alpha_b}_{p_b,q_b}}
+(1+\lambda)^{\frac{\theta_0}2}\|u\|_{\mL^\infty_T\bB^{2+\alpha_b-2/q_b-\theta_0}_{p_b,1}}\lesssim_C\|f\|_{\mL^q_T\bB^{\alpha}_{p,q}}.
$$
The estimate \eqref{S5:AA07} is then obtained by substituting this into \eqref{08:208}.
\end{proof}

\br\label{Re34}
By \eqref{LA4}, it is easy to see that the unique solution in the above theorem also satisfies that for any $\varphi\in C^\infty_c(\mR^d)$,
$$
\p_t \<u^\lambda,\varphi\>=\<u^\lambda, \Delta \varphi\>-\lambda \<u^\lambda,\varphi\>+\<b\odot\nabla u^\lambda-\div b\preceq u^\lambda, \varphi\>+\<f,\varphi\>.
$$
In other words, $u^\lambda$ is a weak solution of PDE \eqref{S1:PDE}.
\er
\subsection{Supercritical case: $\frac{d}{p}+\frac{2}{q}<2+\alpha$}\label{Sub32}
We introduce the following supercritical index set
\begin{align}\label{XX1}
\bI_d:=\Big\{(\alpha,p,q)\in[-1,0]\times [2,\infty]^2,\ \tfrac{d}{p}+\tfrac{2}{q}<2+\alpha\Big\},
\end{align}
and also the following energy space:
\begin{align}\label{VV1}
\sV_T:=\left\{u: \|u\|_{\sV_T}:=\sup_{t\in[0,T]}\|u(t)\|_2+\|\nabla u\|_{\mL^2_T}<\infty\right\}.
\end{align}
We have the following simple lemma.
\bl\label{Le34}
For any $(\alpha,p,q)\in\bI_d$ and $T>0$, there is a constant $C=C(\alpha,p,q,d, T)>0$ such that for any $f\in\mL^q_T\bH^\alpha_p$ and $u,g\in\sV_T$, 
$$
\|\<fu,g\>\|_{L^1_T}\le C \|f\|_{\mL^q_T\bH^\alpha_p}\|u\|_{\sV_T}\|g\|_{\sV_T}.
$$
\el
\begin{proof}
Let $(\alpha,p,q)\in\bI_d$. Consider the case $p\in[2,\infty)$. Let $r\in(1,2]$ be defined by $\frac1p+\frac1r=1$.
Let us choose $r_0\in(1,r]$ so that $-\alpha-\frac dr\leq1-\frac d{r_0}$. By Sobolev's embedding \eqref{Sob11}, we have
$$
|\<fu,g\>|\lesssim \|f\|_{\bH^\alpha_p}\|ug\|_{\bH^{-\alpha}_{r}}\lesssim \|f\|_{\bH^\alpha_p}\|ug\|_{\bH^1_{r_0}}.
$$
Let $\frac1{r_0}=\frac1{r_1}+\frac12$. By  H\"older's inequality, we have
$$
\|ug\|_{\bH^1_{r_0}}\lesssim\|u\|_{r_1}\|g\|_2+\|u\|_{r_1}\|\nabla g\|_2+\|g\|_{r_1}\|\nabla u\|_2.
$$
Hence, for $\frac1q+\frac1{q_1}+\frac12=1$,
$$
\|\<fu,g\>\|_{L^1_T}\lesssim \|f\|_{\mL^q_T\bH^\alpha_p}\Big(\|u\|_{\mL^{q_1}_TL^{r_1}}\|g\|_{\mL^2_T\bH^1_2}
+\|u\|_{\mL^{2}_T\bH^1_2}\|g\|_{\mL^{q_1}_TL^{r_1}}\Big).
$$
Since $p,q\in[2,\infty)$ satisfies $\frac dp+\frac2q<2+\alpha$, one has $r_1,q_1\in[2,\infty)$ and $\frac d{r_1}+\frac2{q_1}>\frac d2$. By \cite[Lemma 2.1]{ZZ21}, we conclude the desired estimate.

Next for $p=\infty$, by definition we have
$$
|\<fu,g\>|\leq \|f\|_{\bH^\alpha_\infty}\|ug\|_{\bH^{-\alpha}_1}\lesssim \|f\|_{\bH^\alpha_\infty}\|ug\|_{\bH^1_1},
$$
and by H\"older's inequality,
$$
\|ug\|_{\bH^1_1}\lesssim\|u\|_2\|g\|_2+\|u\|_2\|\nabla g\|_2+\|g\|_2\|\nabla u\|_2.
$$
Hence,
$$
\|\<fu,g\>\|_{L^1_T}\lesssim \|f\|_{\mL^2_T\bH^\alpha_\infty}\Big(\|u\|_{\mL^\infty_TL^2}\|g\|_{\mL^2_T\bH^1_2}
+\|u\|_{\mL^{2}_T\bH^1_2}\|g\|_{\mL^\infty_TL^2}\Big).
$$
The proof is complete.
\end{proof}

Next we make the following assumption on $b$.
\begin{itemize}
\item[{\bf(H$^{\rm sup}_b$)}]  Let $(\alpha_b,p_b,q_b)\in\bI_d$. Suppose that
$$
\kappa_b:=\|b\|_{\mL^{q_b}_T\bH^{\alpha_b}_{p_b}}<\infty\quad \text{and}\quad \tilde \kappa_b:=\|\div b\|_{\mL^{2}_TL^\infty}<\infty.
$$
\end{itemize}
Recall the parameter set \eqref{Pa1}. We can show the following existence and uniqueness result.
\bt\label{S1:well}
Under {\bf(H$^{\rm sup}_b$)}, for any $(\alpha,p,q)\in\bI_d$ and $f\in \mL^{q}_T\bH^{\alpha}_{p}$,  
there is a weak solution $u$ to PDE \eqref{S1:PDE} in the sense that
for any smooth function $\varphi\in C^\infty_c(\mR^{d})$,
\begin{align}\label{S2:BB00}
\p_t\<u,\varphi\>=\<u,\Delta \varphi-b\cdot\nabla \varphi+\div b\,\varphi\>+\<f,\varphi\>,
\end{align}
and  for some $C=C(\Theta, \alpha,p,q)>0$,
\begin{align}\label{AF1}
\|u\|_{\mL^\infty_T}+\|u\|_{\sV_T}\le C\| f\|_{\mL^{q}_T\bH^{\alpha}_{p}}.
\end{align}
If in addition that $b\in\mL^\infty_T\cB+\mL^2_TL^2$, then the uniqueness holds in the class \eqref{AF1}.

\et
\begin{proof}
{\bf (Existence)} Let $b_n:=b*\phi_n$ and $f_n:=f*\phi_n$ be the mollifying approximation. Clearly,
$$
\kappa_{b_n}\leq \kappa_{b},\ \ \tilde\kappa_{b_n}\leq \tilde\kappa_{b}.
$$
Since $b_n,f_n\in\mL^\infty_TC^\infty_b$,
it is well-known that there is a classical solution $u_n\in \mL^\infty_TC^\infty_b$ to the following PDE
\begin{align*}
\p_t u_n=\Delta u_n+b_n\cdot\nabla u_n+f_n=\Delta u_n+\div (b_n u_n)-(\div b_n) u_n+f_n,\quad u_n(0)=0. 
\end{align*}
By \cite[Theorem 2.2]{ZZ21}, there is a constant $C=C(\Theta, \alpha,p,q)>0$ such that for all $n\in\mN$,
\begin{align}\label{AA00}
\|u_n\|_{\mL^\infty_T}+\|u_n\|_{\sV_T}\le C\| f_n\|_{\mL^{q}_T\bH^{\alpha}_{p}}\le C\| f\|_{\mL^{q}_T\bH^{\alpha}_{p}}.
\end{align}
By \eqref{Sob11} and \eqref{EW2}, it is easy to see that for $\frac1{p_0}:=\frac1p_b+\frac12\leq\frac1d+\frac12$,
$$
    \|fg\|_{\bH^{-2}_2}\lesssim \|fg\|_{\bH^{-1}_{p_0}}\le \|f\|_{\bH^{-1}_{p_b}}\|g\|_{\bH^1_2}.
$$
Thus by \eqref{AA00}, we have for Lebesgue almost all $t$,
\begin{align*}
\|\p_t u_n(t)\|_{\bH^{-3}_2}
&\le \|\Delta u_n(t)+\div((b_n u_n)(t))-\div b_n(t)\, u_n(t)+f_n(t)\|_{\bH^{-3}_2}\\
&\lesssim \|u_n(t)\|_2 +\|(b_n u_n)(t)\|_{\bH^{-2}_2}+\| \div b_n(t) u_n(t)\|_2+\| f_n(t)\|_{\bH^{\alpha}_p}\\
&\lesssim \|u_n\|_{\mL^\infty_TL^2} +\|b_n(t)\|_{\bH^{-1}_{p_b}}\|u_n(t)\|_{\bH^{1}_2}+\| \div b_n(t)\|_\infty\|u_n(t)\|_2+\| f(t)\|_{\bH^{\alpha}_p}\\
&\lesssim \| f\|_{\mL^{q}_T\bH^{\alpha}_{p}}+\|b(t)\|_{\bH^{-1}_{p_b}}\|u_n(t)\|_{\bH^{1}_2}+\| \div b(t)\|_\infty\| f\|_{\mL^{q}_T\bH^{\alpha}_{p}}+\| f(t)\|_{\bH^{\alpha}_p},
\end{align*}
where the implicit constant does not depend on $n$ and $t$.
Hence, for all $0\leq t_0<t_1\leq T$,
\begin{align*}
&\sup_n\|u_n(t_1)-u_n(t_0)\|_{\bH^{-3}_2}\leq \sup_n\int^{t_1}_{t_0}\|\p_t u_n(t)\|_{\bH^{-3}_2}\dif t\\
&\quad\lesssim (t_1-t_0)+\left(\int^{t_1}_{t_0}\|b(t)\|^2_{\bH^{-1}_{p_b}}\dif t\right)^{1/2}+\int^{t_1}_{t_0}\Big(\|\div b(t)\|_\infty+\| f(t)\|_{\bH^{\alpha}_p}\Big)\dif t,
\end{align*}
which in turn implies that
$$
\lim_{|t_1-t_0|\to 0}\sup_n\|u_n(t_1)-u_n(t_0)\|_{\bH^{-3}_2}=0.
$$
By Aubin-Lions' lemma (see \cite[Corollary 4]{Si86} for instance) and the above uniform continuity, there is a function $u\in \mL^\infty_T\cap \sV_T$ and subsequence $\{n_k\}$ such that
 for each $t$ and $\varphi\in C^\infty_c(\mR^d)$,
\begin{align}\label{SW6}
\lim_{k\to\infty}\<u_{n_k}(t),\varphi\>=\<u(t),\varphi\>,\ \ u_{n_k}\stackrel{k\to\infty}{\to} u\mbox{weakly in $\mL^2_T\bH^1_2$},
\end{align}
and for each $R>0$,
$$
\lim_{k\to\infty}\int^T_0 \|u_{n_k}(s)-u(s)\|^2_{L^2(B_R)}\dif s=0.
$$
Thus, for each $t\in[0,T]$,
$$
\|u(t)\|_{\infty}=\sup_{\varphi\in C^\infty_c, \|\varphi\|_1\leq 1}\<u(t),\varphi\>=
\sup_{\varphi\in C^\infty_c, \|\varphi\|_1\leq 1}\lim_{k\to\infty}\<u_{n_k}(t),\varphi\>\leq \sup_n\|u_n(t)\|_{\infty}\leq C\|f\|_{\mL^{q}_T\bH^{\alpha}_{p}}.
$$
Furthermore, for any $\varphi\in L^1(\mR^d)$, we also have
\begin{align}\label{S2:BB06}
\lim_{k\to\infty}\<u_{n_k}(t),\varphi\>=\<u(t),\varphi\>.
\end{align}
By taking weak limits, it is easy to see that $u$ is a weak solution of PDE \eqref{S1:PDE}. Indeed, 
since $\frac{d}{p_b}+\frac{2}{q_b}<2+\alpha$, without loss of generality we may assume $q_b,p_b<\infty$. Otherwise, we may choose $q_b',p_b'<\infty$ so that
$\frac{d}{p'_b}+\frac{2}{q'_b}<2+\alpha$.
Noting that for any $\varphi\in C^\infty_c(\mR^d)$,
\begin{align*}
\<b_{n_k}u_{n_k}, \nabla \varphi\>-\<bu, \nabla \varphi\>=\<(b_{n_k}-b)u_{n_k}, \nabla \varphi\>+\<u_{n_k}-u, b\cdot\nabla \varphi\>,
\end{align*}
by Lemma \ref{Le34} and $p_b,q_b\in[2,\infty)$, we have
$$
\|\<(b_{n_k}-b)u_{n_k}, \nabla \varphi\>\|_{L^1_T}\le C \|(b_{n_k}-b)\chi\|_{\mL^{q_b}_T\bH^{\alpha_b}_{p_b}}\|u_{n_k}\|_{\sV_T}\|\nabla\varphi\|_{\sV_T}\stackrel{k\to\infty}{\to} 0,
$$
where $\chi\in C^\infty_c$ and $\chi=1$ on the support of $\varphi$,
and due to $b\cdot\nabla \varphi\in\mL^2_T\bH^{-1}_2$, by \eqref{SW6},
$$
\|\<u_{n_k}-u, b\cdot\nabla \varphi\>\|_{L^1_T}\stackrel{k\to\infty}{\to} 0.
$$
Similarly, we have
$$
\|\<\div b_{n_k}u_{n_k}-\div b u,\varphi\>\|_{L^1_T}\stackrel{k\to\infty}{\to} 0.
$$

{\bf (Uniqueness)}  Let $u$ be any weak solution of PDE \eqref{S1:PDE} with $f=0$ and satisfy
$$
\|u\|_{\mL^\infty_T}+\|u\|_{\sV_T}<\infty.
$$ 
Define $u_n:=(u*\phi_n)\chi_n$, where $\chi_n=n^{-d}\chi(x/n)$ and $\chi\in C^\infty_c(\mR^d)$ is a cuttoff function with $\chi(x)=1$ on $B_1$. 
By the chain rule and the definition we have
$$
\p_t\<u, u_n\>/2=\<\p_t u, u_n\>=\<u, \Delta u_n\>-\<u, b\cdot\nabla u_n\>+\<\div b u, u_n\>.
$$
In particular,
$$
\frac{\<u(t), u_n(t)\>}{2}=-\int^t_0\<\nabla u, \nabla u_n\>\dif s-\int^t_0\Big[\<u, b\cdot\nabla u_n\>-\<\div b u, u_n\>\Big]\dif s.
$$
Suppose that $b=b_1+b_2$, where $b_1\in\mL^\infty_T\cB$ and $b_2\in\mL^2_TL^2$.
Noting that
\begin{align}\label{CF1}
|\<u, b_1\cdot\nabla u_n\>|\leq\|u\|_{\bH^1_2}\|b_1\cdot\nabla u_n\|_{\bH^{-1}_2}\leq\|u\|_{\bH^1_2}\|b_1\|_\cB\|u_n\|_{\bH^1_2}\leq C\|b_1\|_\cB\|u\|^2_{\bH^1_2},
\end{align}
and
$$
|\<u, b_2\cdot\nabla u_n\>|\leq\|u\|_\infty\|b_2\|_2\|\nabla u_n\|_2\leq\|u\|_\infty\|b_2\|_2\|\nabla u\|_2,
$$
by the dominated convergence theorem, we get
$$
\frac{\|u(t)\|^2_2}{2}\leq-\int^t_0\|\nabla u\|^2_2\dif s+\int^t_0\lim_{n\to\infty}\Big|\<u, b\cdot\nabla u_n\>-\<\div b u, u_n\>\Big|\dif s.
$$
By the definition of $\cB$, it is easy to see that
$$
\lim_{n\to\infty}b_1\cdot\nabla u_n\mbox{exists in $\bH^{-1}_2$},\ \ \lim_{n\to\infty}b_2\cdot\nabla u_n=b_2\cdot\nabla u\mbox{in $L^1$}.
$$
Hence, by $\div b\in L^\infty$ and commutating the order of limits (guaranteed by the above two limits),
\begin{align*}
\lim_{n\to\infty}\<u, b\cdot\nabla u_n\>=\lim_{n\to\infty}\lim_{m\to\infty}\<u_m, b\cdot\nabla u_n\>
&=-\lim_{n\to\infty}\lim_{m\to\infty}\<\div b\, u_m, u_n\>-\lim_{n\to\infty}\lim_{m\to\infty}\<b\cdot \nabla u_m, u_n\>\\
&=-\<\div b\, u, u\>-\lim_{m\to\infty}\<b\cdot \nabla u_m, u\>,
\end{align*}
which implies that
$$
\lim_{n\to\infty}\<u, b\cdot\nabla u_n\>=-\<\div b\, u, u\>/2.
$$
Thus,
$$
\frac{\|u(t)\|^2_2}{2}\leq-\int^t_0\|\nabla u\|^2_2\dif s+\frac12\int^t_0|\<\div b\, u, u\>|\dif s\leq\frac12\int^t_0\|\div b\|_\infty \|u\|_2^2\dif s.
$$
Now by Gronwall's inequality, we obtain $u(t)=0$. The uniqueness is proven.
\end{proof}

\section{Subcritical case: Proof of Theorem \ref{Th1}}

In this section, we establish the weak well-posedness of SDE \eqref{in:SDE} under the subcritical assumption.
To achieve this, we employ the classical Zvonkin transformation (see \cite[Section 5]{ZZ17}), ensuring both the existence and uniqueness of a weak solution.

Throughout this section we fix $T>0$. To introduce a general notion of Krylov's estimate, we start by fixing a Banach space $\mB$ that consists of distributions $f_t: [0,T]\to\sS'(\mR^d)$ so that
$\mL^1_TC_b\cap\mB$ is dense in $\mB$. Moreover, we make the following assumption:
\begin{itemize}
\item
For any bounded measurable function $g:[0,T]\times\mR^d\to\mR$ such that the function $x \to g(t, x)$ is uniformly Lipschitz with respect to $t \in [0, T]$, the product $gf$ is well-defined for all $f \in \mB$,
and there is a constant $C_g > 0$ such that for all $f\in\mB$,
\begin{align}\label{S4:AA00}
\|gf\|_{\mB}\le C_g\|f\|_{\mB}.
\end{align}
\end{itemize}
It's worth noting that the space $\mB:=\mL^q_T\bB^{\alpha}_{p,q}$ with $\alpha\in[-1,0]$ and $p,q\in[1,\infty)$ satisfies the above assumption.
Indeed, \eqref{S4:AA00} follows by Lemma \ref{Le24}. Moreover, let $\phi_n(x)=n^d\phi(nx)$, where  $\phi\in C^\infty_c(\mR^d)$ is a smooth probability density function. For $f\in\mL^q_T\bB^{\alpha}_{p,q}$, if we define
$$
f_n(t,x)=f(t,\cdot)*\phi_n(x),
$$
then by definition, $f_n\in\mL^q_TC^\infty_b$, and by the dominated convergence theorem,
\begin{align}\label{re2}
    \lim_{n\to\infty}\|f_n-f\|^q_{\mL^q_T\bB^{\alpha}_{p,q}}\leq\int^T_0\sum_{j\geq-1}2^{\alpha j q}\lim_{n\to\infty}\|\phi_n*\cR_jf(t)-\cR_j f(t)\|_p^q\dif t=0.
\end{align}
{Throughout most of this section, we limit our analysis to the case $p,q<\infty$ when considering the space $\mB=\mL^q_T\bB^{\alpha}_{p,q}$. The case $p=\infty$ or $q=\infty$ will be addressed later in the proof of Theorem \ref{Th1}.
}

 Now we introduce the following notion used below.
\bd[Krylov's estimate] \label{Kry}
Let $(X_t)_{t\in[0,T]}$ be a $\mR^d$-valued stochastic process on a probability space $(\Omega,\sF,\bP)$. 
For $\theta\in(0,1)$ and $m>2$, 
one says that $(X,\mB)$ satisfies Krylov's estimate with parameters $(m,\theta)$ and constant $C>0$
if for any $f\in \mL^1_TC_b\cap\mB$ and any stopping times 
$0\leq \tau_0\leq \tau_1\leq T$ with $\tau_1-\tau_0\le\delta$,
\begin{align}\label{Kry0}
    \left\|\int^{\tau_1}_{\tau_0}f(s,X_s)\dif s\right\|_{L^m(\Omega)}\leq C\delta^{\frac{1+\theta}2}\|f\|_\mB.
\end{align}
\ed
\br
Here we require $\theta>0$, which is crucial for Young's integral in \eqref{AA0} below.
\er
The following proposition is a direct conclusion of Krylov's estimate and the Kolmogorov continuity theorem. 
\bp\label{Pr88}
Suppose that $(X,\mB)$ satisfies Krylov's estimate with parameter $(m,\theta)\in(2,\infty)\times(0,1)$.
For any $f\in \mB$, let $f_n\in \mL^1_TC_b\cap\mB$ be such that
$$
\lim_{n\to\infty}\|f_n-f\|_{\mB}=0.
$$
Then the limit $A^f_\cdot:=\lim_{n\to\infty}\int^\cdot_0f_n(s,X_s)\dif s$ exists in $L^m(\Omega,C([0,T]))$  and does not depend on the choice of approximation sequence.
Furthermore,  
\begin{align}\label{S4:AA01}
\sup_{t\in[0,T]}\|A^f_t\|_{L^m(\Omega)}+\sup_{s\ne t\in[0,T]}\frac{\|A^f_t-A^f_s\|_{L^m(\Omega)}}{|t-s|^{(1+\theta)/2}}\le C\|f\|_\mB.
\end{align}
\ep

Next we show a substitution formula for Young's integrals that will be used to show the uniqueness by Zvonkin's transformation.
\bp\label{Pro44}
Let $m>2$ and $\theta\in(0,1)$ satisfy $\theta m>2$. Suppose that $(X,\mB)$ satisfies Krylov's estimate with parameters $(m,\theta)$, and for any $\beta\in(0,\frac12)$,
\begin{align}\label{HF1}
\left\|\sup_{s\not=t\in[0,T]}|X_t-X_s|/|t-s|^\beta\right\|_{L^m(\Omega)}<\infty.
\end{align}
Let $g:\mR_+\times\mR^d\to\mR$ be a bounded function satisfying
\begin{align}\label{HF2}
|g(t,x)-g(s,y)|\le C\sqrt{|t-s|}+C|x-y|. 
\end{align}
Then for any space-time function $f\in\mB$,  the integral $\int^t_0g(s, X_s)\dif A^f_s$ is well-defined as Young's integral 
and
\begin{align}\label{AA0}
\int^t_0g(s, X_s)\dif A^f_s=A^{f\cdot g}_t\ \ a.s.,
\end{align}
where $A^f$ is defined in Proposition \ref{Pr88}.
\ep
\begin{proof}
For a function $h:[0,T]\to\mR$ and $\alpha\in(0,1)$, we write
\begin{align*}
[h]_{\alpha}:=\sup_{s\ne t\in[0,T]}\frac{|h(t)-h(s)|}{|t-s|^\alpha}.
\end{align*}
Since $\theta m>2$, one can choose $\alpha\in(\frac12,\frac{1+\theta}2-\frac1m)$ and  
$\beta\in(0,\frac12)$ so that
$$
\alpha+\beta>1.
$$
By \eqref{S4:AA01}, \eqref{HF1} and Kolmogorov's continuity criterion,
there is a constant $C>0$ such that for all $f\in\mB$,
\begin{align}\label{HF5}
\big\|[A^f_\cdot]_{\alpha}\big\|_{L^m(\Omega)}\le C\|f\|_\mB,\ \ \big\|[X_\cdot]_{\beta}\big\|_{L^m(\Omega)}<\infty.
\end{align}
Below we fix  a sample point $\omega$ such that
\begin{align}\label{HF3}
[A^f_\cdot(\omega)]_{\alpha}<\infty,\ \ [X_\cdot(\omega)]_{\beta}<\infty.
\end{align}
For any $0\le s<t\le T$, we define
\begin{align*}
\Gamma_{s,t}(\omega):=g(s,X_s(\omega))(A^f_t(\omega)-A^f_s(\omega)).
\end{align*}
For simplicity, we drop the dependence of $\omega$.
We note that for any $0\le s<u<t\le T$,
\begin{align*}
\delta \Gamma_{s,u,t}:=\Gamma_{s,t}-\Gamma_{s,u}-\Gamma_{u,t}=(g(s,X_s)-g(u,X_u))(A^f_t-A^f_u).
\end{align*}
By \eqref{HF2} and \eqref{HF3}, we have
 \begin{align*}
|\delta \Gamma_{s,u,t}|\lesssim |u-s|^\beta|t-u|^{\alpha}[X_\cdot]_{\beta}[A^f_\cdot]_{\alpha}.
\end{align*}
 Since $\alpha+\beta>1$, by the sewing lemma (cf. \cite[Lemma 4.2]{FH14}), we have
 \begin{align*}
\Gamma_t:=\int_0^t g(r,X_r)\dif A^f_r:=\lim_{|\pi|\to0}\sum_{r,s\in\pi}\Gamma_{r,s}\quad \text{exists},
\end{align*}
 where $\pi$ is any partition of $[0,t]$, and
 \begin{align*}
\left|\int_s^t(g(r,X_r)-g(s,X_s))\dif A^f_r\right|=|\Gamma_t-\Gamma_s-\Gamma_{s,t}|\lesssim |t-s|^{\alpha+\beta}[X_\cdot]_{\beta}[A^f_\cdot]_{\alpha}.
\end{align*}
In particular, taking $s=0$, we get
 \begin{align}\label{S4:AA02}
\left|\int_0^tg(r,X_r)\dif A^f_r\right|&\leq \|g\|_\infty|A^f_t|+C [X_\cdot]_{\beta}[A^f_\cdot]_{\alpha}.
\end{align}
 Now, let $\{f_n\}_{n=1}^\infty \subset \mL^1_TC_b\cap\mB$ be such that $\lim_{n\to\infty}\|f_n-f\|_{\mB}=0$. By definition, one sees that
 \begin{align}\label{S4:AA03}
\int_0^t g(r,X_r)\dif A^{f_n}_r=\int_0^t (gf_n)(r,X_r) \dif r=A^{gf_n}_t\quad a.s.
\end{align}
By \eqref{S4:AA02}, \eqref{HF5} and \eqref{S4:AA00}, we have
  \begin{align*}
&\bE\left(\left| \int_0^t g(r,X_r)\dif (A^{f_n}_r-A^{f}_r)\right|\right)+\bE|A^{gf_n}_t-A^{gf}_t|\\
&\quad=\bE\left(\left| \int_0^t g(r,X_r)\dif A^{f_n-f}_r\right|\right)+\bE|A^{g(f_n-f)}_t|\\
&\quad\lesssim \|f-f_n\|_{\mB}+\|g(f-f_n)\|_{\mB}\to0 \quad\text{as $n\to\infty$}.
\end{align*}
Thus, there is a subsequence $\{n_k\}_{k=1}^\infty$ such that
  \begin{align*}
\lim_{k\to\infty}\left| \int_0^t g(r,X_r)\dif (A^{f_{n_k}}_r-A^{f}_r)\right|+|A^{gf_{n_k}}_t-A^{gf}_t|=0\quad a.s.
\end{align*}
  By taking $n=n_k$ and letting $k\to\infty$ in \eqref{S4:AA03}, we get \eqref{AA0}.
\end{proof}

Below, we always suppose that $b=b_1+b_2$ satisfies
\begin{itemize}
    \item[{($\widetilde{\bf H}^{\rm sub}_{b}$)}]  $b_1$ satisfies {\bf(H$^{\rm sub}_{b_1}$)} with 
    $p_b,q_b<\infty$, 
    and $|b_2(t,x)|\le c_0+c_1|x|$ for some $c_0,c_1\geq0$.
\end{itemize}

By Theorem \ref{S5:thm01}, there is a  unique solution $u^\lambda: [0,T]\times\mR^d\to\mR^d$  to the following $\mR^d$-valued backward PDE system
$$
\p_tu^\lambda+\Delta u-\lambda u^\lambda+{b_1}\odot\nabla u^\lambda+\div {b_1}\preceq u^\lambda+{b_1}=0,\ \ u(T)=0,
$$
with the regularity that there is a $\lambda_0>0$ such that for all $\lambda\geq\lambda_0$ and $\theta\in(0,2+\alpha_b+\frac d{p_b}-\frac{2}{q_b}]$,
\begin{align}\label{Ho2}
\|u^\lambda\|_{\mL^\infty_T\bB^{2+\alpha_b-2/q_b-\theta}_{p_b,1}}\lesssim
(1+\lambda)^{-\frac{\theta}2}\|{b_1}\|_{\mL^{q_b}_T\bB^{\alpha_b}_{p_b,q_b}},
\end{align}
where the implicit constant does not depend on $\lambda$.
In particular, by embedding \eqref{Sob1}, one can choose $\lambda$ large enough so that
\begin{align}\label{Ho22}
\sup_{t\in[0,T]}\|\nabla u^\lambda(t)\|_\infty\leq 1/2.
\end{align}
Let 
\begin{align}\label{Phi1}
\Phi(t,x):=x+u^\lambda(t,x).
\end{align}
Then for each $t\in[0,T]$,
$$
x\mapsto \Phi(t,x)\mbox{forms a $C^1$-diffeomorphism.}
$$
In particular,
\begin{align}\label{Ho12}
\|\nabla\Phi\|_\infty,\ \ \|\nabla\Phi^{-1}\|_\infty\leq 4,
\end{align}
where $\Phi^{-1}$ stands for the inverse of $x\mapsto \Phi(t,x)$.

Below, for simplicity of notations, we write
$$
\bI^{\alpha_b}_{p_b,q_b}:=\Big\{(\alpha,p,q)\in[-1,0]\times[1,\infty)^2: \tfrac{q_b}2\leq q\leq q_b,\ \ p\leq p_b,\ \ \alpha-\tfrac dp-\tfrac2q\geq\alpha_b-\tfrac d{p_b}-\tfrac2{q_b}\Big\}.
$$
We now proceed to show the following crucial Zvonkin transformation.
\bl[Zvonkin's transformation]\label{Th45}
Let $m>2$ and $\theta\in(0,1)$ satisfy $\theta m>2$. 
Let $(\mathfrak{F}, X,W)$ be a weak solution of SDE \eqref{in:SDE} in the sense of Definition \ref{Def1}. 
Suppose that for any $(\alpha,p,q)\in\bI^{\alpha_b}_{p_b,q_b}$, 
$(X, \mL^q_T\bB^\alpha_{p,q})$ satisfies Krylov's estimate with parameter $(m,\theta)$.
Then $Y_t:=\Phi(t,X_t)$ solves the following SDE:
\begin{align}\label{Zo1}
\dif Y_t=\tilde b(t,Y_t)\dif t+\tilde \sigma(t,Y_t)\dif B_t,
\end{align}
where
$$
\tilde b(t,y):=\lambda u(t,\Phi^{-1}(t,y))+(b_2\cdot\nabla\Phi)(t,\Phi^{-1}(t,y)),\ \ \tilde\sigma(t,y):=\sqrt2\nabla\Phi(t, \Phi^{-1}(t,y)).
$$
Moreover, $\tilde b$ and $\tilde\sigma$ are bounded measurable and for some $\gamma\in(0,1)$ and all $t\in[0,T]$, $y,y'\in\mR^d$,
\begin{align}\label{XZ1}
{|\tilde b(t,y)|\leq C(\lambda+c_0+c_1|y|)},\ \ |\tilde \sigma(t,y)-\tilde \sigma(t,y')|\leq C|y-y'|^\gamma,
\end{align}
where $C=C(T,\gamma,\Theta)$,
and for all $\xi\in\mR^d$,
\begin{align}\label{XZ2}
 |\xi|/8\leq |\tilde\sigma(t,y)\xi|\leq 8|\xi|.
\end{align}
\el
\begin{proof}
Let $\phi_n(x)$ be the mollifier as in \eqref{BN1}. Define
$$
u_n(t,x):=u(t,\cdot)*\phi_n(x).
$$
Then by Remark \ref{Re34} with $\varphi=\phi_n(x-\cdot)$, we get
\begin{align}\label{GG1}
\p_s u_n+\Delta u_n-\lambda u_n=-\phi_n*( {b_1}\odot\nabla u -\div {b_1}\preceq u+b)=:g_n,
\end{align}
and
$$
\sup_{t\in[0,T]}\|\nabla u_n(t)\|_\infty\leq \sup_{t\in[0,T]}\|\nabla u(t)\|_\infty\leq 1/2.
$$
Let 
$$
\Phi_n(t,x):=x+u_n(t,x).
$$
Then 
$$
\|\nabla\Phi_n\|_\infty,\ \ \|\nabla\Phi^{-1}_n\|_\infty\leq 4,
$$
and by \eqref{GG1}, 
$$
g_n\in \mL^1_TL^\infty\Rightarrow \p_t\Phi_n\in \mL^1_TL^\infty.
$$
By Krylov's estimate and definition, it is easy to see that
$$
X_t=X_0+A^{b}_t+\sqrt{2}W_t
=X_0+A^{b_1}_t+\int^t_0b_2(s, X_s)\dif s+\sqrt{2}W_t,
$$
where $A^{b_1}_t:=\lim_{n\to\infty}\int^t_0b_{1,n}(s,X_s)\dif s$ by Proposition \ref{Pr88},  satisfies 
$$
\sup_{t\in[0,T]}\|A^{b_1}_t\|_{L^m(\Omega)}+\sup_{s\ne t\in[0,T]}\frac{\|A^{b_1}_t-A^{b_1}_s\|_{L^m(\Omega)}}{|t-s|^{(1+\theta_0)/2}}\le C\|{b_1}\|_{\mL^{q_b}_T\bB^{\alpha_b}_{p_b,q_b}}.
$$
Now by generalized It\^o's formula (see \cite{Fo}),  \eqref{GG1} and \eqref{AA0}, we have
\begin{align}
Y^n_t&:=\Phi_n(t,X_t)=\Phi_n(0,x)+\int^t_0(\p_su_n+\Delta u_n)(s,X_s)\dif s+\int_0^t (b_2\cdot\nabla u_n)(s,X_s)\dif s\no\\
&\qquad+\int^t_0\dif A^{b_1}_s\cdot\nabla\Phi_n(s,X_s)+\sqrt 2\int^t_0 \nabla \Phi_n(s,X_s)\dif W_s\no\\
    &=\Phi_n(0,x)+\lambda\int^t_0u_n(s,X_s)\dif s+A^{b_1\cdot\nabla\Phi_n+g_n}_t+\int_0^t (b_2\cdot\nabla u_n)(s,X_s)\dif s\no\\
&\qquad+\sqrt 2\int^t_0 \nabla \Phi_n(s,X_s)\dif W_s.\label{Ho3}
\end{align}
By \eqref{08:21}, we make the following decomposition:
\begin{align*}
{b_1}\cdot\nabla\Phi_n+g_n&={b_1}+{b_1}\cdot\nabla u_n-\phi_n*( {b_1}\odot\nabla u -\div {b_1}\preceq u+{b_1})=h^n_1+h^n_2,
\end{align*}
where
$$
h^n_1:={b_1}-{b_1}*\phi_n+{b_1}\odot\nabla u_n-\phi_n*({b_1}\odot\nabla u),\ \ h^n_2:=\div {b_1}\preceq u_n-\phi_n*(\div {b_1}b\preceq u).
$$
In particular,
$$
A^{{b_1}\cdot\nabla\Phi_n+g_n}_t=A^{h^n_1}_t+A^{h^n_2}_t.
$$
By Krylov's estimate, we clearly have
\begin{align*}
\bE|A^{h^n_1}_t|\lesssim\|h^n_1\|_{\mL^{q_b}_T\bB^{\alpha_b}_{p_b,q_b}}
&\lesssim\|{b_1}-{b_1}*\phi_n\|_{\mL^{q_b}_T\bB^{\alpha_b}_{p_b,q_b}}+\|{b_1}\odot\nabla (u_n-u)\|_{\mL^{q_b}_T\bB^{\alpha_b}_{p_b,q_b}}\\
&\quad+\|\phi_n*({b_1}\odot\nabla u)-{b_1}\odot\nabla u\|_{\mL^{q_b}_T\bB^{\alpha_b}_{p_b,q_b}}.
\end{align*}
Note that by \eqref{08:20}, \eqref{AB2}, \eqref{Ho1} and \eqref{Ho2}, there is an $\eps>0$ so that for all $n\in\mN$,
\begin{align*}
\|{b_1}\odot\nabla (u_n-u)\|_{\mL^{q_b}_T\bB^{\alpha_b}_{p_b,q_b}}
&\lesssim \kappa_b\|u_n-u\|_{\mL^\infty_T\bB^{1+d/p_b}_{p_b,1}}
\lesssim n^{-\eps}.
\end{align*}
Thus, by the property of convolution, we have
$$
\lim_{n\to\infty}\bE|A^{h^n_1}_t|=0.
$$
Moreover, since $\frac d{p_b}+\frac2{q_b}<1+\alpha_b$, one can choose $\eps$ with
$$
0\leq-1-2\alpha_b<\eps\leq-\big(\alpha_b+\tfrac d{p_b}+\tfrac2{q_b}\big)
$$
so that $(-\eps,p_b/2,q_b/2)\in\bI^{\alpha_b}_{p_b,q_b}$. By Krylov's estimate again, we have
\begin{align*}
\bE|A^{h^n_2}_t|\lesssim\|h^n_2\|_{\mL^{q_b/2}_T\bB^{-\eps}_{p_b,q_b}}
&\lesssim\|\phi_n*(\div {b_1}\preceq u)-\div {b_1}\preceq u\|_{\mL^{q_b/2}_T\bB^{-\eps}_{p_b/2,q_b}}\\
&+\|\div {b_1}\preceq(u_n-u)\|_{\mL^{q_b/2}_T\bB^{-\eps}_{p_b,q_b}}=:I_n+J_n.
\end{align*}
For $J_n$, by \eqref{08:210} we have
$$
\lim_{n\to\infty}J_n\lesssim\lim_{n\to\infty}\|\div {b_1}\preceq(u_n-u)\|_{\mL^{q_b/2}_T\bB^{0}_{p_b,\infty}}
\lesssim \kappa^{b_1}_2\lim_{n\to\infty}\|u_n-u\|_{\mL^{q_b}_T\bB^{2+\alpha_b}_{p_b,q_b}}=0.
$$
For $I_n$, since $\|\div {b_1}\preceq u\|_{\mL^{q_b/2}_T\bB^{-\eps}_{p_b/2,q_b}}\lesssim \|\div {b_1}\preceq u\|_{\mL^{q_b/2}_T\bB^{0}_{p_b/2,\infty}}<\infty$, by the property of convolution, we also have
$$
\lim_{n\to\infty}I_n\to 0.
$$
Combining the above limits, we obtain
$$
\lim_{n\to\infty}\bE|A^{{b_1}\cdot\nabla\Phi_n+g_n}_t|=0.
$$
By taking limits for both sides of \eqref{Ho3}, we conclude the proof of \eqref{Zo1}.
As for \eqref{XZ1} and \eqref{XZ2}, it follows by \eqref{Ho22} and \eqref{Ho12}.
\end{proof}

To show the existence of a weak, 
we consider the following approximation SDE:
\begin{align}\label{NC22}
\dif X^n_t=b_n(t,X^n_t)\dif t+\sqrt2\dif W_t,\ \  X^n_0\sim\mu\in\cP(\mR^d),
\end{align}
where  $b_n:=b_{1,n}+b_2$ and $b_{1,n}(t,x):=b_1(t,\cdot)*\phi_n(x)$.
{Since $\int^T_0\|b_{1,n}(t,\cdot)\|^{q_b}_\infty\dif t<\infty$ and $|b_2(t,x)|\leq c_0+c_1|x|$,
it is well known that the weak well-posedness holds for SDE \eqref{NC22} (see \cite{ZZ17}).}
Note that by the property of convolution,
\begin{align}\label{NC21}
\begin{split}
\kappa_{b_{1,n}}:=\|b_{1,n}\|_{\mL^{q_b}_T\bB^{\alpha_b}_{p_b,q_b}}&\leq \|b_1\|_{\mL^{q_b}_T\bB^{\alpha_b}_{p_b,q_b}}=\kappa_{b_1},\\
\tilde\kappa_{b_{1,n}}:=\|\div b_{1,n}\|_{\mL^{q_b}_T\bB^{-2-\alpha_b}_{p_b,q_b/(q_b-1)}}&\leq \|\div b_1\|_{\mL^{q_b}_T\bB^{-2-\alpha_b}_{p_b,q_b/(q_b-1)}}=\tilde\kappa_{b_1}.
\end{split}
\end{align}
Now we can show the following crucial Krylov estimate.
\bl\label{Le42}
Let $\mu\in\cP(\mR^d)$ have finite $m$-order moment, where $m$ satisfies \eqref{MT}.
For any $(\alpha,p,q)\in\bI^{\alpha_b}_{p_b,q_b}$,  $(X^n, \mL^q_T\bB^\alpha_{p,q})$ satisfies Krylov's estimate \eqref{Kry0} 
with parameters $(m,\theta_0)$ 
and constant $C=C(\Theta,\alpha,p,q,m,\mu, c_0,c_1)$, where $c_0,c_1$ are from {\bf ($\widetilde{\bf H}^{\rm sub}_{b}$)}, and $C=C(\Theta,\alpha,p,q, c_0)$ when $c_1=0$.
\el
\begin{proof}
We divide the proof into two steps.

{\bf (Step 1).} We establish the following uniform moment estimate: for some $C_m = C_m(\Theta, c_0, c_1) > 0$,
\begin{align}\label{1213:00}
    \sup_n \left\|\sup_{t \in [0,T]} |X^n_t|\right\|_{L^m(\Omega)} \leq C_m \big(1 + (\mu(|\cdot|^m))^{1/m}\big) =: \ell_m.
\end{align}
Let $\Phi_n$ be defined as in \eqref{Phi1}, where $b_{1,n}$ is used in place of $b_1$. Define $Y^n_t := \Phi_n(t, X^n_t)$. Using SDE \eqref{Zo1}, BDG's inequality, \eqref{XZ1}, and \eqref{XZ2}, we obtain:
\begin{align*}
    \left\|\sup_{t \in [0,T]} |Y^n_t|\right\|_{L^m(\Omega)} \lesssim \|\Phi_n(0, X^n_0)\|_{L^m(\Omega)} + 1 + \int_0^t \|Y^n_s\|_{L^m(\Omega)} \dif s.
\end{align*}
By Gronwall's inequality and the fact that $|X^n_t| = |\Phi_n^{-1}(t, Y^n_t)| \lesssim 1 + |Y^n_t|$, this implies \eqref{1213:00}.

{\bf (Step 2).} Given  $f\in C^\infty_c([0,T]\times\mR^d)$, let $u_n\in\mL^\infty_TC^\infty_b$ be the unique smooth solution to the following backward PDE
$$
\p_tu_n+\Delta u_n-\lambda u_n+b_{1,n}\cdot\nabla u_n+f=0,\ \ u_n(T)=0.
$$
For any $(\alpha,p,q)\in\bI^{\alpha_b}_{p_b,q_b}$, by \eqref{S5:AA07} and \eqref{NC21}, there is a $\lambda_0\geq 1$ so that for all $\theta\in(0,2+\alpha_b+\frac d{p_b}-\frac{2}{q_b}]$ and $\lambda\geq\lambda_0$,
\begin{align}\label{NC1}
\sup_n\|u_n\|_{\mL^\infty_T\bB^{2+\alpha_b-2/q_b-\theta}_{p_b,1}}\lesssim
(1+\lambda)^{-\frac{\theta}2}\|f\|_{\mL^q_T\bB^{\alpha}_{p,q}}.
\end{align}
By It\^o's formula, we have
$$
u_n(t,X^n_t)=u_n(0,x)+\int^t_0 (\lambda u_n-f)(s,X^n_s)\dif s+{\int_0^t (b_2\cdot\nabla u_n)(s,X^n_s)\dif s}+\sqrt 2\int^t_0 \nabla u_n(s,X^n_s)\dif W_s.
$$
In particular, for any stopping time $\tau_0\leq \tau_1\leq T$ with $\tau_1-\tau_0\leq\delta$,
\begin{align*}
    \left|\int^{\tau_1}_{\tau_0}f(s,X^n_s)\dif s\right|\leq 2\|u_n\|_\infty+\lambda\delta\|u_n\|_\infty+{\|\nabla u_n\|_\infty\!\!\int_{\tau_0}^{\tau_1}|b_2|(s,X^n_s)\dif s}+
\sqrt 2\left|\int^{\tau_1}_{\tau_0} \nabla u_n(s,X^n_s)\dif W_s\right|.
\end{align*}
By BDG's inequality and estimate \eqref{1213:00}, we have 
\begin{align*}
\left\|\int^{\tau_1}_{\tau_0}f(s,X^n_s)\dif s\right\|_{L^m(\Omega)}
&\lesssim(1+\lambda\delta)\|u_n\|_\infty +\|\nabla u_n\|_\infty\delta \left\|\sup_{t\in[0,T]}(c_0+c_1|X^n_t|)\right\|_{L^m(\Omega)}\\
&+m\left\|\left(\int^{\tau_1}_{\tau_0}| \nabla u_n(s,X^n_s)|^2\dif s\right)^{1/2}\right\|_{L^m(\Omega)}\\
\lesssim&(1+\lambda\delta)\|u_n\|_\infty+(m\delta^{1/2}+ \delta(c_0+c_1\ell_m))\|\nabla u_n\|_\infty.
\end{align*}
Here and below, the implicit constant is independent of $m,n,\delta,\lambda$.
Since 
$$
\theta_0=1+\alpha_b-\tfrac{d}{p_b}-\tfrac2{q_b}\leq 2+\alpha_b+\tfrac d{p_b}-\tfrac{2}{q_b},
$$
by \eqref{NC1} and  embedding \eqref{Sob1}, we have
$$
\|u_n\|_\infty\lesssim \|u_n\|_{\mL^\infty_T\bB^{d/p_b}_{p_b,1}}=\|u_n\|_{\mL^\infty_T\bB^{2+\alpha_b-2/q_b-(1+\theta_0)}_{p_b,1}}\lesssim (1+\lambda)^{-\frac{1+\theta_0}2}\|f\|_{\mL^q_T\bB^{\alpha}_{p,q}},
$$
and
$$
\|\nabla u_n\|_\infty\lesssim \|u_n\|_{\mL^\infty_T\bB^{1+d/p_b}_{p_b,1}}=\|u_n\|_{\mL^\infty_T\bB^{2+\alpha_b-2/q_b-\theta_0}_{p_b,1}}\lesssim (1+\lambda)^{-\frac{\theta_0}2}\|f\|_{\mL^q_T\bB^{\alpha}_{p,q}}.
$$
 Thus
\begin{align*}
\left\|\int^{\tau_1}_{\tau_0}f(s,X^n_s)\dif s\right\|_{L^m(\Omega)}
\lesssim&\left[(1+\lambda\delta)(1+\lambda)^{-\frac{1+\theta_0}2}+m\sqrt{\delta}(1+\lambda)^{-\frac{\theta_0}2}\right]
\|f\|_{\mL^q_T\bB^{\alpha}_{p,q}}\\
&+ \delta(c_0+c_1\ell_m)\|f\|_{\mL^q_T\bB^{\alpha}_{p,q}}
\end{align*}
In particular, if we take $\lambda=\lambda_0\vee\delta^{-1}$, then for all $\delta\in(0,1)$,
\begin{align}\label{XX2}
\left\|\int^{\tau_1}_{\tau_0}f(s,X^n_s)\dif s\right\|_{L^m(\Omega)}\lesssim (m+c_0+c_1\ell_m)\delta^{\frac{1+\theta_0}2}
\|f\|_{\mL^q_T\bB^{\alpha}_{p,q}}.
\end{align}
By a standard approximation, the above estimate still holds for any $f\in\mL^1_TC_b\cap \mL^q_T\bB^{\alpha}_{p,q}$. 
\end{proof}
\br
By Stirling's formula,  when $c_1=0$, estimate \eqref{XX2} also implies that for some $c>0$,
$$
\sup_n\bE\exp\left\{c\left|\int^T_0f(s,X^n_s)\dif s\right|/\|f\|_{\mL^q_T\bB^{\alpha}_{p,q}}\right\}<\infty.
$$
\er

Now we can show the following main result of this section.
\bt\label{23thm:48}
Suppose that {\bf ($\widetilde{\bf H}^{\rm sub}_{b}$)} holds.
\begin{enumerate}[(i)]
\item For any $\mu \in \cP(\mR^d)$ having finite $m$-order moment,  
where $m$ satisfies \eqref{MT}, there exists a weak solution $(\mathfrak{F}, X, W)$ to the SDE \eqref{in:SDE} with initial distribution $\mu$. 
Uniqueness holds within the class such that for any $(\alpha, p, q) \in \bI^{\alpha_b}_{p_b, q_b}$, 
$(X, \mL^q_T \bB^\alpha_{p,q})$ satisfies Krylov's estimate with parameters $(m, \theta_0)$, where $\theta_0$ is as defined in \eqref{MT}.
\item If $c_1 = 0$ {in {\bf ($\widetilde{\bf H}^{\rm sub}_{b}$)}, then} the moment assumption on $\mu$ can be dropped; 
and for any $t \in (0, T]$ and $X_0 = x$, $X_t$ admits a density $\rho_t(x, x')$ 
that satisfies the following two-sided Gaussian estimate:
for all $x, x' \in \mR^d$,
\begin{align}\label{HH1}
    C_0 t^{-d/2} \e^{-\gamma_0 |x-x'|^2/t} \leq \rho_t(x, x') \leq C_1 t^{-d/2} \e^{-\gamma_1 |x-x'|^2/t},
\end{align}
where $C_0, C_1, \gamma_0, \gamma_1 > 0$ depend only on the parameters $T, d, \alpha_b, p_b, q_b$, 
and $\|b\|_{\mL^{q_b}_T \bB^{\alpha_b}_{p_b, q_b}}$.
\end{enumerate}
\et
\begin{proof} 
{\bf (Existence).} Consider the approximation SDE \eqref{NC22}. By Lemma \ref{Le42}, the law of $X^n$ in $\mC_T$ is tight. 
Therefore, by Prohorov's theorem, it is relatively compact in $\cP(\mC_T)$. Using the weak convergence method, 
it is standard to derive the existence of a weak solution (see \cite{ZZ18} and also the proof of Theorem \ref{Th5} below).

{\bf (Uniqueness and heat kernel estimate).} By Lemma \ref{Th45}, since the transformed SDE \eqref{Zo1} has 
a linear growth drift and bounded, uniformly non-degenerate H\"older continuous diffusion coefficients, 
it is well known that SDE \eqref{Zo1} admits a unique weak solution $Y_t$ (see \cite{ZZ17}). Consequently, weak uniqueness 
also holds for the original SDE \eqref{in:SDE}. 

When $c_1 = 0$, the Krylov constant in \eqref{Kry0} does not depend on the initial distribution, allowing us to drop the moment assumption. 
Moreover, in this case, the transformed SDE \eqref{Zo1} has a bounded drift. It is well known that for a starting point $Y_0 = y$, 
the process $Y_t$ has a density $\tilde{\rho}_t(y, y')$, which satisfies the following two-sided Gaussian estimate 
(for example, see \cite{CHXZ17}):
$$
\tilde{C}_0 t^{-d/2} \e^{-\tilde{c}_0 |y - y'|^2 / t} \leq \tilde{\rho}_t(y, y') \leq \tilde{C}_1 t^{-d/2} \e^{-\tilde{c}_1 |y - y'|^2 / t}.
$$
The desired estimate \eqref{HH1} then follows by a change of variables (see \cite{XZ21}).
\end{proof}

{Up to this point, we have assumed the condition $p_b,q_b<\infty$ in {\bf ($\widetilde{\bf H}^{\rm sub}_{b}$)}. It must be noticed that for $q=\infty$ or $p=\infty$, the space $C_b(\mR^d)$ is not dense in $L^\infty$, and consequently, the convergence property \eqref{re2} does not holds. 
However, when $q=\infty$ or $p=\infty$, for $f\in \mL^q_T\bB^\alpha_{p,\infty}$, it follows from \eqref{Ho1} that for any $\alpha'<\alpha$,
\begin{align}\label{re201}
    \|f_n-f\|_{\mL^q_T\bB^{\alpha'}_{p,\infty}}\lesssim \|f\|_{\mL^q_T\bB^\alpha_{p,\infty}} n^{\alpha'-\alpha}.
\end{align}
This provides an alternative form of convergence even in cases where $q=\infty$ or $p=\infty$.
}

Now  we can give
\begin{proof}[Proof of Theorem \ref{Th1}]
Since $b_1, \div b_1\in \mL^{q_b}_T\bH^{\alpha_b}_{p_b}$, noting that $\alpha_b>-2-\alpha_b$, by Lemma \ref{lemB1},
$$
\bH^{\alpha_b}_{p_b}\subset\cap_{\eps>0, q_b\in[1,\infty]}\bB^{\alpha_b-\eps}_{p_b,q_b},
$$ 
one sees that {\bf ($\widetilde{\bf H}^{\rm sub}_{b}$)} holds. If $q_b, p_b<\infty$, then
the desired results follow by the above theorem.

For the case $q_b=\infty$ or $p_b=\infty$, since $\frac2{q_b}+\frac{d}{p_b}<1+\alpha_b$, one can choose $\tilde{\alpha}_b\in(-1,\alpha_b)$ and $ \tilde{q}_b<q_b\le\infty$ so that $\frac2{\tilde{q}_b}+\frac{d}{p_b}<1+\tilde{\alpha}_b$. 
In particular, by \eqref{Sob1} we have
$$
\kappa_b=\|b_1\|_{\mL^{\tilde{q}_b}_T\bB^{\tilde{\alpha}_b}_{p_b,\tilde{q}_b}}\leq \|b_1\|_{\mL^{q_b}_T\bB^{\alpha_b}_{p_b,q_b}}<\infty,
$$
and 
$$
\tilde \kappa_b=\|\div b_1\|_{\mL^{\tilde{q}_b}_T\bB^{-2-\tilde{\alpha}_b}_{p_b,1}}\leq \|\div b_1\|_{\mL^{q_b}_T\bB^{\alpha_b}_{p_b,q_b}}<\infty,
$$
and by \eqref{re201},
$$
\|b_{1,n}-b_1\|_{\mL^{\tilde{q}_b}_T\bB^{\tilde{\alpha}_b}_{p_b,\tilde{q}_b}}\leq\|b_{1,n}-b_1\|_{\mL^{\tilde{q}_b}_T\bB^{(\alpha_b+\tilde{\alpha}_b)/2}_{p_b,\infty}}\leq C\|b_1\|_{\mL^{q_b}_T\bB^{\alpha_b}_{p_b,q_b}}n^{(\tilde{\alpha}_b-\alpha_b)/2},
$$
and similarly,
$$
\lim_{n\to\infty}\|\div b_{1,n}-\div b_1\|_{\mL^{\tilde{q}_b}_T\bB^{-2-\tilde{\alpha}_b}_{p_b,1}}=0.
$$
Thus one can check that Lemma \ref{Th45} still holds. The rest is the same as Theorem \ref{23thm:48}. 
\end{proof}

\section{Supercritical case: Proof of Theorem \ref{Th2}}\label{sec:03}

Throughout this section we still fix $T>0$ and always assume that supercritical assumption  {\bf(H$^{\rm sup}_b$)} in Subsection \ref{Sub32} holds.
In Subsection \ref{Sec2}, we focus on establishing the existence of weak solutions. The crucial step in this regard is demonstrating the tightness of approximation solutions. Notably, in the supercritical case, due to the drift being a distribution,  we are not able
to establish a Krylov estimate similar to that in Lemma \ref{Le42} for any moments. Instead, we can only establish a Krylov estimate for the second-order moment, as demonstrated in Lemma \ref{Le54}, which is particularly useful for taking limits. However, this alone is not sufficient to establish tightness. To overcome this difficulty, 
we employ the stopping time technique and  the strong Markov property of approximation SDEs.
In Subsection \ref{Sec3}, we discuss the uniqueness of generalized martingale problem and  ultimately complete  the proof of Theorem \ref{Th2}.
 \subsection{Existence of weak solutions}\label{Sec2}
We still consider the approximation SDE 
\begin{align}\label{S2:SDEn}
 X^n_t =X_0+\int_0^t b_n(s,X_s^n)\dif s+\sqrt2 W_t,
\end{align}
where $b_n=b*\phi_n$.
Recall the definition of $\bI_d$ in \eqref{XX1} and the parameter set $\Theta$ in \eqref{Pa1}. We first show the following uniform Krylov estimate.
\bl\label{Le51}
Under {\bf(H$^{\rm sup}_b$)}, for any $(\alpha,p,q)\in\bI_d$,
there is a constant $C=C(\Theta, \alpha,p,q)>0$ such that for all $t\in(0,T]$, stopping times $0\leq\tau_0\leq \tau_1\leq t$ and $f\in  C_c([0,t]\times\mR^d)$,
$$
\sup_{n\in\mN}\bE\left(\int^{\tau_1}_{\tau_0} f(s,X^n_s)\dif s\Big|\sF_{\tau_0}\right)\le C\|f\|_{\mL^q_t\bH^\alpha_p}.
$$
\el
\begin{proof}
Without loss of generality, we assume $f\in C^\infty_c([0,T]\times\mR^d)$.
Fix $t\in[0,T]$. Let $u_n\in\mL^\infty_TC^\infty_b(\mR^d)$ be the smooth solution of the following backward PDE:
\begin{align*}
\p_s u_n+\Delta u_n+b_n\cdot\nabla u_n=f,\quad u_n(t)=0.
\end{align*}
By \eqref{AF1}, we have
\begin{align}\label{S2:AA01}
\sup_{n}\Big(\|u_n\|_{\mL^\infty_t}+\|u_n\|_{\sV_t}\Big)\le C \|f\|_{\mL^q_t\bH^\alpha_p}.
\end{align}
Using It\^o's formula for $(s,x)\to u_n(s,x)$, one sees that for any stopping times $0\leq\tau_0\leq \tau_1\leq t$,
$$
u_n(\tau_1,X^n_{\tau_1})-u_n(\tau_0,X^n_{\tau_0})=\int_{\tau_0}^{\tau_1} f(s,X_s^n)\dif s+\sqrt2\int^{\tau_1}_{\tau_0}\nabla u_n(s,X^n_s)\dif W_s.
$$
Hence, by the optional stopping time theorem,
\begin{align}\label{BX1}
\bE\left(\int_{\tau_0}^{\tau_1}f(s,X_s^n)\dif s\Big|\sF_{\tau_0}\right)=\bE(u_n(\tau_1,X^n_{\tau_1})|\sF_{\tau_0})-u_n(\tau_0,X^n_{\tau_0}).
\end{align}
By \eqref{S2:AA01}, we immediately obtain the result.
\end{proof}
\br
Note that if $\alpha=0$, then we automatically have
$$
\bE\left|\int^{\tau_1}_{\tau_0} f(s,X^n_s)\dif s\right|\leq \bE\left(\int^{\tau_1}_{\tau_0} |f(s,X^n_s)|\dif s\right)\leq C\|f\|_{\mL^q_tL^p}.
$$
However, for $f\in \mL^q_t\bH^\alpha_p$ with $\alpha<0$ and in the supercritical case, we can not show any absolute moment estimate since 
we do not have the maximum estimate of $\|\nabla u\|_\infty$ as done in the subcritical case. Thus we have to carefully treat the distribution-valued $f$.
\er
We have the following tightness of the law $\mP_n$ of approximation solution $X^n_\cdot$ in $\mC_T$.
\bl\label{in:thm01}
Under {\bf(H$^{\rm sup}_b$)}, the family of probability measures $(\mP_n)_{n\in\mN}$ is tight. 
\el
\begin{proof}
Without loss of generality, we assume that $q_b<\infty$. In the case where $q_b=\infty$ and $d/p_b<2+\alpha_b$, one can choose $\tilde{q}_b<\infty$ 
such that $2/\tilde{q}_b+ d/p_b<2+\alpha_b$ and $\kappa_b=\|b\|_{\mL^{\tilde{q}_b}\bH^{\alpha_b}_{p_b}}<\infty$.

Let $\sT_T$ be the set of all stopping times bounded by $T$.
By Aldous' criterion of tightness, it suffices to prove that
$$
\lim_{\delta\downarrow 0}\sup_{\tau,\eta\in\sT_T,\eta\leq\tau\leq\eta+\delta}\sup_n\bE|X_\tau^n-X_\eta^n|=0.
$$
By the strong Markov property, we only need to show 
$$
\lim_{\delta\downarrow 0}\sup_{t\in[0,T]}\sup_{x_0\in\mR^d}\sup_{\tau\leq\delta}\sup_n\bE|X_\tau^{t,n}(x_0)-x_0|=0,
$$
where $X_\cdot^{t,n}(x_0)$ stands for the solution of SDE \eqref{S2:SDEn} with starting point $x_0$ and drift 
$$
b^t_n(s,\cdot)=b_n(t+s,\cdot).
$$
For $\eps\in(0,1)$, define
$$
h_\eps(x):=\sqrt{\eps^2+|x-x_0|^2}.
$$
By It\^o's formula, we have
\begin{align*}
\bE|X^{t,n}_\tau(x_0)-x_0|\leq \bE h_\eps(X^{t,n}_\tau(x_0))=\eps+\bE\left(\int^\tau_0(\Delta+b_n^t(s)\cdot\nabla)h_\eps(X^{t,n}_s(x_0))\dif s\right).
\end{align*}
Note that for some $C$ independent of $\eps$ and $t,x_0$,
\begin{align}\label{AP0}
|\nabla h_\eps(x)|=\frac{|x-x_0|}{\sqrt{\eps^2+|x-x_0|^2}}\leq 1,\ \  |\nabla^2h_\eps|\leq 2d^2/\eps.
\end{align}
Thus for $\tau\leq\delta$,
$$
\bE\left(\int^\tau_0\Delta h_\eps(X^{t,n}_s(x_0))\dif s\right)\leq\delta\|\Delta h_\eps\|_\infty\leq 2d\delta/\eps.
$$
Moreover, by Lemma \ref{Le51}, (iv) of Lemma \ref{lemB1} and \eqref{AP0}, we have
\begin{align*}
&\bE\left(\int^{\delta\wedge\tau}_0(b^t_n(s)\cdot\nabla h_\eps)(X^{t,n}_s(x_0))\dif s\right)
\leq C_1\| b^t_n\cdot\nabla h_\eps\|_{\mL^{q_b}_\delta\bH^{\alpha_b}_{p_b}} \\
&\qquad\leq C_2 \| b^t_n \|_{\mL^{q_b}_\delta\bH^{\alpha_b}_{p_b}} \|\nabla h_\eps\|_{C^1_b}
\leq C_3 \| b^t\|_{\mL^{q_b}_\delta\bH^{\alpha_b}_{p_b}}/\eps.
\end{align*} 
Combining the above calculations, we obtain
$$
\bE|X^{t,n}_\tau(x_0)-x_0|\leq \eps+2d\delta/\eps+C_3 \| b^t\|_{\mL^{q_b}_\delta\bH^{\alpha_b}_{p_b}}/\eps.
$$
Since 
$
\lim_{\delta\to 0}\sup_{t\in[0,T]}\int^{(t+\delta)\wedge T}_t\|b(s)\|^{q_b}_{\bH^{\alpha_b}_{p_b}}\dif s=0,
$ 
firstly letting $\delta\to 0$ and then $\eps\to 0$, we complete the proof.
\end{proof}

By the above lemma and Prokhorov's theorem, $(\mP_n)_{n\in\mN}$ is relatively compact in $\cP(\mC_T)$.
Next, we want to show that any accumulation point of the sequence $(\mP_n)_{n\in\mathbb{N}}$ is indeed a weak solution of SDE \eqref{in:SDE}. 
However, to take the limit, we need a different type of Krylov's estimate, which, in turn, requires that the initial distribution  has an $L^2$-distribution.

For approximation SDE \eqref{S2:SDEn}, it is well-known that for any $t\in(0,T]$, $X^n_t$ admits a smooth density $\rho_n(t,x)$ that satisfies the Fokker-Planck equation
\begin{align}\label{S2:FPE}
\p_t\rho_n=\Delta\rho_n-\div (b_n\rho_n).
\end{align} 
We need the following simple lemma.
\bl
Suppose $\rho_0\in L^2$ and let $\chi\in C^\infty_c(\mR^d)$ with $\chi=1$ on $B_1$ and $\chi=0$ on $B^c_2$.  We have
\begin{align}\label{S2:AA10}
\|\rho_n\|_{\sV_T}\leq (\tilde \kappa_b+1)\e^{\tilde \kappa_b}\|\rho_0\|_2,
\end{align}
where $\sV_T$ is defined in \eqref{VV1}, and
\begin{align}\label{S2:AA11}
\lim_{R\to\infty}\sup_n\|\rho_n(1-\chi(\cdot/R))\|_{\sV_T}=0.
\end{align}
\el
\begin{proof}
First of all, since $\div b_n\in \mL^2_TC^\infty_b$ and $\rho_0\in L^2$, it is standard to derive that
$$
\rho_n\in \mL^\infty_TL^2.
$$
Let $\varphi\in C^\infty_c(\mR^d)$. By multiplying both sides of \eqref{S2:FPE} by $\varphi$, we get
\begin{align}\label{AP7}
\p_t\rho_n\varphi=\Delta(\rho_n\varphi)-\div (b_n\rho_n\varphi)+g^\varphi_n,
\end{align}
where
$$
g^\varphi_n:=b_n\rho_n\cdot\nabla\varphi-2\nabla\rho_n\cdot\nabla\varphi-\rho_n\Delta\varphi.
$$ 
Multiplying both sides of \eqref{AP7} by $\rho_n\varphi$ and integrating on $\mR^d$, we have
$$
\frac12\p_t\|\rho_n\varphi\|_2^2+\|\nabla(\rho_n\varphi)\|_2^2=\<b_n\rho_n\varphi,\nabla(\rho_n\varphi)\>+\<g^\varphi_n,\rho_n\varphi\>.
$$
Noting that
$$
\<b_n\rho_n\varphi,\nabla(\rho_n\varphi)\>
=-\frac12\<\div b_n,(\rho_n\varphi)^2\>\leq\frac12\|\div b_n\|_\infty\|\rho_n\varphi\|^2_2\leq\frac12\|\div b\|_\infty\|\rho_n\varphi\|^2_2
$$
and
$$
\<g^{\varphi}_n,\rho_n\varphi\>
=\<b_n, \varphi\rho^2_n\nabla\varphi\>+\<\rho^2_n,|\nabla\varphi|^2\>,
$$
we have
\begin{align*}
\|\rho_n(t)\varphi\|^2_2+2\int^t_0\|\nabla(\rho_n\varphi)\|^2_2\dif s
&\leq\|\rho_0\varphi\|^2_2+\int^t_0\|\div b\|_\infty\|\rho_n\varphi\|^2_2\dif s\\
&+\int^t_0\Big(\<b_n, \varphi\rho^2_n\nabla\varphi\>+\<\rho^2_n,|\nabla\varphi|^2\>\Big)\dif s.
\end{align*}
By Gronwall's inequality, we derive
\begin{align}\label{S2:AA00}
\|\rho_n\varphi\|_{\sV_T}\leq (\tilde \kappa_b+1)\e^{\tilde \kappa_b}\left(\|\rho_0\varphi\|_2+\|\<b_n, \varphi\rho^2_n\nabla\varphi\>\|_{L^1_T}
+\|\<\rho^2_n,|\nabla\varphi|^2\>\|_{L^1_T}\right).
\end{align}
Let $\chi_R(x):=\chi(x/R).$
Replacing $\varphi$ in \eqref{S2:AA00}  by $\chi_R$ and letting $R\to\infty$, by the dominated convergence theorem, it is easy to derive 
$$
\|\rho_n\|_{\sV_T}\leq (\tilde \kappa_b+1)\e^{\tilde \kappa_b}\|\rho_0\|_2.
$$
On the other hand, replacing $\varphi$  in \eqref{S2:AA00}  by $\bar\chi_R:=1-\chi_R$, we get
$$
\|\rho_n\bar\chi_R\|_{\sV_T}\leq (\tilde \kappa_b+1)\e^{\tilde \kappa_b}\left(\|\rho_0\bar\chi_R\|_2+\|\<b_n, \bar\chi_R\rho^2_n\nabla\bar\chi_R\>\|_{L^1_T}
+\|\<\rho^2_n,|\nabla\bar\chi_R|^2\>\|_{L^1_T}\right).
$$
By Lemma \ref{Le34} and the definition of $\sV_T$, we have
\begin{align*}
\|\<b_n, \bar\chi_R\rho^2_n\nabla\bar\chi_R\>\|_{L^1_T}&\lesssim\|b_n\|_{\mL^{q_b}_T\bH^{\alpha_b}_{q_b}}\|\rho_n\|_{\sV_T}\|\bar\chi_R\nabla\bar\chi_R\rho_n\|_{\sV_T}\\
&\lesssim\|b\|_{\mL^{q_b}_T\bH^{\alpha_b}_{q_b}}\|\rho_n\|^2_{\sV_T}/R\lesssim\|\rho_0\|^2_2/R
\end{align*}
and
$$
\|\<\rho^2_n,|\nabla\bar\chi_R|^2\>\|_{L^1_T}\lesssim\|\rho_n\|^2_{\mL^2_T}/R\lesssim\|\rho_0\|^2_2/R.
$$
Hence,
$$
\sup_n\|\rho_n\bar\chi_R\|_{\sV_T}\lesssim\|\rho_0\bar\chi_R\|_2+\|\rho_0\|^2_2/R,
$$
which converges to zero as $R\to\infty$. The proof is complete.
\end{proof}
We have the following Krylov estimate.
\bl\label{Le54}
Suppose $\rho_0\in L^2$. Under {\bf(H$^{\rm sup}_b$)}, for any $(\alpha,p,q)\in\bI_d$,
there is a constant $C=C(\Theta, \alpha,p,q)>0$ such that for all $t\in(0,T]$ and $f\in C_c([0,t]\times\mR^d)$,
\begin{align}
\bE\left|\int^t_0 f(s,X^n_s)\dif s\right|^2\le C\|f\|^2_{\mL^q_t\bH^\alpha_p}\|\rho_0\|_2.
\label{SW22}
\end{align}
\el
\begin{proof}
Without loss of generality, we assume $f\in C^\infty_c([0,T]\times\mR^d)$.
Fix $t\in[0,T]$. Let $u$ solve the following backward PDE
$$
\p_s u_n+\Delta u_n+b_n\cdot\nabla u_n+f=0,\quad u_n(t)=0.
$$
By \eqref{BX1} with $\tau_0=s<t$ and $\tau_1=t$, we have
\begin{align*}
\bE\left(\int_s^tf(r,X_r^n)\dif r\Big|\sF_s\right)=u_n(s,X^n_s).
\end{align*}
Therefore,
\begin{align}
\bE\left|\int_0^tf(r,X^n_r)\dif r\right|^2&=2\bE\left[\int_0^t f(s,X_s^n)\left(\int_s^t f(r,X_r^n) \dif r\right) \dif s\right]\no\\
&=2\bE\left[\int_0^t f(s,X_s^n)\bE\left(\int_s^tf(r,X_r^n)\dif r\Big|\sF_s\right) \dif s\right]\no\\
&=2\bE\left[\int_0^t f(s,X_s^n)u_n(s,X_s^n)\dif s\right]\no\\
&=2\int^t_0\<f(s)u_n(s),\rho_n(s)\>\dif s\leq 2\|\<fu_n,\rho_n\>\|_{L^1_t}.\label{SW2}
\end{align}
By Lemma \ref{Le34}, we get
$$
\bE\left|\int_0^tf(s,X^n_s)\dif s\right|^2
\lesssim\|f\|_{\mL^q_t\bH^{\alpha}_{p}}\|\rho_n\|_{\sV_t}\|u_n\|_{\sV_t}.
$$
The desired estimate follows by \eqref{S2:AA10} and \eqref{S2:AA01}.
\end{proof}
\br
We are not able to show the following stronger estimate
$$
\bE\left(\sup_{t\in[0,T]}\left|\int^t_0 f(s,X^n_s)\dif s\right|^2\right)\le C\|f\|^2_{\mL^q_T\bH^\alpha_p}\|\rho_0\|_2
$$
since we are using the estimate \eqref{S2:AA10} of the density and $f$ is a distribution.
\er
Now we can show the following existence result of weak solutions.
\bt\label{Th5}\label{Th57}
Suppose that $\rho_0\in L^2$. Under {\bf(H$^{\rm sup}_b$)}, there is a weak solution to SDE \eqref{in:SDE} 
such that for all $t\in(0,T]$ and $f\in  C_c([0,t]\times\mR^d)$,
\begin{align}\label{SW1}
\bE\left|\int^t_0 f(s,X_s)\dif s\right|^2\le C\|f\|^2_{\mL^q_t\bH^\alpha_p}\|\rho_0\|_2,
\end{align}
where $C=C(\Theta, \alpha,p,q)>0$.
Moreover, for each $t\in(0,T]$, $X_t$ admits a density $\rho_t\in L^2$ with 
\begin{align}\label{1218:001}
    \|\rho\|_{\sV_T}\leq(\tilde \kappa_b+1)\e^{\tilde \kappa_b}\|\rho_0\|_2
\end{align}
and in the distributional sense
\begin{align}\label{AQ1}
\p_t\rho=\Delta\rho-\div (b\rho).
\end{align}
\et
\begin{proof}
By Lemma \ref{in:thm01} and Prokhorov's criterion, there exists a subsequence $n_k$ and a probability measure $\mP \in \cP(\mC_T)$ 
such that $\mP_{n_k}$ weakly converges to $\mP$ as $k \to \infty$. Without loss of generality, we assume that $\mP_n$ weakly converges 
to $\mP$ as $n \to \infty$. By Skorokhod's representation theorem, if necessary, changing the probability space, 
we may assume that there exists a continuous process $(X_t)_{t \in [0,T]}$ such that
$$
X_\cdot^n \to X_\cdot \quad \text{a.s. in $\mC_T$.}
$$
Estimate \eqref{SW1} follows by taking limits in \eqref{SW22}. To show that $X$ is a weak solution of SDE \eqref{in:SDE}, 
it suffices to show that
\begin{align}\label{DZ1}
\lim_{n \to \infty} \bE\left| \int_0^t b_n(r, X^n_r) \dif r - \int_0^t b_n(r, X_r) \dif r \right|^2 = 0.
\end{align}
Notice that for fixed $m \in \mN$, by the dominated convergence theorem,
$$
\lim_{n \to \infty} \bE\left| \int_0^t b_m(r, X^n_r) \dif r - \int_0^t b_m(r, X_r) \dif r \right|^2 = 0,
$$
and once we have shown 
\begin{align}\label{AP9}
\lim_{m, m' \to \infty} \sup_n \bE\left| \int_0^t (b_m - b_{m'})(r, X^n_r) \dif r \right|^2 = 0,
\end{align}
then \eqref{DZ1} follows. Moreover, we also have
$$
\lim_{n \to \infty} \int_0^t b_n(r, X_r) \dif r =: A^b_t \quad \text{in $L^2$}.
$$
Let $\chi\in C^\infty_c(\mR^d)$ with $\chi=1$ on $B_1$ and $\chi=0$ on $B^c_2$. For $R>0$, define cutoff function
$$
\chi_R(x):=\chi(x/R).
$$
Since $\frac{d}{p_b}+\frac2{q_b}<2+\alpha$, one can choose finite $\bar q_b,\bar p_b\in[2,\infty)$ with $\bar q_b\leq q_b, \bar p_b\leq p_b$ and such that $\frac{d}{\bar p_b}+\frac2{\bar q_b}<2+\alpha$.
Then for fixed $R>0$, by Lemma \ref{Le54}, we have
$$
\lim_{m,m'\to\infty}\sup_n\bE\left|\int_0^t((b_m-b_{m'})\chi_R)(r,X^n_r)\dif r\right|^2\lesssim\lim_{m,m'\to\infty}\|(b_m-b_{m'})\chi_R\|^2_{\mL^{\bar q_b}_t\bH^{\alpha_b}_{\bar p_b}}=0.
$$
Thus for proving \eqref{AP9}, it remains to show that for $\bar\chi_R(x):=1-\chi_R(x)$,
\begin{align}\label{KM1}
\lim_{R\to\infty}\sup_{m,m'}\sup_n\bE\left|\int_0^t((b_m-b_{m'})\bar\chi_R)(r,X^n_r)\dif r\right|^2=0.
\end{align}
Let $u^{m,m'}_{n,R}=u$ solve the following backward PDE:
$$
\p_s u+\Delta u+b_n\cdot\nabla u_n+(b_m-b_{m'})\bar\chi_R=0,\quad u(t)=0.
$$
By \eqref{SW2} and Lemma \ref{Le34}, we have
\begin{align*}
&\bE\left|\int_0^t((b_m-b_{m'})\bar\chi_R)(r,X^n_r)\dif r\right|^2
\leq2\|\<(b_m-b_{m'})\bar\chi_Ru^{m,m'}_{n,R},\rho_n\>\|_{L^1_t}\\
&\qquad\lesssim\|b_m-b_{m'}\|_{\mL^{q_b}_t\bH^{\alpha_b}_{p_b}}\|\bar\chi_R\rho_n\|_{\sV_t}\|u^{m,m'}_{n,R}\|_{\sV_t}
\lesssim\|b\|^2_{\mL^{q_b}_t\bH^{\alpha_b}_{p_b}}\|\bar\chi_R\rho_n\|_{\sV_t},
\end{align*}
which implies \eqref{KM1} by   \eqref{S2:AA11}.
Finally, for \eqref{AQ1}, as in the proof of existence in Theorem \ref{S1:well}, it follows by \eqref{S2:FPE} and \eqref{S2:AA10}.
\end{proof}
 \subsection{Generalized martingale problems}\label{Sec3}
Let $(\alpha,p,q)\in\bI_d$ (see \eqref{XX1} for a definition). For any $f\in\mL^q_T\bH^\alpha_p$, by Theorem \ref{S1:well}, there exists a weak solution to the following backward PDE
 \begin{align}\label{in:PDE}
\p_t u+\Delta u+b\cdot\nabla u=f,\quad u(T)\equiv0
\end{align}
with 
$$
\|u\|_{\mL^\infty_T\cap\sV_T}\leq C\|f\|_{\mL^q_T\bH^\alpha_p}.
$$
Note that the weak solution may be not unique  if we do not assume $b\in\mL^\infty_T\cB+\mL^2_T L^2$. However, if we consider the following regularized PDE
 \begin{align}\label{in:PDE0}
\p_t u_n+\Delta u_n+b_n\cdot\nabla u_n=f,\quad u_n(T)\equiv0.
\end{align}
Then for each $n$, there is a unique solution $u_n$ denoted by $\cS_f^n$ so that
\begin{align}\label{1205:00}
\|\cS_f^n\|_{\mL^\infty_T\cap\sV_T}\leq C\|f\|_{\mL^q_T\bH^\alpha_p},
    \end{align}
where $C$ is independent of $n$.
In other words, $f\mapsto\cS_f^n$ is a bounded linear operator from $\mL^q_T\bH^\alpha_p$ to $\mL^\infty_T\cap\sV_T$. 

Now suppose $p,q<\infty$. Let $\{f_m\}_{m=1}^\infty\subset \mL^q_T\bH^\alpha_p$ be a countable dense subset. 
By the Cantor-Hilbert selection and from the proof of Theorem \ref{S1:well}, there is a subsequence $\{n_k\}$ such that for each $f=f_m$, there is a solution $\cS_{f}=:u$ to PDE \eqref{in:PDE}
such that for any $t\in[0,T]$ and $\varphi\in L^1(\mR^d)$,
 \begin{align}\label{SW40}
\lim_{k\to\infty}\<\cS_f^{n_k}(t),\varphi\>=\<\cS_f(t),\varphi\>,
\end{align}
and for each $R>0$,
 \begin{align}\label{SW50}
\lim_{k\to\infty}\int^T_0\|\cS_f^{n_k}(s)-\cS_f(s)\|^2_{L^2(B_R)}\dif s=0.
\end{align}
Moreover, by \eqref{1205:00} and \eqref{SW40}, we have
$$
    \|\cS_{f_m}\|_{\mL^\infty_T\cap\sV_T}\leq C\|f_m\|_{\mL^q_T\bH^\alpha_p},\quad m\in \mN.
$$
Thus $\cS_f$ can be extended as a bounded linear operator  from
$\mL^q_T\bH^\alpha_p$ to $\mL^\infty_T\cap\sV_T$ so that \eqref{SW40} and \eqref{SW50} holds for all $f\in \mL^q_T\bH^\alpha_p$.

In the following we fix such an operator $\cS$ and construct a unique generalized martingale solution associated with $\cS$. 
Of course, its definition depends on the choice of subsequence $n_k$ as well as the mollifiers $\phi_{n_k}$.
We emphasize that when $b\in\mL^\infty_T\cB+\mL^2_T L^2$, $\cS_f$ is nothing, but the unique solution of \eqref{in:PDE}, which is independent of the choice of $n_k$ and the mollifiers.

 \bd\label{in:def}
Let $\mu \in \cP(\mR^d)$. We call a probability measure $\mP \in \cP(\mC_T)$ a generalized martingale solution of SDE \eqref{in:SDE} starting from $\mu$ and associated with the operator $\cS$, 
if $\mP \circ (\ww_0)^{-1} = \mu$, 
and for each $t \in (0,T)$, $\mP \circ w_t^{-1}$ has a density with respect to the Lebesgue measure, 
and for any $f \in C^\infty_c([0,T] \times \mR^d)$,
$$
M_t^f := \cS_f(t, \ww_t) - \cS_f(0, \ww_0) - \int_0^t f(r, \ww_r) \dif r, \quad w_\cdot \in \mC_T,
$$
is an almost surely martingale under $\mP$ with respect to the natural filtration $\sB_s$, in the sense that there is a Lebesgue null set $\sN \subset (0,T)$ such that for all $0 \leq s < t \leq T$ with $s,t \notin \sN$,
$$
\mE(M^f_t | \sB_s) = M^f_s.
$$
\ed
\br
Although $t \mapsto M^f_t$ is defined pointwise, we lack information about the path properties of $t \mapsto M^f_t$.  
For instance, we have no knowledge of its c\'adl\'ag property since we only have $\cS_f \in \mL^\infty_T$.
\er

Now we show the following main result.
\bt\label{in:main}
Assume   {\bf(H$^{\rm sup}_b$)} holds. If the initial distribution $\mu$ has an $L^2$-density w.r.t. the Lebesgue measure, then there is a unique 
generalized martingale solution to SDE \eqref{in:SDE} associated with the operator $\cS$ in the sense of Definition \ref{in:def}.
\et
\begin{proof}

{\bf (Existence)} 
Let $\mP$ be any accumulation point of $(\mP_n)_{n\in\mN}$. We want to show that $\mP$ is a martingale solution  in the sense of Definition \ref{in:def}.
Up to taking a subsequence, without loss of generality we assume that
\begin{align*}
\mP_n\to \mP \quad\text{as $n\to\infty$ in $\cP(\mC_T)$ weakly.}
\end{align*}
 In particular, by \eqref{S2:AA10}, for each $t\in(0,T)$, $\mP\circ(w_t)^{-1}$ has a density w.r.t. the Lebesgue measure.

For $f\in C^\infty_c([0,T]\times\mR^d)$, let $u_n=\cS^n_f$ solve the following backward PDE:
$$
\p_t u_n+\Delta u_n+b_n\cdot\nabla u_n=f,\ \ u_n(T)=0.
$$
By It\^o's formula, we have
\begin{align}\label{SS1}
u_n(t,X^n_t)=u_n(0,X^n_0)+\int_0^t f(s,X^n_s)\dif s+\sqrt2\int^t_0\nabla u_n(s,X^n_s)\dif W_s.
\end{align}
Define
\begin{align*}
M^n_t:=u_n(t,\ww_t)-u_n(0,\ww_0)-\int_0^t f(r,\ww_r)\dif r.
\end{align*}
Since $\mP_n=\bP\circ(X^n_\cdot)^{-1}$, by \eqref{SS1} one sees that $M^n_t$ is a martingale under $\mP_n$ with respect to  $\sB_s$. 
More precisely, for any $s<t$ and bounded continuous $\sB_s$-measurable functional $G_s$,
$$
\mE^{\mP_n}((M^n_t-M^n_s) G_s)=0,
$$
equivalently,
$$
\mE^{\mP_n}\Big(((u_n(t,\ww_t)-u_n(s,\ww_s))G_s\Big)=\mE^{\mP_n}\left( G_s\int_s^t f(r,\ww_r)\dif r\right).
$$
Clearly, we have
\begin{align}\label{S3:AA00}
\lim_{n\to\infty}\mE^{\mP_n}\left( G_s\int_s^t f(r,\ww_r)\dif r\right)=\mE^{\mP}\left( G_s\int_s^t f(r,\ww_r)\dif r\right).
\end{align}
By \eqref{SW50}, there is a Lebesgue null set $\sN\subset[0,T]$ depending on $f$ such that for all $t\notin\sN$ and $R>0$,
\begin{align}\label{AF2}
\lim_{k\to\infty}\|\cS^{n_k}_f(t)-\cS_f(t)\|_{L^2(B_R)}=0.
\end{align}
We want to show that for $s,t\notin\sN$ with $s<t$,
\begin{align}\label{AF3}
\lim_{n\to\infty}\mE^{\mP_n}\Big((u_n(t,\ww_t)-u_n(s,\ww_s))G_s\Big)=\mE^{\mP}\Big((\cS_f(t,\ww_t)-\cS_f(s,\ww_s))G_s\Big).
\end{align}
For fixed $m\in\mN$, since $x\mapsto u_m(t,x)$ is bounded continuous, we have
$$
\lim_{n\to\infty}\mE^{\mP_n}\Big((u_m(t,\ww_t)G_s\Big)=\mE^{\mP}\Big(u_m(t,\ww_t)G_s\Big).
$$
On the other hand, noting that for $R>0$,
\begin{align*}
\mE^{\mP_n}\big|(u_m(t,\ww_t)-\cS_f(t,\ww_t))\b1_{B_R}(\ww_t)\big|
&=\int_{B_R}\big|u_m(t,x)-\cS_f(t,x)\big|\rho_n(t,x)\dif x\\
&\leq\|u_m(t,\cdot)-\cS_f(t,\cdot)\|_{L^2(B_R)}\|\rho_n(t)\|_{L^2(B_R)},
\end{align*}
by \eqref{AF2} and \eqref{S2:AA00} we have
$$
\lim_{m\to\infty}\sup_n\mE^{\mP_n}\big|(u_m(t,\ww_t)-\cS_f(t,\ww_t))\b1_{B_R}(\ww_t)\big|=0
$$
and
$$
\lim_{m\to\infty}\mE^{\mP}\big|(u_m(t,\ww_t)-\cS_f(t,\ww_t))\b1_{B_R}(\ww_t)\big|=0.
$$
Moreover, we also have
\begin{align*}
\mE^{\mP_n}\big|(u_m(t,\ww_t)-\cS_f(t,\ww_t))\b1_{B^c_R}(\ww_t)\big|
&=\int_{B^c_R}\big|u_m(t,x)-\cS_f(t,x)\big|\rho_n(t,x)\dif x\leq\|u_m(t)\|_{\infty}\mu^n_t(B^c_R),
\end{align*}
where $\mu^n_t=\mP^n\circ w^{-1}_t$.
By the tightness of $(\mu^n_t)_{n\in\mN}$, we have
$$
\lim_{R\to\infty}\sup_{m,n}\mE^{\mP_n}\big|(u_m(t,\ww_t)-\cS_f(t,\ww_t))\b1_{B^c_R}(\ww_t)\big|\leq C\lim_{R\to\infty}\sup_{n}\mu^n_t(B^c_R)=0.
$$
Thus we get \eqref{AF3}. Hence, for all $s,t\notin\sN$ with $s<t$ and any bounded continuous $\sB_s$-measurable $G_s$,
$$
\mE^{\mP}\Big((\cS_f(t,\ww_t)-\cS_f(s,\ww_s))G_s\Big)=\mE^{\mP}\left( G_s\int_s^t f(r,\ww_r)\dif r\right),
$$
which gives
$$
\mE(M^f_t|\sB_s)=M^f_s.
$$
In particular, by noticing that $u_n(T)=\cS_f(T)=0$, \eqref{AF3} holds for $t=T$ and $s\notin \sN$. Now, we show \eqref{AF3} holds for $s=0$ and $t\notin \sN$ as well. 
We only need to show that for any  bounded continuous function $g:\mR^d\to\mR$,
\begin{align*}
\lim_{n\to\infty}\mE^{\mP_n}u_n(0,\ww_0)g(\ww_0)=\mE^{\mP}\cS_f(0,\ww_0)g(\ww_0).
\end{align*}
We note that $\mP_n\circ(\ww_0)^{-1}=\mP\circ(\ww_0)^{-1}=\mu_0$ and by \eqref{SW40}, 
\begin{align*}
\lim_{n\to\infty}\<u_n(0),\rho_0g\>=\<\cS_f(0),\rho_0g\>.
\end{align*}
Therefore, we have $\mE(M^f_t|\sF_0)=0$.
The existence of a generalized martingale solution is proven.

{\bf (Uniqueness)} Let $\mP_1,\mP_2$ be two generalized martingale solutions with initial  distribution $\mu_0(\dif x)=\rho_0(x)\dif x$ in the sense of Definition \ref{in:def}. 
By the definition, for any $f\in C_c([0,T]\times\mR^d)$, there is a Lebesgue full measure set $\sA_{f}\subset[0,T]$ containing $0$ and $T$ such that for any $s,t\in \sA_f$ and bounded $\sB_s$-measurable $G_s$,
\begin{align}\label{S3:CC00}
\mE^{\mP_i} ((M^f_t-M^f_s)G_s)=0,\ \ i=1,2.
\end{align}
Since $C_c([0,T]\times\mR^d)$ is separable, we can choose $\sA=\sA_{f}$ independent of $f$.

Firstly,  taking $(s,t)=(0,T)$, $G\equiv1$ in \eqref{S3:CC00}, we have
 \begin{align*}
\mE^{\mP_1}\left(\int_0^T f(t,\ww_t)\dif t\right)=-\<\cS_f(0),\rho_0\>=\mE^{\mP_2}\left(\int_0^T f(t,\ww_t)\dif t\right).
\end{align*}
which implies that $\mP_1$ and $\mP_2$ have the same one dimensional time marginal distribution. In fact, taking $f(t,x)=h(t)g(x)$ in the above equality, 
where $h\in C([0,T])$ and $g\in C_c(\mR^d)$, by Fubini's theorem, one finds that
for Lebesgue almost all $t\in(0,T)$ and any $g\in C_c(\mR^d)$,
$$
\mE^{\mP_1} g(\ww_t)=\mE^{\mP_2} g(\ww_t).
$$
Since $t\mapsto \ww_t$ is continuous, we further have for any $t\in[0,T]$ and $g\in C_c(\mR^d)$,
$$
\mE^{\mP_1} g(\ww_t)=\mE^{\mP_2} g(\ww_t)\Rightarrow \mP_1\circ \ww_t^{-1}=\mP_2\circ \ww_t^{-1}.
$$

Below we show that $\mP_1$ and $\mP_2$ have the same finite  dimensional time marginal distribution by induction.
We assume that for some $n\in\mN$, $\mP_1$ and $\mP_2$ have the same $n$-dimensional time marginal distribution, that is, for any 
$0\leq t_1<t_2<\cdots<t_n\leq T$,
$$
\mP_1\circ(\ww_{t_1},\cdots, \ww_{t_n})^{-1}=\mP_2\circ(\ww_{t_1},\cdots, \ww_{t_n})^{-1}. 
$$
Let $t_0=0$ and $0\leq t_1<t_2<\cdots<t_{n-1}\leq T$ be fixed. Let $g\in C_c(\mR^{nd})$ and $t_n\in (t_{n-1},T)\cap\sA$.
By \eqref{S3:CC00} with $(s, t)=(t_n,T)$ and $G_s=g(\ww_{t_1}\cdots,\ww_{t_n})$, we have
\begin{align}\label{0112:00}
\mE^{\mP_i}\left(\int_{t_n}^T f(s,\ww_s)\dif s\,g(\ww_{t_1}\cdots,\ww_{t_n})\right)=-\mE^{\mP_i}\Big(\cS_f(t_n,\ww_{t_n})g(\ww_{t_1}\cdots,\ww_{t_n})\Big),\quad i=1,2,
\end{align}
which implies by the induction hypothesis that
\begin{align*}
\mE^{\mP_1}\left(\int_{t_n}^T f(s,\ww_s)\dif s\,g(\ww_{t_1}\cdots,\ww_{t_n})\right)
=\mE^{\mP_2}\left(\int_{t_n}^T f(s,\ww_s)\dif s\,g(\ww_{t_1}\cdots,\ww_{t_n})\right).
\end{align*}
Therefore, for any $t_n\in (t_{n-1},T)\cap\sA$ and Lebesgue almost all $t_{n+1}\in(t_{n},T]$,
$$
\mE^{\mP_1}\left(g_1(\ww_{t_{n+1}})g(\ww_{t_1}\cdots,\ww_{t_n})\right)
=\mE^{\mP_2}\left(g_1(\ww_{t_{n+1}})g(\ww_{t_1}\cdots,\ww_{t_n})\right).
$$
By the continuity of $t\mapsto \ww_t$,  we get for all $t_{n}, t_{n+1}\in(t_{n-1}, T]$ with $t_n<t_{n+1}$,
$$
\mP_1\circ(\ww_{t_1},...,\ww_{t_n}, \ww_{t_{n+1}})^{-1}=\mP_2\circ(\ww_{t_1},...,\ww_{t_n}, \ww_{t_{n+1}})^{-1}.
$$
Namely, $\mP_1=\mP_2$. The proof is complete.
\end{proof}

Now we can give
\begin{proof}[Proof of Theorem \ref{Th2}]
Let $n_k$ be any subsequence and $n_k'$ be the subsubsequence of $n_k$ used in the definition of $\cS_f$.
By Theorems \ref{in:main} and  \ref{Th5}, any accumulation point of $\mP_{n'_k}$ is a generalized martingale solution, and meanwhile, a weak solution.
By the uniqueness of generalized martingale solutions, for any initial value $\mu$, without further selecting a subsequence, $\mP_{n'_k}$ weakly converges to the weak solution of SDE \eqref{in:SDE}.
If $b=b_1+b_2$, where $b_1\in \mL^\infty_T\cB$ and $b_2\in \mL^2_TL^2$, then PDE \eqref{in:PDE} has a unique weak solution. Thus, no necessary to select a subsequence, $\cS_f$ is well-defined, and 
as above, the whole sequence $\mP_n$ weakly converges to the weak solution of SDE \eqref{in:SDE}. 

Finally, we show the Markov property by showing it for the generalized martingale solution.  Since $C^\infty_c([0,T]\times\mR^d)$ is separable, by Definition \ref{in:def},
we can find a Lebesgue zero set $\sN\subset(0,T)$ such that for any $f\in C^\infty_c([0,T]\times\mR^d)$ and $s\in [0,T]\backslash\sN$,
$$
\mE^{\mP_\mu}\Big[M^f_T-M^f_s\big|\sB_s\Big]=0{=\mE^{\mP_\mu}\left[\mE^{\mP_\mu}\left[M^f_T-M^f_s\big|\sB_s\right]\big|w_s\right]}=\mE^{\mP_\mu}\Big[M^f_T-M^f_s\big|w_s\Big],
$$
which in turn implies by $\cS_f(T)=0$ that
$$
\mE^{\mP_\mu}[\cS_f(s,w_{s})|\sB_s]+\mE^{\mP_\mu}\left[\int_{s}^T f(r,w_r)\dif r\Big|\sB_s\right]=\mE^{\mP_\mu}[\cS_f(s,w_{s})|w_s]+\mE^{\mP_\mu}\left[\int_{s}^T f(r,w_r)\dif r\Big|w_s\right].
$$
Hence, for any $f\in C^\infty_c([0,T]\times\mR^d)$ and $s\in [0,T]\backslash\sN$,
$$
\mE^{\mP_\mu}\left[\int_{s}^T f(r,w_r)\dif r\Big|\sB_s\right]=\mE^{\mP_\mu}\left[\int_{s}^Tf(r,w_r)\dif r\Big|w_s\right].
$$
In particular, for fixed $t\in(s,T]$ and $f\in C_c(\mR^d)$, by mollifying approximation, we further have
$$
\mE^{\mP_\mu}\left[\int_{(t-1/n)\vee s}^t f(w_r)\dif r\Big|\sB_s\right]=\mE^{\mP_\mu}\left[\int_{(t-1/n)\vee s}^tf(w_r)\dif r\Big|w_s\right],\ \ n\in\mN.
$$
Multiplying both sides by $n$ and letting $n\to\infty$, by the dominated convergence theorem, we obtain 
\begin{align*}
\mE^{\mP_\mu}[f(w_t)|\sB_s]=\mE^{\mP_\mu}[f(w_t)|{w_s}].
\end{align*}
The proof is complete.
\end{proof}

\section{Applications to diffusion in random environment}\label{sec6}

In this section, we introduce certain Gaussian noise $b$ so that our results can be applied to them. In particular, we construct the diffusion process in random environment using this noise.

\subsection{Vector-valued Gaussian fields}
In this section we introduce the notion of vector-valued Gaussian field.
Fix $m\in\mN$ and let $\mu(\dif\xi):=\{\mu_{ij}(\dif\xi)\}_{i,j=1}^m$ be a matrix-valued signed Radon measure on $\mR^{d}$. We assume that
\begin{enumerate}
\item[{\bf (A)}]  For each Borel measurable set $A\subset\mR^d$, $(\mu_{ij}(A))_{i,j=1}^m$ is a symmetric semi-positive definite Hermitian matrix, that is, for any complex numbers $(a_i)^m_{i=1}$,
\begin{align}\label{AG32}
\mu_{ij}(A)=\mu_{ji}(A),\  \ \sum_{i,j=1}^ma_i\bar a_j\mu_{ij}(A)\geq 0.
\end{align}
Moreover, $\mu(\dif\xi)=\mu(-\dif\xi)$ and for some $\ell\in\mR$,
\begin{align}\label{AG3}
\sup_{i,j=1,\cdots, m}\int_{\mR^d}(1+|\xi|)^{\ell}|\mu_{ij}|_{\rm var}(\dif\xi)<\infty,
\end{align}
where $|\cdot|_{\rm var}$ stands for the total variation measure of a signed measure.
\end{enumerate}
For two $\mR^m$-valued Schwartz functions $f,g\in\sS(\mR^{d};\mR^m)$ and $k=0,1,2,\cdots$, we define
$$
\<f,g\>_{\mu; k}:=\sum_{i,j=1}^m\int_{\mR^{d}}(1+|\xi|^{2k})\widehat{f_i}(\xi)\overline{\widehat{g_j}(\xi)}\mu_{ij}(\dif\xi).
$$
Since $\xi\to \mu_{ij}(\dif\xi)$ is symmetric, one sees that $\<f,g\>_{\mu; k}=\<g,f\>_{\mu; k}\in\mR$ is an inner product on $\sS(\mR^{d};\mR^m)$,
and by \eqref{AG32}, for any $f\in\sS(\mR^{d};\mR^m)$,
\begin{align*}
\|f\|_{\mu; k}:= \left(\sum_{i,j=1}^m\int_{\mR^{d}}(1+|\xi|^{2k})\widehat{f_i}(\xi)\overline{\widehat {f_j}(\xi)}\mu_{ij}(\dif\xi)\right)^{1/2}\geq 0.
\end{align*}
Hence, $\|\cdot\|_{\mu;k}$ is a seminorm on linear $\sS(\mR^{d};\mR^m)$. Clearly,
$$
\|\cdot\|_\mu:=\|\cdot\|_{\mu;0}\leq\|\cdot\|_{\mu;k}\leq \|\cdot\|_{\mu;k+1}.
$$ 
It must be noticed that $\|f\|_{\mu}=0$ does not imply $f=0$.
Let $E_0$ be the null space of $\|\cdot\|_{\mu}$, i.e.,
$$
E_0:=\{f\in \sS(\mR^{d};\mR^m): \|f\|_{\mu}=0\}.
$$
We define $\mH_k$ being the completion of the quotient space $\sS(\mR^{d};\mR^m)/E_0$ with respect to the  seminorm $\|\cdot\|_{\mu; k}$.
Thus we obtain a Hilbert space and 
$$
\mH_{k+1}\subset\mH_k\subset\mH_0.
$$ 
\br\label{Re60}
Assume $\ell\geq 0$ in \eqref{AG3}, then the Dirac measure $\delta^{\otimes m}_0$ belongs to $\mH_0$. Indeed, for $f_n(x)=n^{d}f(nx)$, where $f\in\sS(\mR^d;\mR^m)$ satisfies 
$\int_{\mR^d} f_i(x)\dif x=1$ for each $i=1,\cdots,m$, we have
$$
\|f_n-f_{n'}\|^2_\mu=\sum_{i,j=1}^m\int_{\mR^{d}}(\widehat{f_i}(\xi/n)-\widehat{f_i}(\xi/n'))\overline{(\widehat {f_j}(\xi/n)-\widehat{f_j}(\xi/n'))}\mu_{ij}(\dif\xi).
$$
Since $\ell\geq 0$ in \eqref{AG3} and $\widehat{f_i}(0)=1$, by the dominated convergence theorem, we get 
$$
\lim_{n\to\infty}\|f_n-f_{n'}\|_\mu=0.
$$
In general, elements in $\mH_0$ might not be Schwartz distribution.
\er

For any $x\in\mR^d$, $k\in\mN$ and $f\in\sS(\mR^d;\mR^m)$, define
$$
\tau_x f(\cdot):=f(\cdot-x)\in \sS(\mR^d;\mR^m),\ \ \nabla^k f(\cdot):=\p^k_x f(\cdot)\in\sS(\mR^d;\mR^m\otimes(\mR^d)^{\otimes k}).
$$
Noting that
$$
\widehat{\tau_x f}(\xi)=\e^{-{\rm i}x\cdot\xi}\hat{f}(-\xi),\ \ \widehat{\nabla^k f}(\xi)={\rm i}^{k}\xi^{\otimes k} \hat{f}(\xi),
$$
we clearly have
$$
\|\tau_xf\|_{\mH_0}\leq\|f\|_{\mH_0},\ \ \|\nabla^k f\|_{\mH_0}\leq \|f\|_{\mH_k}.
$$
In particular, we can extend $\tau_x: \mH_0\to\mH_0$ and $\nabla^k:\mH_k\to\mH_0$ being bounded linear operators.

We have the following simple lemma.
\bl\label{Le61}
For any $f\in\mH_0$, the function $x\mapsto\tau_x f\in\mH_0$ is continuous, and
for any $f\in\mH_k$, $x\mapsto \tau_x f\in\mH_0$ is $k$-order Fr\'echet differentiable and
\begin{align}\label{1127:00}
\nabla_x^k\tau_xf=\tau_x\nabla^k f\ \  \text{in $\mH_0$}.
\end{align}
\el
\begin{proof}
By the boundedness  $\|\tau_xf\|_{\mH_0}\leq\|f\|_{\mH_0}$, we may assume $f\in\sS(\mR^d;\mR^m)$. Thus, by definition, 
$$
\|\tau_x f-\tau_yf\|^2_{\mu}=\sum_{i,j=1}^m\int_{\mR^d}|\e^{-{\rm i}\xi\cdot x}-\e^{-{\rm i}\xi\cdot y}|^2\widehat  f_i(\xi)
\overline{\widehat  f_j(\xi)}\mu_{ij}(\dif \xi),
$$
which implies by the dominated convergence theorem,
$$
\lim_{x\to y}\|\tau_x f-\tau_yf\|^2_{\mu}=0.
$$
For \eqref{1127:00},  it suffices to show it for $k=1$. For $h,x\in\mR^d$ with $h\not=0$, by definition we have
\begin{align*}
&\left\|\frac{\tau_{x+h} f-\tau_x f}{h}-\tau_x\nabla  f\right\|^2_{\mu}= 
\sum_{i,j=1}^m\int_{\mR^d}\left|\frac{\e^{-{\rm i}\xi\cdot h}-1}{h}+{\rm i}\xi\right|^2\widehat  f_i(\xi)
\overline{\widehat  f_j(\xi)}\mu_{ij}(\dif \xi).
\end{align*}
Noting that
\begin{align*}
\sum_{i,j=1}^m\int_{\mR^d}\sup_h\left|\frac{\e^{-{\rm i}\xi\cdot h}-1}{h}+{\rm i}\xi\right|^2|\widehat  f_i(\xi)\overline{\widehat  f_j(\xi)}\mu_{ij}(\dif \xi)
\le  \sum_{i,j=1}^m\int_{\mR^d}|\xi|^{2}\widehat  f_i(\xi)\overline{\widehat  f_j(\xi)}\mu_{ij}(\dif \xi)<\infty,
\end{align*}
by the dominated convergence theorem, we have
$$
\lim_{h\to 0}\left\|\frac{\tau_{x+h} f-\tau_x f}{h}-\tau_x\nabla  f\right\|_{\mu}=0.
$$
The proof is complete.
\end{proof}

Now we introduce the  following notion of vector-valued Gaussian random field (see \cite{SV97}). 
\bd
Let $U$ be the real valued Gaussian random field on $\mH_0$, that is, $U$ is a continuous linear operator from $\mH_0$ to $L^2(\Omega,\bP)$,
and for each $f\in\mH_0$,
$U(f)$ is a real-valued Gaussian random variable with mean zero and variance $\|f\|^2_{\mH_0}$. In particular,
\begin{align}\label{ISO}
\bE\big(U(f)U(g)\big)=\<f,g\>_{\mH_0}.
\end{align}
We shall call $U$ an $\mR^m$-valued Gaussian noise over $\mR^d$ with matrix-valued
spectral measure $\mu$.
\ed

We call $U$ an $\mR^m$-valued Gaussian noise. The reason is that 
 for given real valued Schwartz function $h\in\sS(\mR^d;\mR)$ and $i=1,\cdots,m$, if we define
\begin{align}\label{EE9}
f_i:=(0,\cdots, 0, h, 0,\cdots,0)\in\sS(\mR^d;\mR^m),\ \ U_i(h):=U(f_i),
\end{align}
then it is easy to see that
$$
U=(U_1,\cdots, U_m).
$$
Moreover, for any $f,g\in\sS(\mR^d;\mR)$, we have
$$
\bE(U_i(f)U_j(g))=\int_{\mR^{d}}\hat f(\xi)\overline{\hat g(\xi)}\mu_{ij}(\dif\xi),
$$
which gives the covariance between different components.
In particular, for $f=(f_1,\cdots,f_m)\in\sS(\mR^d;\mR^m)$, we shall write
$$
U(f)=\<U,f\>=U_1(f_1)+\cdots+U_m(f_m).
$$

On the other hand, in the above definition, the Gaussian field $U$ is defined as a bounded linear operator from $\mH_0$ to $L^2(\Omega,\bP)$. A natural question arises: can we view $U$ as an $m$-dimensional random distribution of the spatial variable $x\in\mR^d$? We shall provide an affirmative answer and show its Besov regularity in the variable $x$ below.

To achieve this, we first introduce the convolution of $U$ with $\varphi\in\mH_0$. For a given $\varphi\in\mH_0$, we define
$$
(\varphi*U)(x):=U_\varphi(x):=U(\tau_x\varphi)\in L^2(\Omega,\bP).
$$
\bt
For any $\varphi\in\mH_0$, the function $x\mapsto U_\varphi(x)$ is continuous in $L^2(\Omega,\bP)$,
and for any  $k\in\mN$ and $\varphi\in\mH_{k}$,  $x\mapsto U_\varphi(x)$ is $k$-order Fr\'echet differentiable, and for each $x\in\mR^d$,
\begin{align}\label{AQ11}
\nabla^k_xU_\varphi(x)=\<U, \tau_x\nabla^k\varphi\>,\ \ \bP-a.s.
\end{align}
Moreover, for any $p\geq 2$, there is a constant $C=C(d,p)>0$ such that for all $\varphi\in\mH_{k}$,
\begin{align}\label{AG2}
\sup_{x\in\mR^d}\|\nabla^kU_\varphi(x)\|_{L^p(\Omega)}\leq C\|\varphi\|_{\mu;k}.
\end{align}
\et
\begin{proof}
Note that
$$
\bE|U_\varphi(x)-U_\varphi(y)|^2=\bE|U(\tau_x\varphi-\tau_y\varphi)|^2=\|\tau_x\varphi-\tau_y\varphi\|_{\mH_0}^2.
$$
The continuity of $x\mapsto U_\varphi(x)$ now follows by Lemma \ref{Le61}.
Similarly, by \eqref{1127:00}, we have \eqref{AQ11}.
Moreover, by  the hypercontractivity of Gaussian random variables, we have for any $p\ge2$,
\begin{align*}
\|\nabla^kU_\varphi(x)\|_{L^p(\Omega)}&\lesssim \|U_{\nabla^k\varphi}(x)\|_{L^2(\Omega)}=
\left(\sum_{i,j=1}^m\int_{\mR^d}\widehat{\nabla^k\varphi_i}(\xi)
\overline{\widehat{\nabla^k\varphi_j}(\xi)}\mu_{ij}(\dif \xi)\right)^{1/2}\\
&=\left(\sum_{i,j=1}^m\int_{\mR^d}|\xi|^{2k}\widehat{\varphi_i}(\xi)
\overline{\widehat{\varphi_j}(\xi)}\mu_{ij}(\dif \xi)\right)^{1/2}\leq \|\varphi\|_{\mu;k}.
\end{align*}
This completes the proof.
\end{proof}

Now for $f\in\sS(\mR^d;\mR)$, we define
$$
U_\varphi(f):=\int_{\mR^d}U_\varphi(x)f(x)\dif x.
$$
\bl
For any $\varphi\in \mH_0$, $U_\varphi$ is an $\mR$-valued Gaussian noise with spectral measure 
$$
\mu_\varphi(\dif\xi)=\sum_{i,j=1}^m\hat\varphi_i(\xi)\hat\varphi_j(-\xi)\mu_{ij}(\dif \xi).
$$
\el
\begin{proof}
For any $f,g\in\sS(\mR^d;\mR)$, by definition and Fubini's theorem, we have
\begin{align*}
\bE\big(U_\varphi(f)U_\varphi(g)\big)
&=\int_{\mR^d}\int_{\mR^d}\bE\big(U_\varphi(x)U_\varphi(y)\big)f(x)g(y)\dif x\dif y\\
&=\sum_{i,j=1}^m\int_{\mR^d}\int_{\mR^d}\left[\int_{\mR^d}\widehat{\tau_x\varphi_i}(\xi)\overline{\widehat{\tau_y\varphi_j}(\xi)}\mu_{ij}(\dif\xi)\right]f(x)g(y)\dif x\dif y\\
&=\sum_{i,j=1}^m\int_{\mR^d}\int_{\mR^d}\left[\int_{\mR^d}\e^{{\rm i}(\xi\cdot x-\xi\cdot y)}\widehat{\varphi_i}(\xi)\overline{\widehat{\varphi_j}(\xi)}
\mu_{ij}(\dif\xi)\right]f(x)g(y)\dif x\dif y\\
&=\sum_{i,j=1}^m\int_{\mR^d}\hat f(-\xi)\hat g(\xi) \widehat{\varphi_i}(\xi)\widehat{\varphi_j}(-\xi)
\mu_{ij}(\dif\xi)=\<f,g\>_{\mu_\varphi}.
\end{align*}
The result now follows.
\end{proof}
Below, for a Banach space $\mB$, let $\bB^{s}_{p}(\mB):=\bB^{s}_{p,\infty}(\mB)$ be the $\mB$-valued Besov space. 
For $\kappa\in\mR$, let $\rho_\kappa(x):=(1+|x|)^{-\kappa}$. Let $\bB^{s}_{p}(\rho_\kappa;\mB)$ be the weighted Besov space, which is defined by (see \cite{HZZZ21})
\begin{align}\label{Bes9}
\|f\|_{\bB^{s}_{p}(\rho_\kappa;\mB)}:=\|f\rho_k\|_{\bB^{s}_{p}(\mB)}.
\end{align}
We have the following regularity result.
\bt\label{Th65}
Under {\bf (A)}, for any $p\in[1,\infty)$, $\kappa>d/p$ and $\beta<\ell/2$, it holds that
	\begin{align}\label{1207:00}
		U\in \bB^{\ell/2}_{\infty}(L^p(\Omega;\mR^m))\subset\bB^{\ell/2}_{p}(\rho_\kappa; L^p(\Omega;\mR^m))\subset L^p(\Omega;\bB^{\beta}_{p}(\rho_\kappa; \mR^m)).
	\end{align}
\et
\begin{proof}
First of all, recalling the definition of $(\phi_j)_{j\geq -1}$ in \eqref{Phj},
one sees that for any $\gamma\in\mR$,
\begin{align}\label{SpRI}
	|\phi_j(\xi)|\lesssim 1\wedge\big(2^{\gamma j}(1+|\xi|)^{-\gamma}\big),\ \ j\geq -1,\ \ \xi\in\mR^{d}.
	\end{align}
Indeed, for $j=-1$, since supp$\phi_{-1}\subset \bB_1$, it is obvious. For $j\geq 0$, noting that
	$$
	K_j:=\text{supp}\phi_j\subset \{\xi: 2^{j-1}\leq|\xi|\leq{2^{j+1}}\},
	$$
we have for any $\gamma\in\mR$, 
	\begin{align*}
	|\phi_j(\xi)|\le\1_{K_j}(\xi)\lesssim2^{j \gamma}/(1+|\xi|)^{\gamma}.
	\end{align*}
For each $i=1,\cdots, m$, under \eqref{AG3}, by definition \eqref{EE9} and the hypercontractivity of Gaussian random variables, we have for any $p\in[1,\infty)$ and $j=0,1,\cdots$,
	\begin{align}
	\|\cR_jU_i(x)\|_{L^p(\Omega)}&\lesssim\|\cR_jU_i(x)\|_{L^2(\Omega)}=\|U_{i}(\tau_x\check\phi_j)\|_{L^2(\Omega)}
	=\left(\int_{\mR^{d}}|\phi_j(\xi)|^2\mu_{ii}(\dif \xi)\right)^{1/2}\no\\
	&\lesssim 2^{-\ell j/2}\left(\int_{\mR^{d}}(1+|\xi|)^{\ell}\mu_{ii}(\dif \xi)\right)^{1/2}
	\lesssim 2^{-\ell j/2},\label{AG4}
	\end{align}
and for $\kappa>d/p$,
$$
\|\rho_\kappa\cR_j U\|^p_{L^p(\mR^d;L^p(\Omega))}=\int_{\mR^d}(1+|x|)^{-p\kappa}\|\cR_jU(x)\|^p_{L^p(\Omega)}\dif x
\lesssim 2^{-p\ell j/2}\int_{\mR^d}(1+|x|)^{-p\kappa}\dif x.
$$
Hence, by \cite[Theorem 2.7]{HZZZ21},
$$
U\in \bB^{\ell/2}_{\infty}(L^p(\Omega;\mR^m))\subset\bB^{\ell/2}_{p}(\rho_\kappa; L^p(\Omega;\mR^m)).
$$
Finally, for $\beta<\ell/2$, by Fubini's theorem, we have
$$
\Big\|\sup_{j\geq 0}2^{\beta j}\|\rho_\kappa\cR_j U\|_p\Big\|_{L^p(\Omega)}\leq\sum_{j\geq 0}2^{\beta j}\|\rho_\kappa\cR_j U\|_{L^p(\mR^d\times \Omega)}
\leq \sum_{j\geq 0}2^{(\beta-\ell/2) j}\|U\|_{\bB^{\ell/2}_{p}(\rho_\kappa; L^p(\Omega;\mR^m))}.
$$
The proof is completed .
\end{proof}
\subsection{Examples: Divergence-free Gaussian field}
In this section we  present  concrete examples to illustrate our results.
The following example is standard.
\bx\rm\label{KE1}
Let $m=1$. For  $\gamma\in[0,d)$, let
$$
\mu(\dif\xi):=|\xi|^{-\gamma}\dif\xi.
$$
Obviously, \eqref{AG3} holds.
In this case, for $\gamma>0$, it is well known that
for some $c_{d,\gamma}>0$ (see \cite[p117, Lemma 2]{St}),
\begin{align*}
\hat\mu(x)=c_{d,\gamma}|x|^{\gamma-d}.
\end{align*}
In particular, for any $f,g\in\sS(\mR^{d})$,
\begin{align*}
\bE\big(U(f)U(g)\big)&=\int_{\mR^{d}}\hat f(\xi)\hat g(-\xi)\mu(\dif\xi)
=c_{d,\gamma}\int_{\mR^{d}}\int_{\mR^{d}}
	\frac{f(x)g(y)}{|x-y|^{d-\gamma}}\dif x\dif y.
\end{align*}
For $\gamma=0$, we have
$$
\hat\mu(x)=\delta_0(\dif x)
$$
and
\begin{align*}
\bE\big(U(f)U(g)\big)=\int_{\mR^d}f(x)g(x)\dif x.
\end{align*}
Namely, $X$ is a space white noise on $\mR^d$.
As in \eqref{AG4} and by the change of variables, we have
	\begin{align*}
	\|\cR_jU(x)\|_{L^p(\Omega)}&\lesssim
	\left(\int_{\mR^{d}}|\phi_j(\xi)|^2\mu(\dif \xi)\right)^{1/2}
	=2^{j(d-\gamma)/2}\left(\int_{\mR^{d}}|\phi_0(\xi)|^2\frac{\dif \xi}{|\xi|^\gamma}\right)^{1/2}.
	\end{align*}
	Hence, for any $p\in[1,\infty)$, $\kappa>d/p$ and $\beta<(\gamma-d)/2$,
\begin{align}\label{AG6}
	U\in \bB^{(\gamma-d)/2}_{\infty}(L^p(\Omega))\subset \bB^{(\gamma-d)/2}_{p}(\rho_\kappa; L^p(\Omega))
	\subset L^p(\Omega,\bB^{\beta}_{p}(\rho_\kappa)).
\end{align}
In particular, bigger $\gamma$ means that Gaussian field $U$ has better regularity.
\ex
\bx[Divergence free vector-valued Gaussian field]\rm\label{KE2}
For given $\gamma\in[0,d)$, define
\begin{align}\label{AG87}
\mu^{(\gamma)}(\dif\xi):=|\xi|^{-\gamma}\Big(\mI_{d\times d}-\frac{\xi\otimes\xi}{|\xi|^2}\Big)\dif\xi.
\end{align}
It is easy to see that \eqref{AG3} holds for $\ell<\gamma-d$, and for $a:=(a_i)_{i=1}^d\subset\mC$ and $\xi\in\mR^d{\backslash\{0\}}$,
$$
a^T\Big(\mI_{d\times d}-\frac{\xi\otimes\xi}{|\xi|^2}\Big)\bar a=|a|^2-\sum_{i,j}(a_i\xi_i\bar a_j\xi_j)/|\xi|^2
=|a|^2-\Big|\sum_{i}(a_i\xi_i)\Big|^2/|\xi|^2\geq 0.
$$
Thus {\bf (A)} holds.
Let $b^{(\gamma)}$ be the $\mR^d$-valued Gaussian random field with matrix-valued spectral measure $\mu^{(\gamma)}(\dif x)$.
Then by \eqref{1207:00}, for any $p\in[1,\infty)$, $\kappa>\frac dp$ and $\beta<\frac{\gamma-d}{2}$,
$$
b^{(\gamma)}\in L^p(\Omega;\bB^{\beta}_{p}(\rho_\kappa; \mR^d)).
$$
In particular,  for any $\kappa>0$ {arbitrarily small} and $\beta<\frac{\gamma-d}2$, one can choose $p$ large enough so that $\kappa>\frac dp$ and $\beta+\frac{d}{p}<\frac{\gamma-d}{2}$, and for 
$\bP$-almost all $\omega$,
\begin{align}\label{1210:00}
b^{(\gamma)}(\omega,\cdot)\in \bB^{\beta+d/p}_{p}(\rho_\kappa; \mR^d)\subset 
\bB^{\beta}_{\infty}(\rho_\kappa; \mR^d)=:\bC^\beta_\kappa(\mR^d),
\end{align}
and
$$
\div b^{(\gamma)}(\omega,\cdot)=0.
$$
In fact, for any $f\in\sS(\mR^d)$,
\begin{align}
\bE|\<\div b^{(\gamma)}, f\>|^2&=\bE|\< b^{(\gamma)},\nabla f\>|^2
=\int_{\mR^d}\widehat{\nabla f}(\xi)\Big(\mI_{d\times d}-\frac{\xi\otimes\xi}{|\xi|^2}\Big)\big[\overline{\widehat{\nabla f}(\xi)}\big]^*|\xi|^{-\gamma}\dif\xi\no\\
&=\int_{\mR^d}|\hat f(\xi)|^2\xi\Big(\mI_{d\times d}-\frac{\xi\otimes\xi}{|\xi|^2}\Big)\xi^*|\xi|^{-\gamma}\dif\xi=0,\label{Div8}
\end{align}
where $*$ stands for the transpose of a row vector.

{In particular, when $\gamma\in(d-2,d)$, we can take $\beta>-1$. In this case, for 
$\bP$-almost all $\omega$,} by Lemma \ref{App:00} in appendix, there exists 
a decomposition $b^{(\gamma)}=b^{(\gamma)}_1+b^{(\gamma)}_2$ with
\begin{align*}
    b^{(\gamma)}_1{(\omega,\cdot)}\in 
    \bC^{\beta-\kappa}_0(\mR^d)),\ \ \div b^{(\gamma)}_1{(\omega,\cdot)}\in \bC^{\beta-\kappa}_0(\mR)\quad \text{and}\quad 
    b^{(\gamma)}_2{(\omega,\cdot)}\in \bC^{\kappa}_{-\beta+2\kappa}(\mR^d),
\end{align*}
where $\kappa\in(0,(1+\beta)/2)$. In particular, $b^{(\gamma)}_2{(\omega,\cdot)}$ is continuous and
\begin{align*}
    |b^{(\gamma)}_2({\omega, } x)|\lesssim (1+|x|)^{-\beta+2\kappa}\lesssim 1+|x|.
\end{align*}
Therefore, one can apply Theorem \ref{Th1} to solve the following SDE 
$$
\dif X_t=b^{(\gamma)}(\omega, X_t)\dif t+\sqrt2\dif W_t,
$$
where $X_t$ can be thought of the diffusion in random environment (see \cite{FA98}).
\ex

\br
Consider $d=3$. Let ${\bf U}^{(\gamma)}:=(U_1,U_2,U_3)$ be the 3-dimensional Gaussian field with spectral measure matrix $\mu(\dif \xi)=|\xi|^{-\gamma}\mI_{3\times 3}\dif\xi$, 
where $\gamma\in(1,3)$.  
Let $b^{(\gamma)}$ be the $\mR^3$-valued Gaussian field in the above example.
It is easy to see that (see \cite{FA98})
\begin{align*}
b^{(\gamma)}=\nabla\times (-\Delta)^{-\frac12}{\bf U}^{(\gamma)},\quad \nabla\times:=\begin{pmatrix}
0&-\p_3&\p_2\\
\p_3&0&-\p_1\\
-\p_2&\p_1&0
\end{pmatrix}.
\end{align*}
\er

\br
Consider $d=2$. Let $U$ be the two-dimensional spatial white noise with spectral measure $\mu(\dif \xi)=\dif\xi$ (see Example \ref{KE1}).  
Let $\sS_0(\mR^2)$ be the class of Schwartz function $f\in \sS(\mR^2)$ with $0\notin{\rm supp}(\hat f)$.
Let $b^{(0)}$ be the $\mR^2$-valued Gaussian field in Example \ref{KE2}. Then
it is easy to see that  for  any $f,g\in \sS_0(\mR^2)$,
\begin{align*}
\bE\big( (\nabla^{\perp}(-\Delta)^{-\frac12}U)(f)(\nabla^{\perp}(-\Delta)^{-\frac12})U(g)\big)
=\bE(b^{(0)}(f)b^{(0)}(g)),
\end{align*}
where $\nabla^\perp:=(\p_{x_2},-\p_{x_1})$.
In other words, 
\begin{align*}
b^{(0)}=\nabla^{\perp}(-\Delta)^{-\frac12}U \mbox{on $\sS_0(\mR^2)$}.
\end{align*}
However, $\sS_0(\mR^2)$ is not a  determining class since for $f\in\sS_0(\mR^2)$, $\int_{\mR^2}f(x)\dif x=0$.
So we can not claim $b^{(0)}=\nabla^{\perp}(-\Delta)^{-\frac12}U$ as distributions.
\er

\subsection{Stochastic heat equations}\label{sec63}
Let $\gamma\in(-\infty,d)$ and $\eta^{(\gamma)}$ be an $\mR^d$-valued time-white Gaussian noise with matrix-valued spectral measure $\mu^{(\gamma)}(\dif \xi)\dif t$,
where $\mu^{(\gamma)}$ is given by \eqref{AG87}.
More precisely, $\eta^{(\gamma)}$ is a Gaussian field over $\mR_+\times\mR^2$ with covariance given by that 
for any $f,g\in L^2([0,\infty);\sS(\mR^d;\mR^d))$,
\begin{align*}
\bE \big(\eta^{(\gamma)}(f)\eta^{(\gamma)}(g)\big)=\sum_{i,j=1}^d\int_0^\infty\!\!\!\int_{\mR^2}\hat{f}_i(t,\xi)\hat{g}_j(t,-\xi)\mu^{(\gamma)}_{ij}(\dif \xi)\dif t.
\end{align*}
Now, we consider the following $\mR^d$-valued stochastic heat equation (abbreviated as SHE) on $\mR^d$:
\begin{align}\label{SPDE0}
\p_t u^{(\gamma)}=\Delta u^{(\gamma)}+\eta^{(\gamma)},\quad u^{(\gamma)}_0=0.
\end{align}
Let
$$
p_t(x)=(4\pi t)^{-1}\exp\{-|x|^2/(4t)\}.
$$
By  Duhamel's formula, the solution is determined by the following formula:
\begin{align*}
u^{(\gamma)}(t,x):=\int_0^t\int_{\mR^d}p_{t-s}(x-y)\eta^{(\gamma)}(\dif y,\dif s):=\eta^{(\gamma)}(\1_{(0,t)}\tau_x p_{t-\cdot}).
\end{align*}
Noting that for any $j\geq 0$ and $i=1,\cdots,m$,
$$
\cR_ju^{(\gamma)}_i(t,x)=\eta^{(\gamma)}_i(\1_{(0,t)}\tau_x\cR_jp_{t-\cdot}),
$$
as in Theorem \ref{Th65}, by  \eqref{Phj} and the change of variables, we have 
\begin{align*}
\|\cR_ju^{(\gamma)}_i(t,x)\|^2_{L^p(\Omega)}&\lesssim \|\cR_ju^{(\gamma)}_i(t,x)\|^2_{L^2(\Omega)}
=\int_0^t\!\!\!\int_{\mR^d}|\tau_x\cR_jp_{t-s}(\xi)|^2\mu^{(\gamma)}_{ii}(\dif \xi)\dif s\\
&=\int_0^t\!\!\!\int_{\mR^d}|\phi_j(\xi)|^2\e^{-2(t-s)|\xi|^2}|\xi|^{-\gamma}(1-|\xi_i|^2/|\xi|^2)\dif \xi\dif s\\
&\leq 2^{j(d-\gamma)}\int_0^t\!\!\!\int_{\mR^d}|\phi_0(\xi)|^2\e^{-2^{2j+1}s|\xi|^2}|\xi|^{-\gamma}\dif \xi\dif s\\
&\lesssim 2^{j(d-\gamma-2)}\int_{\mR^d}|\phi_0(\xi)|^2|\xi|^{-2-\gamma}\dif \xi,
\end{align*}
and also,
$$
\|\cR_0u^{(\gamma)}_{-1}(t,x)\|^2_{L^p(\Omega)}\lesssim \int_{\mR^d}|\phi_{-1}(\xi)|^2|\xi|^{-\gamma}\dif \xi<\infty.
$$
Hence, for $\gamma<d-2$,
$$
u^{(\gamma)}\in L^\infty(\mR_+; \bB^{1+\frac{\gamma-d}{2}}_\infty(L^p(\Omega;\mR^d))).
$$
Moreover, as \eqref{Div8}, we also have
$$
\div u^{(\gamma)}=0.
$$
In particular, we can use Theorem \ref{Th1} and Lemma \ref{App:00}, as  in Example \ref{KE2}, to obtain the weak well-posedness of SDE with {\it time-dependent} vector field $b(t,x)=u^{(\gamma)}(t,x)$ with  $\gamma\in(d-4,d-2)$.

\begin{appendix}
\renewcommand{\thetable}{A\arabic{table}}
\numberwithin{equation}{section}

\section{Auxiliary lemmas}

In the appendix, we recall a useful decomposition result for a distribution $f\in\bC^\alpha_\kappa$ that is taken from \cite[Lemma 2.4]{GH19}.
Let $\kappa,\alpha\in\mR$ and $m\in\mN$. Recalling $\rho_\kappa(x)=(1+|x|)^{-\kappa}$, we introduce the weighted H\"older space by (see \eqref{Bes9})
$$
    \bC^\alpha_\kappa(\mR^m):=
    \left\{
    \begin{aligned}
    &\bB^\alpha_\infty(\rho_\kappa;\mR^m),\ \alpha\notin\mN_0,\\
    &\bH^\alpha_\infty(\rho_\kappa;\mR^m),\ \alpha\in\mN_0.
    \end{aligned}
    \right.
$$
For any $f\in \bC^\alpha_\kappa(\mR^m)$, we define
\begin{align*}
    \sU_\ge(f):=\sum_{k=-1}^\infty\phi_k (\mI-S_k)f,\quad    \sU_\le(f):=\sum_{k=-1}^\infty\phi_k S_kf,
\end{align*}
where $\phi_k$ is defined in \eqref{Phj} and $S_k$ is defined in \eqref{EM9}. 
We shall use the following easy fact:
\begin{align}\label{Fac}
\phi_j=\phi_j\tilde\phi_j\mbox{with }\tilde{\phi}_j:=\phi_{j-1}+\phi_j+\phi_{j+1}.
\end{align}
We first introduce the following result.
\bl
For any $\kappa,\alpha\in\mR$, it holds that
\begin{align}\label{1213:equ}
    \|f\|_{\bC^\alpha_\kappa}\asymp\sup_{j\geq -1}\|\phi_jf\|_{\bC^\alpha_\kappa}.
\end{align}
\el
\begin{proof}
{First of all, by the definition and scaling of function $\phi_j$, for any fixed $m\in\mN$, 
$$
\sup_{j\ge-1}\|\phi_j\|_{\bC^m_0}<\infty,
$$
which implies that for $m\ge |\alpha|$,
  $$
  \sup_{j\geq -1}\|\phi_jf\|_{\bC^\alpha_\kappa}{\lesssim}\sup_{j\geq -1}\|\phi_j\|_{\bC^m_0}\|f\|_{\bC^\alpha_\kappa}\lesssim \|f\|_{\bC^\alpha_\kappa}.
  $$
  Thus it remains to show $\|f\|_{\bC^\alpha_\kappa}\lesssim\sup_{j\geq -1}\|\phi_jf\|_{\bC^\alpha_\kappa}.$ 
}
 
  By interpolation, without loss of generality, we may assume $\alpha\in\mZ$. 
  Let $\alpha\in\mN_0$. 
  Thanks to $\mathrm{supp}\,\phi_j\subset B_{2^{j+2}/3}\setminus B_{2^{j-1}}$,
  $$
  \|f\|_{\bC^\alpha_\kappa}=\Big\|\sum_{j\geq -1}\phi_j f\Big\|_{\bC^\alpha_\kappa}\lesssim\sup_{j\geq -1}\|\phi_j f\|_{\bC^\alpha_\kappa}.
  $$
  Thus we get \eqref{1213:equ} for $\alpha\in\mN_0$. On the other hand, by duality and \eqref{Fac}, we have
      \begin{align*}
          \|f\|_{\bC^{-\alpha}_\kappa}&=\|\rho_\kappa f\|_{\bC^{-\alpha}}
          =\sup_{\|g\|_{\bH^{\alpha}_1}\leq 1}|\<\rho_\kappa f,g\>|=\sup_{\|g\|_{\bH^{\alpha}_1}\leq 1}\Big|\sum_{j}\<\rho_\kappa\tilde\phi_j f, {\phi}_jg\>\Big|\\
  &        \leq\sup_{\|g\|_{\bH^{\alpha}_1}\leq 1}\sum_{j}\Big(\|\phi_jg\|_{\bH^{\alpha}_1}\|\rho_\kappa\tilde\phi_jf\|_{\bC^{-\alpha}}\Big)
  \lesssim\sup_{\|g\|_{\bH^{\alpha}_1}\leq 1}\sum_{j}\|\phi_jg\|_{\bH^{\alpha}_1}\sup_j\|\phi_jf\|_{\bC^{-\alpha}_\kappa}.
      \end{align*}
      We note that 
      \begin{align*}
          \sum_{j}\|\phi_jg\|_{\bH^{\alpha}_1}&\lesssim  \sum_{j}\sum_{k\le\alpha}\|\nabla^{\alpha-k}\phi_j\nabla^kg\|_1\lesssim \sum_{j}\sum_{k\le\alpha}\|\nabla^{\alpha-k}\phi_j\tilde{\phi}_j\nabla^kg\|_1\\
          &{\lesssim} \sum_{j}\sum_{k\le\alpha}2^{(k-\alpha)j}\|\tilde{\phi}_j\nabla^kg\|_1\lesssim  \sum_{k\le\alpha}\sum_{j}\|{\phi}_j\nabla^kg\|_1\lesssim \sum_{k\le\alpha}\|\nabla^kg\|_1
          \lesssim\|g\|_{\bH^{\alpha}_1}.
      \end{align*}
Then combining the above two estimates, we get $\|f\|_{\bC^{-\alpha}_\kappa}\lesssim\sup_j\|\phi_jf\|_{\bC^{-\alpha}_\kappa}$.
\end{proof}

Then we have the following result.
\bl\label{App:00}
Let $\alpha,\kappa\in\mR$ and $f\in \bC^\alpha_\kappa(\mR^d)$. 
Then $f=\sU_\ge(f)+\sU_\le(f)$ and for any $\beta,\gamma>0$,
\begin{align}\label{1213:ge}
    \|\sU_\ge(f)\|_{\bC^{\alpha-\beta}_{\kappa-\beta}}\lesssim\|f\|_{\bC^\alpha_\kappa},\quad   \|\div\sU_\ge(f)\|_{\bC^{\alpha-\beta}_{\kappa-\beta}}\lesssim\|f\|_{\bC^\alpha_{\kappa+1}}+\|\div f\|_{\bC^\alpha_{\kappa}},
\end{align}
and
\begin{align}\label{1213:le}
    \|\sU_\le(f)\|_{\bC^{\alpha+\gamma}_{\kappa+\gamma}}\lesssim\|f\|_{\bC^\alpha_\kappa}.
\end{align}
\el
\begin{proof}
By \eqref{1213:equ}, we have
    \begin{align*}
        \|\sU_\ge(f)\|_{\bC^{\alpha-\beta}_{\kappa-\beta}}\lesssim \sup_{j}\|\phi_j\sU_\ge(f)\|_{\bC^{\alpha-\beta}_{\kappa-\beta}}\le 
        \sup_{j}\left\| \sum_{k=-1}^\infty{\phi}_j\phi_k(\mI-S_k)f\right\|_{\bC^{\alpha-\beta}_{\kappa-\beta}}\le \sup_{k}\left\| \phi_k(\mI-S_k)f\right\|_{\bC^{\alpha-\beta}_{\kappa-\beta}},
    \end{align*}
    where for $m\ge |\alpha|$
   \begin{align*}
       \left\| \phi_k(\mI-S_k)f\right\|_{\bC^{\alpha-\beta}_{\kappa-\beta}}\lesssim \|\phi_k\|_{\bC^m_{-\beta}}\|(\mI-S_k)f\|_{\bC^{\alpha-\beta}_{\kappa}}{\lesssim }2^{\beta k}\sum_{i=k}^\infty\|\cR_if\|_{\bC^{\alpha-\beta}_{\kappa}}.
    \end{align*} 
It follows from \cite[(2.22) and (2.14)]{HZZZ21} that
\begin{align*}
    \|\cR_if\|_{\bC^{\alpha-\beta}_{\kappa}}&\lesssim\sup_j2^{(\alpha-\beta)j}\|\cR_j\cR_if\|_{L^\infty(\rho_\kappa)}\lesssim \sup_{j\sim i}2^{(\alpha-\beta)j}\|\cR_j\cR_if\|_{L^\infty(\rho_\kappa)}\\
    &\lesssim 2^{(\alpha-\beta)i}\|\cR_if\|_{L^\infty(\rho_\kappa)}\lesssim2^{-\beta i}\|f\|_{\bC^\alpha_\kappa}, 
\end{align*}
    which implies that
    \begin{align*}
       \|\sU_\ge(f)\|_{\bC^{\alpha-\beta}_{\kappa-\beta}}&\lesssim  \sup_{k}2^{\beta k}\sum_{i=k}^\infty\|\cR_if\|_{\bC^{\alpha-\beta}_{\kappa}}
       \lesssim \sup_{k}2^{\beta k}\sum_{i=k}^\infty2^{-\beta i}\|f\|_{\bC^\alpha_\kappa}\lesssim \|f\|_{\bC^\alpha_\kappa}.
    \end{align*}
    Similarly, we have
   \begin{align*}
        \|\div\sU_\ge(f)\|_{\bC^{\alpha-\beta}_{\kappa-\beta}}&\le \left\| \sum_{k=-1}^\infty\nabla\phi_k\cdot(\mI-S_k)f\right\|_{\bC^{\alpha-\beta}_{\kappa-\beta}}+\|\sU_\ge(\div f)\|_{\bC^{\alpha-\beta}_{\kappa-\beta}},
  \end{align*}
where 
\begin{align*}
    \|\sU_\ge(\div f)\|_{\bC^{\alpha-\beta}_{\kappa-\beta}}\lesssim \| \div f\|_{\bC^{\alpha}_{\kappa}}
\end{align*}
  and
 \begin{align*}
       \left\| \sum_{k=-1}^\infty\nabla\phi_k\cdot(\mI-S_k)f\right\|_{\bC^{\alpha-\beta}_{\kappa-\beta}} &\lesssim \sup_{j}\left\| \sum_{k=-1}^\infty{\phi}_j\nabla\phi_k\cdot(\mI-S_k)f\right\|_{\bC^{\alpha-\beta}_{\kappa-\beta}}\\
       &\le \sup_{k} \|\nabla\phi_k\|_{\bC^m_{-\beta-1}}\|(\mI-S_k)f\|_{\bC^{\alpha-\beta}_{\kappa+1}}\\
       &{\lesssim }\sup_{k} 2^{\beta k}\|(\mI-S_k)f\|_{\bC^{\alpha-\beta}_{\kappa+1}}\lesssim \|f\|_{\bC^{\alpha-\beta}_{\kappa+1}}.
    \end{align*} 
  The estimate for  \eqref{1213:le} is the same.
\end{proof}

\end{appendix}
\vspace{5mm}

{\bf Acknowledgement: 
\rm
We would like to thank Rongchan Zhu and Xiangchan Zhu for their invaluable contributions, particularly for their insights into the background of the critical cases and the uniqueness of generalized martingale solutions.  
We are also deeply thankful to the referees for their meticulous review of the manuscript and their constructive suggestions.
}

\end{document}